\documentclass[letterpaper,11pt,oneside,reqno]{article}

\usepackage[margin=1.0in]{geometry}
\usepackage{amsmath,amsthm,amsfonts,amssymb}
\usepackage[breaklinks,colorlinks,linktoc=page]{hyperref}
\hypersetup{linkcolor=blue,citecolor=blue,urlcolor=blue}
\usepackage{environ,xspace,array,longtable,booktabs}

\usepackage[msc-links,non-compressed-cites]{amsrefs}

\makeatletter
\renewcommand*{\toclevel@part}{1} %

\def\tjm@abbrev#1{\new@ifnextchar.{\tjm@abbrev@end{#1}}{#1\@\xspace}}
\def\tjm@abbrev@end#1.{#1}
\newcommand{\as}{\ifmmode\text{a.s.}\else\tjm@abbrev{a.s.}\fi}
\newcommand{\aev}{\ifmmode\text{a.e.}\else\tjm@abbrev{a.e.}\fi}
\newcommand{\ie}{\ifmmode\text{i.e.}\else\tjm@abbrev{i.e.}\fi}
\newcommand{\eg}{\ifmmode\text{e.g.}\else\tjm@abbrev{e.g.}\fi}
\newcommand{\iid}{\ifmmode\text{i.i.d.}\else\tjm@abbrev{i.i.d.}\fi}
\newcommand{\inid}{\ifmmode\text{i.n.i.d.}\else\tjm@abbrev{i.n.i.d.}\fi}
\newcommand{\pmp}{\ifmmode\text{p.m.p.}\else\tjm@abbrev{p.m.p.}\fi}
\newcommand{\wolog}{\ifmmode\text{w.l.o.g.}\else\tjm@abbrev{w.l.o.g.}\fi}
 
\def\MR#1{\relax\ifhmode\unskip\spacefactor\@m\space\fi \MRhref{MR#1}{#1}}
\renewcommand{\eprint}[1]{\href{https://arxiv.org/abs/#1}{arXiv:#1}}
 
\numberwithin{equation}{section}

\renewcommand{\newtheorem}{\@ifstar{\@tjm@newthm@star}{\@tjm@newthm}}
\def\@tjm@newthm@star{\@ifnextchar[{\@tjm@newthm@star@}{\@tjm@newthm@star@[]}}
\def\@tjm@newthm{\@ifnextchar[{\@tjm@newthm@}{\@tjm@newthm@[]}}
\def\@tjm@newthm@star@[#1]{\@xnthm *}
\def\@tjm@newthm@[#1]{\@xnthm \relax}
 
\newenvironment{bareproof}{\par\pushQED{\qed}}{\popQED\par}
 
\newcommand{\monthname}{%
  \ifcase\month
  \or January%
  \or February%
  \or March%
  \or April%
  \or May%
  \or June%
  \or July%
  \or August%
  \or September%
  \or October%
  \or November%
  \or December%
  \fi
}

\newtheoremstyle{hyp}{}{}{\upshape}{\parindent}{\scshape}{}{0.5em}%
  {(\thmnumber{#2})\thmnote{\normalfont{} (#3)}}
\makeatother
 
\usepackage[capitalize,noabbrev]{cleveref}
\AtBeginDocument{\RemoveFromHook{label}[firstaid/cleveref]}

\usepackage{graphicx}
\usepackage{mathtools}

\usepackage{xspace}
\usepackage{xfrac}
\usepackage{xcolor}
\usepackage{dsfont}
\usepackage{tikz}
\usetikzlibrary{arrows.meta, positioning, calc, hobby}
\usepackage{esint}
\usepackage{mathrsfs}
\usepackage{enumitem}

\usepackage{booktabs}

\numberwithin{equation}{section}

\makeatletter
\long\def\@part[#1]#2{%
	\ifdefined\NR@gettitle \NR@gettitle{#1}\fi
	\ifnum \c@secnumdepth >\m@ne
	\refstepcounter{part}%
	\addcontentsline{toc}{part}{\thepart\hspace{1em}#1}%
	\else
	\addcontentsline{toc}{part}{#1}%
	\fi
	{\centering \interlinepenalty \@M \normalfont \Large \bfseries
		\ifnum \c@secnumdepth >\m@ne \partname~\thepart.\enspace \fi
		#2\par}%
	\markboth{}{}%
	\nobreak \vskip 3ex
	\@afterheading}
\long\def\@spart#1{%
	\ifdefined\NR@gettitle \NR@gettitle{#1}\fi
	{\centering \interlinepenalty \@M \normalfont \Large \bfseries #1\par}%
	\nobreak \vskip 3ex
	\@afterheading}

\newcommand{\subjclass}[2][2020]{%
  \gdef\tjm@subjclassname{\textup{#1} Mathematics Subject Classification}%
  \gdef\tjm@subjclass{#2}}
\newcommand{\keywords}[1]{\gdef\tjm@keywords{#1}}
\newcommand{\keywordsname}{Key words and phrases}
\let\tjm@subjclass\@empty
\let\tjm@keywords\@empty
\def\tjm@mscfootnotes{%
  \begingroup
    \let\@thefnmark\@empty
    \ifx\tjm@subjclass\@empty\else
      \@footnotetext{{\itshape\tjm@subjclassname.}\enspace
        \tjm@subjclass\@addpunct.}%
    \fi
    \ifx\tjm@keywords\@empty\else
      \@footnotetext{{\itshape\keywordsname.}\enspace
        \tjm@keywords\@addpunct.}%
    \fi
  \endgroup}
\g@addto@macro\@maketitle{\g@addto@macro\@thanks{\tjm@mscfootnotes}}
\makeatother

\newtheoremstyle{definition}{}{}{\normalfont}{}{\itshape}{.}%
{5pt plus 1pt minus 1pt}%
{\thmname{#1}\thmnumber{ #2}\thmnote{\normalfont{} (#3)}}
\newtheoremstyle{hyp}{}{}{\upshape}{\parindent}{\bfseries}{}{0.5em}%
{(\thmnumber{#2})\thmnote{\normalfont{} (#3)}}

\newtheorem[{}\it]{theorem}{Theorem}[section]
\newtheorem{corollary}[theorem]{Corollary}
\newtheorem{lemma}[theorem]{Lemma}
\newtheorem{proposition}[theorem]{Proposition}

\newtheorem{conjecture}[theorem]{Conjecture}

\theoremstyle{definition}
\newtheorem{definition}{Definition}[section]
\newtheorem{remark}{Remark}[section]

\newtheorem{step}{Step}

\theoremstyle{hyp}
\newtheorem{assumption}{Assumption}

\crefname{axiom}{Axiom}{Axioms}
\crefname{assumption}{}{}
\crefformat{assumption}{(#2#1#3)}
\crefname{conjecture}{Conjecture}{Conjectures}

\Crefname{axiom}{Axiom}{Axioms}
\Crefname{assumption}{}{}
\Crefformat{assumption}{(#2#1#3)}
\Crefname{conjecture}{Conjecture}{Conjectures}

\Crefname{lemma}{Lemma}{Lemmas}
\Crefname{proposition}{Proposition}{Propositions}
\Crefname{corollary}{Corollary}{Corollaries}
\crefname{lemma}{Lemma}{Lemmas}
\crefname{proposition}{Proposition}{Propositions}
\crefname{corollary}{Corollary}{Corollaries}
\crefname{appendix}{Appendix}{Appendices}

\def\N{{\mathbb{N}}}

\def\E{{\mathbb{E}}}
\def\P{{\mathcal{P}}}
\def\R{{\mathbb{R}}}

\def\dd{{\,d}}

\newcommand{\GG}[1]{}

\DeclareMathOperator*{\essinf}{ess\,inf}
\DeclareMathOperator*{\esssup}{ess\,sup}

\def\rest{\hskip 1pt{\hbox to 10.8pt{\hfill
\vrule height 7pt width 0.4pt depth 0pt\hbox{\vrule height 0.4pt width 7.6pt depth 0pt}\hfill}}}

\allowdisplaybreaks

\definecolor{ibm0}{RGB}{255,176,0}
\definecolor{ibm1}{RGB}{254,97,0}
\definecolor{ibm2}{RGB}{220,38,127}
\definecolor{ibm3}{RGB}{120,94,240}
\definecolor{ibm4}{RGB}{100,143,255}

\definecolor{ito1}{RGB}{35,145,106}
\definecolor{ito2}{RGB}{0,114,178}
\definecolor{ito3}{RGB}{86,180,233}
\definecolor{ito4}{RGB}{240,228,66}
\definecolor{ito5}{RGB}{230,159,0}
\definecolor{ito6}{RGB}{213,94,0}
\definecolor{ito7}{RGB}{156,51,106}

\title{Asymptotics of Gibbs Posterior Measures on Infinite--Dimensional Hypothesis Spaces}
\author{Terje Mildner\thanks{Department of Statistics, University of Oxford, Oxford, United Kingdom. \texttt{\href{mailto:Terje.Mildner@stats.ox.ac.uk}{Terje.Mildner@stats.ox.ac.uk}}}}
\date{\the\day\ \monthname\ \the\year}

\begin{document}

    \setcounter{footnote}{1}
    \maketitle

\begin{abstract}
We prove almost sure contraction rates and consistency for Gibbs posterior measures on arbitrary hypothesis spaces, in general infinite dimensional, under dependent data, nonadditive and time--dependent losses and learning rates, misspecified models, and unattainable minimisers of the limiting objective. In particular, we develop a novel proof approach in contraction rate theory by studying the asymptotics of Gibbs measures through Lorentz spaces, real interpolation theory, and Laplace transforms dispensing with testing or variational formulations which are unavailable in the considered generality. As a consequence we obtain almost sure contraction rates for dominated Bayes and fractional Bayes posterior measures in the same generality, extending existing theory on misspecified Bayesian inference. Furthermore, we demonstrate our theoretical results on novel applications to Gibbs posteriors for Dirichlet process priors with a nonadditive Wasserstein metric loss function under ergodic, non--Markovian observables.
\end{abstract}

\tableofcontents

\section{Introduction}

The Gibbs posterior measure over any measurable space $(\Theta,\mathcal T)$ is a probability measures formed by combining a prior measure $\Pi_0$ with a non--negative function $R_t$, connecting observables $X_1^t$ and hypotheses $P$, typically through an exponentiated negative loss function and learning rate as $R_t(X_1^t,P)=\exp(-t\beta_tL_t(X_1^t,P))$, replacing the likelihood of Bayes posteriors,
\begin{equation}
    \Pi_t(A|X_1^t)=\int_A R_t(X_1^t,P)\,\Pi_0(dP)\Big/Z_t(\Theta),
\end{equation}
where $Z_t(\Theta)$ is the normalising constant.
These posterior measures, also referred to as generalised Bayes posteriors, have seen explosive interest in the last decade following the seminal work of Bissiri, Holmes, and Walker \cite{bissiri2016}, who provide these as a principled posterior inference framework, mitigating misspecified models and hypothesis spaces even when Bayes fails (cf. \cites{jewson2026}). 

Precisely because the choices of loss functions and learning rates allow for more control of the objective are we required to establish that these posteriors behave desirably.
However, \cites{bissiri2016} formulate their approach in the parametric setting based on a coherence axiom, whereas the literature has moved beyond this classical and parametric approach to consider non--additive loss functions, dependent observables, and infinite dimensional hypothesis spaces
\cites{altamirano2024,mildner2025a,laplante2025,matsubara2022,matsubara2024,lee2026,dellaporta2022,dellaporta2026,baraud2024,winter2026}. Their theoretical study has, however, largely been restricted to rather limited settings with additive loss functions and require the existence of minimisers of the loss, or been confined to particular instantiations of the Gibbs posterior \cites{grunwald2020,miller2021,syring2023,mcgoff2022,matsubara2022}. 

Asymptotic consistency of Gibbs posteriors is precisely what licenses that these concentrate on hypotheses that are of interest and close to optimal so that they may be used to mitigate model misspecification, and hence we cannot perform Gibbs posterior inference without understanding their asymptotics. Moreover, establishing contraction rates provides that consistency does not merely take place arbitrarily slowly but rather at a useful speed.

In this paper we will rigorously extend the applicability of Gibbs posteriors by studying the asymptotics of Gibbs posterior measures when the observables are dependent, the hypothesis space may be arbitrary and infinite--dimensional, the model misspecified, the loss function nonadditive, and no minimiser of the limiting objective may exist. To achieve this, we establish novel proof techniques in posterior asymptotics through Lorentz function spaces and real interpolation and rearrangement theory, as existing approaches to posterior asymptotics, including testing \cite{ghosal2017} and variational formulations \cites{alquier2016,katsevich2024,mildner2025b}, are generally unavailable.

\subsection{The Gibbs posterior measure}
Let $X_1,X_2, \dots$ be a sequence of observables taking values in the measurable space $(\Xi,\mathcal{X})$ generated from the probability space $(\Omega,\mathcal{F},P_0)$ where $P_0$ is the infinite process measure of the sequence $(X_t)_{t\in\N}$. Denote a finite subsequence $X_1,\dots,X_t\equiv X_1^t$ and let $\sigma(X_1^t)$ be its natural filtration.

Suppose $(\Theta,\mathcal{T})$ is the hypothesis space, in general misspecified \ie $P_0\notin\Theta$, with $\Pi_0$ a prior probability measure on this space. Let $L_t:\Xi^{(t)}\times\Theta\to[-\infty,+\infty]$ be a sequence of extended real valued functions connecting observables and hypotheses. We call these loss functions although these need not be bounded below, and we naturally assume this definition extends to include $P_0\notin\Theta$. Further, let $\beta_t>0$ be a deterministic sequence, and denote $s_t:=t\beta_t$. We assume throughout that for all $\delta>0$, $\sum_te^{-\delta s_t}<\infty$, a sufficient condition is $s_t/\log t\to+\infty$.

The Gibbs posterior measure is thus defined as
\begin{equation}\label{eqn:gibbs-measure-L_t}
\Pi_t(A|X_1^t)=\int_A\exp\bigl(-t\beta_tL_t(X_1^t;P)\bigr)\Pi_0(dP)\Big/ Z_t(\Theta)
\end{equation}
for all $A\in\mathcal T$, where 
\begin{equation}\label{eqn:Z_t(A)}
    Z_t(A):=\int_A e^{-s_tL_t(X_1^t;P)}\,\Pi_0(dP), \qquad s_t:= t\beta_t.
\end{equation}
We prefer the term Gibbs posterior \cites{catoni2007,jiang2008} to ``generalised Bayes'' to illustrate that we do not necessitate the coherence of Bissiri, Holmes, and Walker \cite{bissiri2016} but rather extend their applicability and may nevertheless perform principled posterior inference. Moreover, this further highlights the intricate relationship to the use of ``Gibbs measures'' in PDE theory and mathematical physics (see for instance Barashkov and Gubinelli \cite{barashkov2020}). %

Equivalently, \eqref{eqn:gibbs-measure-L_t} may be written using $\Delta_t(X_1^t;P)=L_t(X_1^t;P)-L_t(X_1^t;P_0)$, where we assume that $L_t(X_1^t;P_0)\in\R$ \as, and $R_t(P):=\exp\{-s_t\Delta_t(X_1^t;P)\}$ as
\begin{equation}\label{eqn:gibbs-measure-R_t}
    \Pi_t(A|X_1^t):=\cfrac{\int_AR_t(P)\,\Pi_0(dP)}{\int_\Theta R_t(P)\,\Pi_0(dP)}\,.
\end{equation}
We assume throughout for $R_t:\Xi^{(t)}\times\Theta\to[0,\infty]$ to be an extended real valued function that is $\sigma(X_1^t)\otimes\mathcal{T}$ measurable for all $t$ and for which the denominator is \as positive so that the object of study is well defined. 

While we work with $\Delta_t(X_1^t;P)$ for familiarity of the objective in Gibbs posterior inference, we stress that this is not necessary for our results. We may likewise work with $R_t$ directly as a function mapping to $[0,+\infty]$, and further have this be a random sequence of functions without observables present, $R_t(\omega,P)$, or fully deterministic $R_t(P)$. 

\begin{remark}
We have not assumed any topology on any of the measurable spaces, and these are in general assumed to be infinite dimensional. The observables can be arbitrarily dependent, and we do not in general suppose that the model is well specified in that $P_0\in\Theta$, that is we are in the $M$--open case of \cite{bernardo2000}. 
The misspecified setting itself motivated the development of the Gibbs posteriors, and hence any such assumption would be counterproductive.
Neither do we assume that there exists some common dominating measure, nor that the hypotheses are probability measures, they likewise could be (generalised) functions in some function space. The notation $P\in\Theta$ has been adapted to relate to the canonical setting of the Generalised Bayes literature of utilising probability distributions in the loss, albeit in the parametric setting, and to highlight the nonparametric generality.
\end{remark} 
\begin{remark}
Our definition of $\Delta_t$ utilising $P_0$ is motivated by not requiring the existence of a minimiser of the limiting objective, which we make precise below, and hence we may be unable to work with a minimiser as in the existing literature on asymptotics of Gibbs measures \cites{miller2021,syring2023} and misspecified Bayes \cites{kleijn2006,walker2013}. 
\end{remark}

\subsection{Contraction rates for Gibbs posterior measures}

Let $\epsilon(t)$ be a positive rate function converging to zero, and $P\mapsto\zeta(P)$ a function that quantifies the optimality of the hypotheses $P\in\Theta$ in some sense which we make precise later. Consider the decreasing sequence of sets 
\begin{equation}
    N_t:=\{P\in\Theta:\zeta(P)\lesssim \epsilon(t)\}.
\end{equation}
Contraction rates ask for $\epsilon(t)$ decaying sufficiently slowly for the corresponding $N_t$ to retain sufficient prior mass, while decaying sufficiently quickly for the rate to remain meaningful. In this paper we will primarily investigate sufficient, and in some cases necessary, conditions under which the following conjecture holds true.

\begin{conjecture}
For $\epsilon(t)\to0$ and $t\epsilon(t)\to+\infty$, it holds, under certain conditions on the prior hypothesis space, data generating process, and functions $L_t$, that
\begin{equation}
    \Pi_t(\{P\in\Theta:\zeta(P)\lesssim \epsilon(t)\}|X_1^t)\to 1\qquad P_0\mathrm-\as
\end{equation}
\end{conjecture}
We answer the above in the affirmative in \cref{thm:rates,thm:rates-endpoint-control} where the conditions are made precise. For now, we leave these conditions abstract. 

Much of the difficulty in establishing the above conjecture lies in the considered generality of arbitrary hypothesis spaces, dependent observables, non--additive and time--dependent loss functions, and unattained infima of $P\mapsto \zeta(P)$. 

Whereas the Bayesian setting provides natural choices for $P\mapsto\zeta(P)$ even under misspecification (see \cites{kleijn2006,kleijn2012}), the Gibbs posterior requires more care in identifying a sequence of functions $L_t$ with desirable asymptotic properties. It turns out that the contraction rate $\epsilon(t)$ is intricately related to the rate of convergence of $\Delta_t$ to some limiting function in $L^1(P)$. However, it is well known that the speed of such convergence may indeed be arbitrarily slow \cite{krengel1978}. Moreover, under non--additive and time--varying functions $L_t$, accounting also for prior misspecification where optimal parts of the hypothesis space may be outside the support of the prior, it is non--obvious what $\zeta(P)$ should be. 

In this paper, we utilise, adapt, and extend ideas of \cite{schwartz1965} to consider an analogue of the Kullback--Leibler neighbourhoods for Gibbs posteriors where these instead are informed by the limiting function of $\E_{P_0^{(t)}}\Delta_t$ instead. To account for model and prior misspecification, we subtract the essential infimum of this limit under the prior, even when this may not be attained.

One may further ask, in relation to the results of Miller \cite{miller2021}, whether we require some pointwise almost sure convergence of $\Delta_t(P)$ to some limiting function for asymptotic consistency or contraction rates to hold. We answer this in the negative.

\subsection{Contributions}
The main contributions of this work are the establishment of the first contraction rates for Gibbs posteriors for the nonparametric and arbitrary hypothesis spaces setting, nonadditive loss functions, unattained loss minimisers, misspecified models, and the almost sure case. Our results provide a novel proof technique in posterior asymptotics, building on ideas from real interpolation spaces, the Marcinkiewicz interpolation theorem, Calder\'on--Zygmund theory, and the Laplace transform.
We demonstrate our results in novel settings that were previously unavailable to existing theory, including Dirichlet process %
priors, as well as observables that are not Markov at any order. 

As dominated Bayes posteriors are Gibbs posteriors, our results extend to consistency and contraction rates for (fractional) Bayes under the same desired generality---dependent observables, misspecified models, arbitrary hypothesis spaces, unattained minimisers, and almost sure results. 

Moreover, in the course of the development of this paper, we have explored the work of Shalizi \cite{shalizi2009} as one possible approach, but have found the theory to be incorrect. We expand on this in \cref{apx:on-shalizi} where we provide an explicit counterexample to the contraction rate theorem of \cite{shalizi2009}.

\subsection{Outline of the paper}

In this part, we formally introduce the Gibbs posterior and present our main results on almost sure contraction rates of nonparametric Gibbs posteriors. \cref{sec:notation} develops the necessary notation, and standing assumptions, required when stating the main results on contraction rates and posterior consistency in \cref{sec:main-results}. This includes their relation to the Bayes posterior, and demonstrating applicability of our results to tempered Bayes posterior inference in the same generality. 
We further present an example of our theory on Dirichlet process priors under dependent observables. To contextualise the need for novel proof approaches, we discuss the existing theoretical work on posterior asymptotics in \cref{sec:related-work}, including why previous proof approaches in Bayesian nonparametric contraction rates are unavailable in the considered setting. In \cref{sec:main-ideas}, we outline the rather lengthy proof arguments and explain the main ideas behind these. The full proofs can be found in \cref{part:Gibbs-measure-properties,part:denominator-analysis,part:contraction-analysis}. A discussion on the main results and open problems follows in \cref{sec:discussion}.
Moreover, during the completion of this paper, we have found the results of Shalizi \cite{shalizi2009} to be incorrect, and provide an explicit counterexample and a discussion on this in \cref{apx:on-shalizi}.

\subsection{Acknowledgements}
The author is indebted to Geoff Nicholls for his guidance and support, and for providing invaluable suggestions and corrections to an earlier draft of the paper. 
Further, the author acknowledges support from an EPSRC Industrial Doctoral Landscape Award [UKRI2002] and GSK.

\section{Notation, and standing assumptions}\label{sec:notation}
In the following arguments, for positive integers $l,r\in\N:=\{1,2,\dots\}$, $l\le r$ we denote $[l;r]:=\{l,l+1,\dots,r\}$. Further, we say $a_t\gtrsim b_t$ eventually if there exists some unspecified generic constant $c>0$ such that $a_t\ge c\cdot b_t$ for all $t$ sufficiently large, and denote $a_t\sim b_t$ iff $a_t\gtrsim b_t$ and $b_t\gtrsim a_t$. The maximum is denoted $a\vee b=\max(a,b)$ and the minimum $a\wedge b=\min(a,b)$. The complement of a set $A\in\mathcal T$ is denoted $A^c:=\Theta\setminus A$.

To avoid confusion, we say a function converges $P_0$--almost surely when referring to $(\Omega,P_0)$ and $\Pi_0$--almost everywhere when meaning $(\Theta,\Pi_0)$.

\subsection{General definitions}
It will be convenient to consider, instead of the classical Lebesgue function spaces, a generalisation of these which will allow us to study the types of functions for which the Gibbs posterior is consistent and the corresponding contraction rates more accurately. 
\begin{definition}
We denote Lorentz spaces on a $\sigma$--finite measure space $(X,\mu)$ by $L^{p,\infty}(X,\mu)$, also known as weak $L^p$ spaces, which, for $p\in[1,\infty)$ are equipped with the Lorentz quasi--norm 
\begin{equation}
||u||_{L^{p,\infty}(\mu)}=\sup_{\lambda>0} \lambda\cdot[\mu(\{x\in X: |u(x)|>\lambda\})]^{1/p}.
\end{equation}
For $p=\infty$ the Lorentz space is naturally identified with the $L^\infty$ space.\footnote{It is worth noting that $L^p\subsetneq L^{p,\infty}$, that is Lorentz spaces are strictly larger than their $L^p$ counterparts for $p\in[1,\infty)$ (see e.g. \cite{grafakos2014}).} 
\end{definition}

It is immediate that all $L^p$ functions are in $L^{p,\infty}$, however the reverse need not hold.

\begin{definition}
The average integral is $\fint_E f\dd\mu=\mu(E)^{-1}\int_E f\dd\mu$ for $f\in\mathcal L(X)$ and $E\subseteq X$ with $\mu(E)\in(0,\infty)$.
\end{definition}

\begin{definition}
For a $\mathcal T$--measurable $g:\Theta\to(-\infty,\infty]$, we also define the increasing rearrangements of $g$
\begin{equation}
g_\star(\lambda):=\inf\{v\in[-\infty,+\infty]:\Pi_0(g\le v)\ge \lambda\},\quad\lambda\in(0,1],
\end{equation}
and its corresponding minimal Hardy average
\begin{equation}\quad g_{\star\star}(v):=\frac1v\int_0^vg_\star\dd\,\mathrm{Leb}, \quad v\in (0,1].
\end{equation}
\end{definition}
Further, we put forward the following rigorous notion by what we mean that a statement holds eventually almost surely.
\begin{definition}
A statement holds eventually $P_0$ almost surely, if  $\exists\Omega'\subseteq\Omega$ with $P_0(\Omega')=1$ on which the statement holds for all sufficiently large $t$, being precise here to avoid confusion with the quantifier order.
\end{definition}
\subsection{The limiting function and sieves}

We introduce the limiting function $h:\Theta\to(-\infty,\infty]$, a $\mathcal{T}$--measurable function, defined pointwise as
\begin{equation}\label{asp:1}
    \lim_{t\to\infty}s_t^{-1} \E_{P_0^{(t)}}[-\log R_t(P)]=h(P),\quad\mathrm{with}\quad\inf_{P\in\Theta}h(P)>-\infty
\end{equation}
which we assume throughout to exist, supposing that $\E_{P_0^{(t)}}[-\log R_t(P)]$ is well defined as an extended real value. As $h$ is assumed uniformly lower bounded, we may assume without loss of generality that it is bounded below by 0. 
\begin{remark}
This convergence \eqref{asp:1} should be understood as requiring that $\Delta_t$ converges pointwise in $L^1(\Omega,P_0)$. And crucially we do not require the stronger notion of pointwise $P_0$--almost sure convergence of $\Delta_t$ to $h$. 
\end{remark}

Our results allow for $h$ to be $+\infty$, and we assume for the objective to be well posed that the set of hypotheses where $h$ is infinite does not have full prior measure, as otherwise asymptotic consistency would be impossible. That is, $\Pi_0(I)<\Pi_0(\Theta)=1$ where $$I:=\{P\in\Theta:h(P)=+\infty\}.$$

We work with $h$ rather than a converged loss $L(P)$, as the latter may be infinite \aev when $h$ need not be.
We further introduce notation adapted from large deviation theory,
\begin{equation}
    \mathfrak{h}(A) =\essinf_{P\in A} h(P)
\end{equation}
The essential infimum is taken with respect to $\Pi_0$ throughout, and we define $$\mathfrak{h}^\star:=\mathfrak{h}(\Theta),$$ noting that this requires no topology. By assumption that $\Pi_0(I)<1$ we have, by \cref{lem:h*<infty=>Pi_0(I)<1} that $\mathfrak{h}^\star$ is finite. And to mitigate $\Pi_0(I)\in[0,1)$, we work with a generic constant $M_0>\mathfrak{h}^\star$, which may be taken arbitrary as long as the inequality is strict, throughout.

Further, we assume throughout the paper the existence of a sieve, which for our purposes is a sequence of $\mathcal{T}$--measurable sets $G_t\subseteq\Theta$, such that there exist constants $\alpha,\gamma>0$ such that eventually
\begin{equation}
\Pi_0(G_t)\ge 1- \alpha e^{-s_t\gamma}.
\end{equation}

Lastly, the following assumption allows for a broad class of Gibbs posteriors.
\begin{assumption}\label{asp:lorentz}
There exists a $p\in(1,+\infty]$ and $\kappa\in\R$ for which 
\begin{equation}\label{eqn:lorentz-asp}
\limsup_{t\to\infty}s_t^{-1}\log||c_t||_{L^{p,\infty}(\Pi_0)}\le \kappa\qquad\mathrm{where}\qquad c_t(P):=\E_{P_0^{(t)}}[R_t(P)],
\end{equation}
and we denote the conjugate index $q=p/({p-1})$ throughout.
\end{assumption}
As a consequence of \cref{asp:lorentz}, we have that $R_t\in L^1(\Pi_0)$ eventually \as, by \cref{lem:finite-Z_t-Theta}, and we explicitly allow $R_t\in L^1\setminus L\log L (\Pi_0)$ (\as), in which case $$D_\mathrm{KL}(\Pi_t(\cdot|X_1^t):\Pi_0)=+\infty.$$
\begin{remark}
This assumption is a significant weakening of standard assumptions in the literature of Bayesian and Generalised Bayesian asymptotics. In Bayes, $c_t\le 1$, see \cref{prop:Bayes-c_t<=1}, and likewise the central condition of Syring and Martin \cite{syring2023} supposes that $c_t\le 1$  for all $t$ uniformly over $P$, who further impose a stronger structural condition on $c_t$ and that a minimiser exists. This is also the central condition in Gr\"unwald and Mehta \cite{grunwald2020} who additionally require $R_t\in L\log L(\Pi_0)$ which likewise need not hold.
\end{remark}
\begin{assumption}\label{asp:sieve-gamma}
The sieve $G_t$ can be chosen so that $\gamma>q(\mathfrak{h}^\star+\kappa)$.
\end{assumption}
\begin{remark}
This assumption connects the allowance of \cref{asp:lorentz} to the misspecified setting where $\mathfrak{h}^\star>0$. In fact, for some contraction rate $\epsilon_t$, we may choose $\gamma\equiv\gamma_t= q(\mathfrak{h}^\star+\kappa)+C\epsilon_t$ for some $C>0$ sufficiently large.
\end{remark}

We assume \cref{asp:lorentz} and \cref{asp:sieve-gamma} hold throughout the entire paper unless specifically indicated otherwise.

\section{Main results}\label{sec:main-results}

The purpose of this section is the establishment of contraction rates for Gibbs posterior measures. In this, we discuss assumptions under which the Gibbs measure allows for contraction rates, and assumptions under which we merely have consistency. Attention is also paid to an exploration of the precise minimising set of $h$, including unique minimisation, and how we may nevertheless allow for unattained infima of $h$, as well as the case when the limiting function $h$ is related to (pre--)metrics. Furthermore, we establish how our results transfer to the dominated (tempered) Bayes posterior measures, and provide a main example utilising Dirichlet Process priors and non--additive loss functions in Gibbs posterior measures under highly dependent observables.

\subsection{Contraction rates and consistency for Gibbs posterior measures}

\subsubsection*{Contraction rates for the Gibbs posterior}
Throughout we consider only the rates $\epsilon_t>0$ satisfying
\begin{equation}\label{eqn:epsilon_t-condition}
    \epsilon_t\to0\qquad\mathrm{and}\qquad s_t\epsilon_t\to\infty.
\end{equation}
\begin{remark}
As we establish contraction in ``neighbourhoods'' of $h$, this $\epsilon_t$ corresponds to the $\epsilon_t^2$, and the corresponding condition $t\epsilon_t^2\to\infty$, found in contraction rates for Bayes posteriors (cf. Ghosal and van der Vaart \cite{ghosal2017}), we expand on this in the discussion surrounding \cref{prop:pre-metric-contraction}.
\end{remark}
Given such a rate $(\epsilon_t)_t$ we may define the shells
\begin{equation}\label{eqn:shell-definition-intro}
    S_t^{(j)}:=\{P:j\epsilon_t\le h(P)\wedge M_0-\mathfrak{h}^\star<(j+1)\epsilon_t\}\quad\forall j\in\N.
\end{equation}
\begin{assumption}\label{asp:lorentz-sieve-shell-interaction}
There exists $p_1\in(1,\infty]$, $\alpha_1\in(0,\infty)$, and $c_1\in(0,1)$ such that $P_0$--almost surely it holds for all $t$ sufficiently large and for all $j\in\N$
\begin{equation}\label{eqn:hypothesis-lorentz-sieve-eqn}
    \Pi_0(S_t^{(j)}\cap G_t\cap\{h\wedge M_0>u+\Delta_t\})\le\alpha_1 e^{-p_1s_tu}\quad\forall u\ge c_1j\epsilon_t.
\end{equation}
\end{assumption}
Condition \cref{asp:lorentz-sieve-shell-interaction} assumes that we may decompose the hypothesis space $\Theta$ into disjoint sets $S_t^{(j)}\cap G_t$, on which the deviation of $\Delta_t$ from the limiting function is controlled. We particularly allow variability of the bound across shells, which fundamentally forms a Lorentz type assumption.

Recall also that $M_0>\mathfrak{h}^\star$ is any arbitrary constant that is introduced to control the divergence of $\Delta_t$ on $I$, but the results are independent of the specific choice of $M_0$. The case $p_1<\infty$ is especially suited to the tempered regime where $\beta_t<1$.
\begin{remark}
The endpoint $p_1=+\infty$ is to be understood as an $L^\infty$ bound, where the right hand side of \labelcref{eqn:hypothesis-lorentz-sieve-eqn} is zero. We understand in the remainder, for $c\in\R$, $k\in\N$ and $p_1=+\infty$, that $c/(p_1-k)=0$.
\end{remark}
Moreover, we work also with a sequence of $\mathcal{T}$--measurable, deterministic functions $g_t:\Theta\to(-\infty,+\infty]$ for which there exists a $\bar{b}\ge 0$ such that eventually 
\begin{equation}\label{eqn:g_t-condition-a-e}
    g_t\ge \mathfrak{h}^\star-\bar{b}\epsilon_t
\end{equation}
$\Pi_0$--almost everywhere. 
\begin{remark}
Such a function $g_t$ always exists, take $$\E_{P_0^{(t)}}[\Delta_t]\vee (\mathfrak{h}^\star-\bar{b}\epsilon_t),$$ but the specific choice will depend on the example considered.
\end{remark}
\begin{assumption}\label{asp:averaged-g_t}
Let $(g_t)_t$ satisfy \labelcref{eqn:g_t-condition-a-e} with $\bar{b}\ge 0$, and suppose $\exists a,b\ge0$ and $a'>0$ such that
\begin{equation}\label{eqn:overline-Phi-asp-theorem-1-1}
\overline{\Phi}_{g_t}^{\mathfrak{h}^\star}(a\epsilon_t):=-\log\left(\sup\{v\in(0,1]:g_{t\,\star\star}(v)\le \mathfrak{h}^\star+a\epsilon_t\}\right)\le bs_t\epsilon_t,
\end{equation}
and,
\begin{equation}\label{eqn:summability-asp-F_t(E_t)-theorem-1-1}
    \sum_{t\in\N}P_0\left(\fint_{E_t}(\Delta_t-g_t)\dd\Pi_0>a'\epsilon_t\right)<\infty.
\end{equation}
where $E_t:=\{g_t\le \mathfrak{h}^\star+(2a+\bar{b})\epsilon_t\}$.
\end{assumption}

The bound \labelcref{eqn:overline-Phi-asp-theorem-1-1} is related via the minimal Hardy average of $g_t$ as defined previously, and replaces the more traditional prior mass condition on $\Pi_0$ but satisfies ultimately a similar purpose, namely that the prior does not vanish too quickly on $N_{a\epsilon_t}:=\{h<\mathfrak{h}^\star+a\epsilon_t\}$. The condition \labelcref{eqn:summability-asp-F_t(E_t)-theorem-1-1} is a condition on the fluctuations of $\Delta_t$ over the sets $E_t$. Notice here that we have not assumed an existential quantifier for $E_t$ but have given it. 

The canonical choice of $g_t$ on $E_t$ is $K_t:=\E_{P_0^{(t)}}[\Delta_t]$ as this allows $$\E_{P_0^{(t)}}\fint_{E_t}(\Delta_t-K_t)\dd\Pi_0=0$$ by \cref{lem:exchanging-integration-for-K_t-and-Delta_t}, which is not immediate, but $g_t$ may vary off $E_t$.

\begin{theorem}\label{thm:rates}
Let $(\epsilon_t)_t$ satisfy \labelcref{eqn:epsilon_t-condition}. Suppose that \cref{asp:lorentz-sieve-shell-interaction} holds with $c_1\in(0,1)$ and that \cref{asp:averaged-g_t} is satisfied with $a,b,\bar{b}\ge0$, $a'>0$.

Then, $P_0$--almost surely, for all $M>0$ sufficiently large
\begin{equation}\label{eqn:contraction-rates-conclusion-1}
\Pi_t\left(\{P\in\Theta: h(P)< \mathfrak{h}^\star+M\epsilon_t\}\big|X_1^t\right)\to 1.
\end{equation}
The number $M>0$ is bounded below, that is \labelcref{eqn:contraction-rates-conclusion-1} holds for each $c'\in(c_1,1)$ whenever
\begin{equation}\label{eqn:rates-theorem-M-lower-bound}
    M\ge j_0(c'):=\left\lfloor \frac{\mathfrak{L}}{c'-c_1}\right\rfloor+ 1,
\end{equation}
where $\mathfrak{L}:=2a+a'+b+\bar{b}$. Moreover, eventually $P_0$--almost surely the left hand side of \labelcref{eqn:contraction-rates-conclusion-1} is $1-\mathcal{O}(e^{-s_t\epsilon_t\lfloor M\rfloor(1-c')})$.
\end{theorem}
\cref{thm:rates} demonstrates contraction of the Gibbs posterior around optimal hypotheses at rate just slower than $s_t^{-1}$.

Our next result demonstrates that the choice of $E_t$ in hypothesis \cref{asp:averaged-g_t} is optimal and cannot be improved, up to redefinition of constants, to yield a better rate.  
\begin{theorem}\label{thm:optimal-choice-E_t}
Let $(\epsilon_t)_t$ and $(g_t)_t$ satisfy \eqref{eqn:epsilon_t-condition} and \eqref{eqn:g_t-condition-a-e} respectively.
Consider the following
\begin{gather}
\Pi_0(\{g_t\le \mathfrak{h}^\star+a\epsilon_t\})\ge e^{-bs_t\epsilon_t}\label{eqn:condition-A-thm-optimality-g_t-E_t}\\
\exists E_t\in\mathcal{T}\;\;\mathrm{with}\;\;\Pi_0(E_t)\ge e^{-bs_t\epsilon_t}\;\;\mathrm{and}\;\;\fint_{E_t}g_t\dd\Pi_0\le \mathfrak{h}^\star+a\epsilon_t\label{eqn:condition-B-thm-optimality-g_t-E_t}\\
\overline{\Phi}_{g_t}^{\mathfrak{h}^\star}(a\epsilon_t)\le bs_t\epsilon_t\label{eqn:condition-C-thm-optimality-g_t-E_t}.
\end{gather}
It holds that
\begin{equation}
\mathrm{\labelcref{eqn:condition-A-thm-optimality-g_t-E_t}}\implies \mathrm{\labelcref{eqn:condition-B-thm-optimality-g_t-E_t}}\implies \mathrm{\labelcref{eqn:condition-C-thm-optimality-g_t-E_t}} \overset{((2a+\bar{b})\epsilon_t,\,bs_t\epsilon_t+\log 2)}{\Longrightarrow}\mathrm{\labelcref{eqn:condition-A-thm-optimality-g_t-E_t}}.
\end{equation}
That is, conditions \labelcref{eqn:condition-A-thm-optimality-g_t-E_t,eqn:condition-B-thm-optimality-g_t-E_t,eqn:condition-C-thm-optimality-g_t-E_t} are equivalent up to redefinition of constants under the map $(a\epsilon_t,bs_t\epsilon_t)\mapsto((2a+\bar{b})\epsilon_t,\,bs_t\epsilon_t+\log 2)$.
\end{theorem}
Moreover, any of the above conditions may be used as a verification tool for \cref{asp:averaged-g_t}, including the classical prior mass condition \eqref{eqn:condition-A-thm-optimality-g_t-E_t} of Ghosal and van der Vaart \cite{ghosal2017}, or Bhattacharya, Pati, and Yang \cite{bhattacharya2019} for fractional Bayes posteriors.
We utilise these later in our example.

\subsubsection*{An alternative to \texorpdfstring{\cref{asp:averaged-g_t}}{(H4)}}
Although \cref{thm:rates} will be the main verification tool for the examples we have in mind, the condition \cref{asp:averaged-g_t} may be replaced by the following assumption. 
\begin{assumption}\label{asp:rates-endpoint-ver}
There exists a $p_2\in[0,\infty)$, $\alpha_2\in(0,\infty)$, and $a>0$ for which
\begin{equation}\label{eqn:asp-rates-endpoint}
\sum_{t\in\N}P_0\left(\left\{\omega:\Pi_0(\Delta_t(\omega, \cdot)<\mathfrak{h}^\star+a\epsilon_t)<\alpha_2 e^{-p_2as_t\epsilon_t}\right\}\right)<+\infty.
\end{equation}
\end{assumption}
Hypothesis \cref{asp:rates-endpoint-ver} imposes that the prior eventually places not too quickly vanishing mass on near optimal hypotheses, and complements \cref{asp:averaged-g_t}, with neither implying the other.

\begin{theorem}\label{thm:rates-endpoint-control}
The conclusion \labelcref{eqn:contraction-rates-conclusion-1} of \cref{thm:rates} remains valid under hypothesis \cref{asp:rates-endpoint-ver} instead of \cref{asp:averaged-g_t} with $\mathfrak{L}:=(1+p_2)a$ in the lower bound \labelcref{eqn:rates-theorem-M-lower-bound}.
\end{theorem}

\subsubsection*{Consistency of the Gibbs posterior}
\cref{thm:rates,thm:rates-endpoint-control} have established quantitative contraction results, which necessarily imply for the associated Gibbs posterior measures to be asymptotically consistent. The hypotheses \cref{asp:lorentz-sieve-shell-interaction}, \cref{asp:averaged-g_t}, and \cref{asp:rates-endpoint-ver} are naturally stronger than required for consistency and we may instead establish this under a different set of assumptions where we may weaken \cref{asp:lorentz-sieve-shell-interaction}, and replace \cref{asp:averaged-g_t} or \cref{asp:rates-endpoint-ver} by a rather weak assumption which is satisfied in most practical cases.

\begin{theorem}\label{thm:consistency}
Suppose that there exists $p_1\in(1,\infty]$, $\alpha_1\in(0,\infty)$, $c_1\in(0,1)$ and some nonnegative (random) sequence $\lambda_t\to0$ (\as) satisfying
\begin{equation}\label{eqn:asp-concentration-sieve}
    \Pi_0(G_t\cap \{h(P)\wedge M_0>u+\Delta_t\})\le\alpha_1 e^{-p_1s_tu},\qquad\forall u\ge c_1\lambda_t,
\end{equation}
eventually almost surely, and that for $\Pi_0$--$\aev$ $P$
\begin{equation}\label{eqn:asp-consistency}
P_0\left(\omega:\limsup_{t\to\infty}\Delta_t(\omega;P)\le h(P)\right)=1.
\end{equation}
Pick any set $A\in\mathcal{T}$ for which $\mathfrak{h}(A)>\mathfrak{h}^\star$, then it is $P_0$--almost sure that
\begin{equation}
    \Pi_t(A|X_1^t)\to0.
\end{equation}
\end{theorem}
This naturally presents a notion of frequentist consistency, as $\forall\epsilon>0$  
$$
N_\epsilon^c:=\{P:h(P)\ge\mathfrak{h}^\star+\epsilon\}
$$
satisfies
$$
\mathfrak{h}(N_\epsilon^c)\ge \mathfrak{h}^\star +\epsilon>\mathfrak h^\star,
$$
but presents the stronger result.

It is straightforward that assumption \labelcref{eqn:asp-consistency} is quite general and satisfied in most settings assumed in the literature and in practice, for instance under ergodic theorems.

\subsubsection*{On the assumptions}
In light of the broad generality of the assumptions imposed for Gibbs posterior contraction rates, we further posit some stricter assumptions which satisfy the conditions of the main hypotheses.
\begin{lemma}\label{lem:pointwise-a-s-convergence-to-h}
Suppose that for almost all $P\in\Theta$,
\begin{equation*}
\frac{1}{s_t}\log R_t(P)\to-h(P)\qquad P_0-\as,
\end{equation*}
then assumption \labelcref{eqn:asp-consistency} holds and \cref{thm:consistency} follows under the standing assumptions and \labelcref{eqn:asp-concentration-sieve}.
\end{lemma}
The following imposes control on one--sided uniform deviations of $\Delta_t$ from $h\wedge M_0$, not quite imposing uniform convergence.
\begin{lemma}\label{lem:uniform-convergence-on-sieve}
Suppose that $\exists c_1\in(0,1)$, $M_0>\mathfrak{h}^\star$, and a nonnegative sequence $\lambda_t\to0$ ($\as$) such that $P_0$--almost surely for all $t$ sufficiently large
\begin{equation}\label{eqn:uniform-convergence-on-sieve}
    \esssup_{P\in G_t}\left(\frac1{s_t}\log R_t(P)+(h(P)\wedge M_0)\right)\le c_1\lambda_t
\end{equation}
then \eqref{eqn:asp-concentration-sieve} is satisfied with $c_1,p_1=\infty$, and this $\lambda_t$.
Moreover, if $\lambda_t\le\epsilon_t$ then \cref{asp:lorentz-sieve-shell-interaction} is satisfied with $c_1$ and $p_1=\infty$.
\end{lemma}

\subsection{On attainment of \texorpdfstring{$\mathfrak{h}^\star$}{h*} and pre--metric contraction}

We introduce $\Theta_\star\subset\Theta$, the set of hypotheses on which $\mathfrak{h}^\star$ is attained, that is we define
\begin{equation}
    \Theta_\star:=\{P\in\Theta:h(P)=\mathfrak{h}^\star\}
\end{equation}
which may generally be empty or outside the ``support'' of the prior. This definition is not vacuous, as we have explicitly defined $P\mapsto h(P)$ as a limiting functional, where each $P\in\Theta$ is assigned a precise value rather than $h(P)$ being defined only almost everywhere. For instance, in the Bayesian setting under \iid observables this is the Kullback--Leibler divergence between $P_0$ and the hypotheses $P\in\Theta$, $D_\mathrm{KL}(P_0:P)$.

\subsubsection{Attainment of \texorpdfstring{$\mathfrak{h}^\star$}{h*}}
\subsubsection*{Unattained infimum} 
In key difference to existing work on asymptotics of Gibbs measures, as well as Bayes posteriors, we do not require for some element of $\Theta$ to attain the infimum $\mathfrak{h}^\star$ of $h$, or alternatively this may be attained within $\Theta$ but not within the ``support'' of the prior \cites{kleijn2006,syring2023,miller2021,mcgoff2022}. Our results make no usage of a potential attainment, specifically because we are working with $P_0$ directly rather than some minimiser. In fact, all of our main results are stated in this generality under unattained $\mathfrak{h}^\star$, and occur naturally in the misspecified setting, see for instance \cite{shalizi2009} for such an example.

\subsubsection*{Unique attainment of \texorpdfstring{$\mathfrak{h}^\star$}{h*} and \texorpdfstring{weak$^\star$}{weak*} convergence}
Naturally, to consider Dirac delta measures at some minimiser we require a notion of topology, which we have not assumed so far. We hence suppose this in the following result.
\begin{proposition}\label{thm:weak-star-convergence-unique-min}
Suppose the setting of \cref{thm:consistency}. Suppose also that $\Theta$ admits a topology whose Borel sigma algebra is contained within $\mathcal{T}$, which is first--countable at $\nu^\star\in\Theta$ and is such that the corresponding basis $(U_m)_{m\in\N}$ satisfies
\begin{equation}\label{eqn:asp-separating-basis-dirac-delta-prop}
\mathfrak{h}(\Theta\setminus U_m)>\mathfrak{h}^\star\qquad\forall m\in\N.
\end{equation}
Then it holds $P_0$--\as that
\begin{equation}
\Pi_t(\cdot|X_1^t)\;{\overset{w^\star}{\longrightarrow }}\;\delta_{\nu^\star}.
\end{equation}
\end{proposition}
Here convergence takes place in the weak$^\star$ sense in $\mathcal{P}(\Theta)$, see \cref{def:weak-star_convergence}. Notice that we have not in fact assumed that $\Theta$ itself is first countable, merely that this is satisfied at a single minimising point. 
\begin{remark}
The assumptions here are quite general, and $\nu^\star$ technically need not be unique without a topological structure, but in common settings considered in posterior inference $\Theta$ is regular Hausdorff in which case $\nu^\star$ is unique.
\end{remark}
\subsubsection*{Attainment of the infimum under multiple minimisers}
Naturally, it may be the case that $\Theta_\star$ is not a singleton nor empty, but contains multiple, and possibly an infinite, countable or uncountable, number of minimisers. However, as shown in \cite{coates2026}, even the case of a finite number of minimisers can prove a significant challenge for $k\ge3$ minimisers, and the naive conjecture that weak$^\star$ convergence to elements in a $k-1$ dimensional simplex occurs is non--trivial and a significant open problem under non--atomic priors. %
This presents a central challenge in misspecified Bayes, see the 1966 seminal work of Berk \cite{berk1966}.

\subsubsection{Pre--metric contraction}
It is often desirable to reason with metrics, or more generally pre--metrics, instead of the sublevel sets considered in the main theorems. Under attained infima we may take the pre--metric as the infimal distance of the element $P\in\Theta$ to any of the minimisers in $\nu^\star\in\Theta_\star$. 

However, for unattained infima, this becomes more intricate to formulate, but arises in natural settings, \eg $\Theta$ open but a minimiser may be attained on the boundary. To make this precise without a topology, we instead work with limits of the sublevel sets.
\begin{proposition}\label{prop:pre-metric-contraction}
Suppose there exists some pre--metric $d:\Theta\times\Theta\to[0,\infty]$ and some constant $w>0$ that satisfies
\begin{equation}\label{eqn:pre-metric-asp}
    J(P):=h(P)-\mathfrak{h}^\star\gtrsim d(P,\tilde{\Theta}_\star)^w,
\end{equation}
where 
\begin{equation}
    d(P,\tilde{\Theta}_\star)=
    \begin{cases}
        \inf\limits_{\nu^\star\in\Theta_\star}d(P,\nu^\star), &\Theta_\star\ne \emptyset\\
        \lim\limits_{\delta\downarrow0} \inf\limits_{\nu\in N_\delta}d(P,\nu), &\Theta_\star=\emptyset
    \end{cases}
\end{equation}
then \cref{thm:rates,thm:rates-endpoint-control} hold on the set $\{P:d(P,\tilde{\Theta}_\star)\lesssim (M\epsilon_t)^{1/w}\}$ under their respective assumptions.
\end{proposition}
\begin{remark}
The slightly unwieldy definition of $d(P,\tilde{\Theta}_\star)$ covers the settings discussed previously on the attainment of $\mathfrak{h}^\star$ where our definition allows for unattained infima, and uncountable $\Theta_\star$.
\end{remark}
This proposition makes also precise the relationship of Bayesian inference and concentration in Hellinger neighbourhoods \cite{walker2007} from the Kullback--Leibler divergence and the reason why these neighbourhoods require $t\epsilon_t^2\to\infty$, following the well known relationship between the KL divergence and the Hellinger metric, $D_\mathrm{KL}(\mu:\nu)\gtrsim d_\mathrm H^2(\mu,\nu)$.

Relating this also to the previous proposition on the weak$^\star$ convergence to a Dirac measure for unique minimisers allows us to quantify the rate of this convergence, but we shall not make this precise.

\subsection{The Bayes posterior and Gibbs measure relation}
The relation between the Gibbs posterior and the Bayes posterior is well--known, every dominated Bayes posterior is a Gibbs posterior. This relation immediately licenses the applicability of \cref{thm:rates} to the Bayes setting where we may hence obtain \as contraction rates under misspecified models, dependent observables, unattained minimisers, and arbitrary hypothesis spaces.

In particular, we have that $R_t(P)$ is defined as 
\begin{equation}\label{eqn:R_t-Bayes-def}
    R_t(P):= \frac{dP^{(t),ac}}{dP_0^{(t)}}(X_1^t)
\end{equation}
where $P^{(t),ac}$ is the absolutely continuous part of $P\in\Theta$ at time $t$ with respect to $P_0$. The definition is made rigorous via their Radon--Nikodym derivative. $P^{(t)}$ is understood as the finite image of the measure $P$ at time $t$ as usual, and to make \labelcref{eqn:R_t-Bayes-def} precise, we assume for Bayes that both the model $\Theta$ and the underlying $P_0$ are dominated at each time $t$ by some common $\sigma$--finite reference measure $\mu_t$ on $(\Xi^{(t)},\mathcal{X}^{(t)})$ (cf. Chapters 6 and 8 in \cite{ghosal2017}). Recall here the standing assumption that $R_t$ is assumed $\sigma(X_1^t)\otimes\mathcal{T}$ measurable for all $t$, and we may also note that $dP_0^{(t)}/d\mu_t(X_1^t)$ is almost surely positive implying that $R_t(P)$ is almost surely finite $\forall t$.
It is immediate that $\mathfrak{h}^\star<+\infty$ only if a positive prior measure subset of $\Theta$ has non--null absolutely continuous parts.

The Bayes posterior is then defined as, 
\begin{equation}
    \Pi_t^{(B)}(A|X_1^t)=\frac{\int_A \frac{dP^{(t)}}{d\mu_t}(X_1^t)\,\Pi_0(dP)}{\int_\Theta \frac{dP^{(t)}}{d\mu_t}(X_1^t)\,\Pi_0(dP)}=\frac{\int_A R_t(P)\,\Pi_0(dP)}{\int_\Theta R_t(P)\,\Pi_0(dP)}.
\end{equation}
\begin{remark}
Domination isn't necessary, and we could argue through prior predictives instead \cite{kleijn2021}, however this is standard in the literature under infinite dimensional (misspecified) models with dependent observables \cites{kleijn2006,ghosal2007,walker2007} and makes the relation to tempered Bayes immediate.
\end{remark}
\begin{theorem}\label{thm:Bayes}
For $R_t$ as in \labelcref{eqn:R_t-Bayes-def}, \cref{asp:lorentz} is satisfied with $p=\infty$ and $\kappa=0$, and \cref{asp:sieve-gamma} may be satisfied with any $\gamma>\mathfrak{h}^\star$.

Suppose that assumptions \labelcref{eqn:asp-concentration-sieve,eqn:asp-consistency} hold for $R_t$. Then, for all $A\in\mathcal{T}$ with $\mathfrak{h}(A)>\mathfrak{h}^\star$,
\begin{equation}
    \Pi_t^{(B)}(A|X_1^t)\to0\qquad P_0-\as
\end{equation}

Moreover, let $\epsilon_t\to0$ and $t\epsilon_t\to\infty$, and suppose instead that assumptions \cref{asp:lorentz-sieve-shell-interaction} and either \cref{asp:averaged-g_t} or \cref{asp:rates-endpoint-ver} hold. Then
\begin{equation}
    \Pi_t^{(B)}(\{h<\mathfrak{h}^\star+M\epsilon_t\}|X_1^t)\to 1\qquad P_0-\as
\end{equation}
for all $M>0$ sufficiently large.
\end{theorem}
\begin{remark}
Our results may be used as replacements for those of \cite{shalizi2009}, see \cref{apx:on-shalizi} for a counterexample to the main theorem of \cite{shalizi2009}. 
\end{remark}
\subsubsection*{Fractional Bayes posteriors}
As we allow for varying learning rates $\beta_t$ we may further investigate the relation of our theory to fractional or tempered Bayes posteriors,
\begin{equation}
    \Pi_t^{(B),\beta_t}(A|X_1^t)=\frac{\int_A \bigl(\frac{\dd P^{(t),ac}}{\dd P_0^{(t)}}(X_1^t)\bigr)^{\beta_t}\,\Pi_0(\dd P)}{\int_\Theta \bigl(\frac{\dd P^{(t),ac}}{\dd P_0^{(t)}}(X_1^t)\bigr)^{\beta_t}\,\Pi_0(\dd P)}.
\end{equation}
It is especially noteworthy that the setup yields that the limiting object $h$ agrees in fractional and ordinary Bayes posteriors, that is
$$
h(P)=\lim_{t\to\infty}\frac{1}{t\beta_t}\E_{P_0^{(t)}}\left[-\log \left(\left(\frac{\dd P^{(t),ac}}{\dd P_0^{(t)}}\right)^{\beta_t}\right)\right]=\lim_{t\to\infty}\frac1t\E_{P_0^{(t)}}\left[\log \frac{\dd P_0^{(t)}}{\dd P^{(t),ac}}\right],
$$
in accordance with much of the literature on contraction rates of these \cites{bhattacharya2019,alquier2020,lhuillier2023}. 
Hence, the fractional and Bayes posterior accumulate, if they are consistent, at the same elements of $\Theta$, even in the misspecified setting, and regardless of the learning rate or whether $\beta_t\to0$ or $\beta_t\equiv\beta$. The learning rate will only impact the rate of contraction if $\beta_t$ contracts or diverges, or when assumptions fail for $\beta_t=1$ but hold for some small $\beta_t<1$, in particular \cref{asp:lorentz-sieve-shell-interaction}, (cf. \cites{grunwald2017,bhattacharya2019}).
Furthermore, for $\beta_t\le 1$ \cref{asp:lorentz} holds with $p=\infty$ and $\kappa=0$, and thus \cref{thm:Bayes} applies verbatim, with $s_t$ replacing $t$, to $\Pi_t^{(B),\beta_t}(dP|X_1^t)$.

\subsection{Illustrative example: Gibbs posterior measures under a Dirichlet process prior and ergodic observables}
Suppose $(K,\rho)$ is a compact metric space, and $(\Omega,\mathcal{F},\nu_0,T)$ an ergodic probability measure preserving system generating the sequence of observables, that is for $X_1:\Omega\to K$, and $X_i:=X_1\circ T^{i-1}$ (see \eg \cite{kerr2016}). The stationary measure is the pushforward of $\nu_0$ under $X_1$, $\mu:={X_1}_\sharp \nu_0$, and the finite dimensional marginal $P_0^{(t)}$ is naturally identified as the pushforward of $\nu_0$ under $(X_1,X_1\circ T^1,\dots,X_1\circ T^{t-1})$. 
\begin{remark}
    This approach allows crucially for a multitude of different examples. We may allow for \iid observables by considering $T$ as the one--sided Bernoulli shift, Markovian observables for $T$ the Markov shift, or even observables that are not--Markov at any order, such as the even process considered in \cite{shalizi2009}. We return to this later.
\end{remark}

Set $\Theta=\mathcal{P}(K)$, as compactness of $K$ immediately yields that all first moments of probability measures on $K$ exist, equipped with the Borel $\sigma$--algebra generated by the $1$--Wasserstein metric\footnote{As $K$ is compact and metric, this definition coincides with the classical cylindrical $\sigma$--algebra.} \cite{villani2009}
\begin{equation}
W_1(\nu_1,\nu_2):=\inf_{\pi\in\Gamma(\nu_1,\nu_2)}\int_{K^2}\rho(x,y)\dd\pi(x,y),
\end{equation}
where the infimum is taken over all couplings of $\nu_1,\nu_2\in\mathcal{P}(K)$.

Further, we utilise a Dirichlet process (DP) prior, $\Pi^\mathrm{(DP)}_0=\mathrm{DP}(\alpha, G_0)$ on $\Theta$ \cite{ferguson1974}, which to the best of our knowledge, presents a novel application of Gibbs posteriors in the nonparametric setting. 
In particular, we construct the probability measure $G_0$ on $K$ to be Ahlfors--David regular: for some $C\ge 1$ and $d\ge0$, it holds for all $0<r\le \mathrm{diam}(K)$ and $x\in K$ (see Fraser \cite{fraser2021} or Mattila \cite{mattila2023})
\begin{equation}
    C^{-1}r^d\le G_0(B(x,r)\cap K)\le C r^d.
\end{equation}
Note, we do not impose this condition on $\nu_0$ which can be arbitrarily behaved.

The loss is the $1$--Wasserstein metric to the empirical measure $\mathbb P_t:=\frac1t\sum_{i=1}^t\delta_{X_i}$,
\begin{equation}
\Delta_t(X_1^t;P):=W_1(P,\mathbb P_t)-W_1(\mu,\mathbb P_t).
\end{equation}
The Gibbs posterior measure naturally is defined as
\begin{equation}
    \Pi_t^{(\mathrm{DP)}}(A|X_1^t)=\cfrac{\int_A\exp\{-s_tW_1(P,\mathbb P_t)\}\,\Pi^\mathrm{(DP)}_0(dP)}{\int_\Theta\exp\{-s_tW_1(P,\mathbb P_t)\}\,\Pi^\mathrm{(DP)}_0(dP)},
\end{equation}
where we have picked $\beta_t>0$ satisfying the standing condition. 
Under the above construction, we have the following result.
\begin{theorem}\label{thm:DP-W_1-example}
It holds $\nu_0$--almost surely
\begin{equation}\label{eqn:DP-W_1-example-theorem-consistency}
    \Pi_t^{(\mathrm{DP})}(dP|X_1^t)\overset{w^\star}{\to} \delta_\mu(dP).
\end{equation}
Suppose that for some deterministic $r_t\downarrow0$, eventually $W_1(\mu,\mathbb P_t)\le r_t$ \as and take, for some $A>2$,
\begin{equation}
    \epsilon_t:=A\left(r_t\vee\left(\frac{\log s_t}{s_t}\right)^{1/(d+1)}\right)
\end{equation}
where $d=\mathrm{dim}_H(K)$ is the Hausdorff dimension of $K$, then it holds $\nu_0$--almost surely that
\begin{equation}
    \Pi_t^{(\mathrm{DP})}(\{P\in\Theta: W_1(P,\mu)< M\epsilon_t\}|X_1^t)\to 1
\end{equation}
for all $M>0$ sufficiently large.
\end{theorem}

In proving this, we utilise the following result, which has been significantly extended for our purposes from Lemma 5 in Nguyen \cite{nguyen2013}. 
\begin{proposition}\label{lem:DP-W_1-small-ball-estimates}
Suppose $(K,\rho)$ is a compact metric space, admitting an Ahlfors--David regular probability measure $G_0$. Let $\Pi_0^{(\mathrm{DP})}=\mathrm{DP}(\alpha,G_0)$ for some $\alpha>0$. Then there exists $C_0,\epsilon_0>0$ such that uniformly over $\nu\in\mathcal{P}(K)$ it holds
\begin{equation}
    -\log \Pi_0^{(\mathrm{DP})}(\{P:W_1(P,\nu)\le \epsilon\})\le C_0\,\epsilon^{-d}\log\left(\frac1\epsilon\right),\qquad 0<\epsilon<\epsilon_0.
\end{equation}
Here, $d=\mathrm{dim}_H(K)$ is the Hausdorff dimension of $K$ and equals its Assouad dimension.
\end{proposition}

The last assertion on the Hausdorff dimension being equal to the Assouad dimension is a direct application of Theorem 6.4.1 in Fraser \cite{fraser2021}.

\begin{remark}
The assumption that $G_0$ is AD regular is not necessary, consistency can be proved under $G_0$ having full support on $K$, but results in nicer estimates for the rates and produces an inherent definition of $d$. But while this imposes regularity on the space $(K,\rho)$ it neither reduces the complexity of $\nu_0$, which may carry arbitrary dependencies, nor prevents $\Theta$ being infinite dimensional.
\end{remark}

In particular, our results directly relate the regularity of the underlying compact metric space $K$ to the contraction rate at which the convergence to optimal hypotheses takes place. 

\subsubsection*{The even process}

The even process is an example of a sofic shift process (due to Weiss \cite{weiss1973}, but not to be confused with sofic ergodic theory cf. \cite{kerr2016}), that has a finite hidden Markov representation, but is not Markov at any order \cite{shalizi2009}.

In particular, this process can be represented by a two hidden state, $\{A,B\}$, system with transition matrix
\begin{equation}
    Q_\mathrm{even}=\begin{pmatrix}
        0 & 1\\
        1/2 & 1/2
    \end{pmatrix}
\end{equation}
which produces sequences of $0$s (when changing states) and $1$s (when remaining) where an even number of $0$s is separated by arbitrary sequences of $1$s, \eg $\dots10000111001\dots$
but is itself not--Markov of any order as the current state depends on the entire history.

The metric space $K=\{0,1\}^{\N_0}$ with $\rho(x,y)=2^{-N(x,y)}$ where $N(x,y):=\min\{i\in\N_0: x_i\ne y_i\}$ is compact. Letting $\Omega=K$, we define $X_1:=\mathrm{id}_K$ the identity map on $K$ and $T=\sigma_l$ the left shift. The measure $\mu\in\mathcal{P}(K)$ is then the entire process law under the usual cylinder set $\sigma$--algebra and yields an ergodic \pmp system. %

For the prior, we may take $\alpha>0$ and $G_0:=\mathcal H^1\rest K$, the one--dimensional Hausdorff measure restricted to $K$, which is $1$--AD regular by definition \cite{fraser2021} and coincides with the infinite Bernoulli($1/2$) product measure in this case.

\cref{thm:DP-W_1-example} hence applies, yielding concentration at the entire process measure of the even process.

\section{Related work, and why previous approaches are unavailable}\label{sec:related-work}
In this section, we discuss related work in proving Gibbs and Bayes posterior contraction rates, and why these do not transfer to the settings we consider in this paper. 

\subsection{Comparing results for nonparametric Gibbs posterior measures}
First results for asymptotic consistency, but not contraction rates, were recently established in \cite{miller2021},\cite{mcgoff2022}, and \cite{lopes2022}. In particular, in \cite{lopes2022} Lopes, Lopes, and Varandas established consistency under almost--additive loss functions, Miller \cite{miller2021} allows for general deterministic losses with a global minimiser, under strong assumptions on the asymptotics of this loss sequence, and McGoff, Mukherjee, and Nobel \cite{mcgoff2022} consider general stochastic process observables under compact metric measure hypothesis spaces and additive loss functions. 
Baraud \cite{baraud2024} further studies a restricted subset of Gibbs posteriors, a ``Bayes--like'' conditional posterior, under independent but not necessarily identically distributed observables and a particular nonadditive loss construction under some restrictive assumptions on this loss. \cite{baraud2024} however does not make the connection to Gibbs posteriors. We also mention \cite{winter2026} who utilise the work of Miller \cite{miller2021} for parametric sequential Gibbs posteriors in proving BvM results in a restricted setting on some finite dimensional manifolds.

The results on contraction rates of Gibbs posteriors in the nonparametric setting are significantly sparser. The results established by Gr\"unwald and Mehta \cite{grunwald2020} are only in probability and require strong assumptions on the loss, including lower boundedness, additivity, unique attainment, \iid observables, and $\beta_t\equiv \beta$ fixed, although they allow for $c_t\le e^{\epsilon\beta t}$ uniformly for $\epsilon\ge0$ sufficiently small which presupposes $p=\infty$ in \cref{asp:lorentz}.
Further, Syring and Martin initially establish in probability contraction rates in an applied setting in \cite{syring2020} on which they later build in \cite{syring2023}. In particular, they require \iid observables, although they mention this may be extended to \inid settings, additivity and  attainment under a fixed loss, fixed learning rates, and metric compatibility of the loss.

\subsection{Comparing methods to Bayesian nonparametric contraction rates}

Nonparametric posterior inference requires careful analysis of the asymptotics and the establishment of posterior consistency is necessitated by the seminal converse results of Diaconis and Freedman \cite{diaconis1986} in proving posteriors need not be consistent, which is exacerbated by the nonparametric setting. However, the classical approaches pursued in nonparametric Bayesian inference are unavailable to us in the desired generality we establish in \cref{thm:rates,thm:rates-endpoint-control,thm:consistency}. 

Testing based approaches \cites{ghosal2000,shen2001}, building on the seminal work of Schwartz \cite{schwartz1965}, have enjoyed much development in Bayesian nonparametrics, motivated by the ``necessity of testing'' (we refer to \cite{ghosal2017} for an excellent survey and \cite{kleijn2021} for a more recent account). However, this is unavailable to us as the misspecified setting significantly hinders this approach when considering almost sure rates \cite{kleijn2006} and the nonadditive structure of $L_t$, including this not being a likelihood ratio, prevents this adoption.

Martingale approaches \cite{walker2007} exploit that, in Bayes, $R_t$ is a (super--)martingale bounded in expectation by $1$, \ie $c_t\le 1$, which allows the development of contraction in Hellinger neighbourhoods at $P_0$ in the well specified setting. One can ask for $R_t$ to be a supermartingale satisfying the same properties as in Bayes, \cites{grunwald2020,syring2023} both assume $c_t\le 1$ uniformly, however precluding interesting applications. \cref{asp:lorentz} is strictly more general and prohibits this convenience.

Utilisation of the Donsker--Varadhan variational principle \cites{zhang2006,bhattacharya2019}, following Zellner \cite{zellner1988}, or PAC Bayes \cite{alquier2020}, requires sufficient integrability of $R_t$ by imposing $D_\mathrm{KL}(\Pi_t:\Pi_0)<\infty$, which presupposes $R_t\in L\log L(\Pi_0)$, but we permit $R_t\in L^1\setminus L\log L(\Pi_0)$.

More recent work connecting Wasserstein metrics to contraction rates \cites{dolera2024,rousseau2024}, allowing also for non--dominated Bayes posteriors \cite{camerlenghi2026}, is restricted to particular infinite dimensional probability measures and insufficient for dealing with the generality of our approach.

Lastly, the perhaps most immediately direct comparator to our generality in Bayesian asymptotics is Shalizi \cite{shalizi2009}. However, during the preparation of this paper, we have uncovered that the results of \cite{shalizi2009} do not hold in the claimed generality, and the assumptions, and indeed the proof approaches, are unfortunately insufficient for establishing almost sure contraction rates for Bayes posteriors in the desired setting; we refer the reader to \cref{apx:on-shalizi} for a detailed account, including an explicit counterexample to the contraction rate theorem of \cite{shalizi2009}.

\section{Main ideas, and an outline of the proofs}\label{sec:main-ideas}
Our aim in this section is to present the main ideas of the proofs of the main results developed in the paper. In this, we will pose additional and ancillary results which are of individual interest but primarily serve to illuminate and structure the main argument.
The discussion of the proofs will require us to gloss over some measure theoretic technicalities required for the results, as well as the difficulties provided in dealing with the extended real values. These details are filled in in the subsequent sections.

\subsection{Elements in the proof of \texorpdfstring{\cref{thm:rates}}{Theorem \ref*{thm:rates}}: the sieve decomposition}

As the classical methods to posterior contraction rates developed in the context of Bayesian nonparametric inference are unavailable to us, we employ a different approach which we detail below.

To establish control of the Gibbs posterior on 
$$
N_{M\epsilon_t}:=\{h<\mathfrak{h}^\star+M\epsilon_t\}\in\mathcal{T}
$$
for the rate $\epsilon_t$ we partition the parameter space according to the ``good'' sets $G_t$, the ``bad'' sets $B_t:=\Theta\setminus G_t$, into
\begin{equation}
\begin{aligned}
    \Pi_t(N_{M\epsilon_t}^c|X_1^t)&=\Pi_t(N_{M\epsilon_t}^c\cap G_t|X_1^t)+\Pi_t(N_{M\epsilon_t}^c\cap B_t|X_1^t)\\
    &\le \Pi_t(N_{M\epsilon_t}^c\cap G_t|X_1^t)+\Pi_t(B_t|X_1^t).
\end{aligned}
\end{equation}
where we call $G_t$ ``good'' as we assume explicit control of these through \cref{asp:lorentz-sieve-shell-interaction} as typical in sieve decompositions \cite{barron1999}, and where we employ ideas from Calder\'on--Zygmund theory \cites{calderon1952,grafakos2014}.
We control the convergence behaviour of the Gibbs measure via this decomposition. %

Recalling
\begin{equation}
    Z_t(A):=\int_AR_t\dd\Pi_0,\qquad\forall A\in\mathcal{T}
\end{equation}
which are random variables, the Gibbs posterior is
\begin{equation}
    \Pi_t(A|X_1^t)=\frac{Z_t(A)}{Z_t(\Theta)}.
\end{equation}
\subsubsection{The bad sets and necessity of prior mass control}
We establish that for all $\epsilon>0$ that eventually $(\as)$
\begin{equation}\label{eqn:proof-idea-Bad-set-numerator}
Z_t(B_t)\le \exp(-s_t(\gamma/q-\kappa-\epsilon))
\end{equation}
under just \cref{asp:lorentz} in \cref{prop:numerator-bad-sets}. Appealing to a Hardy--Littlewood inequality, \cref{prop:hardy-littlewood},
\begin{equation}
\int_{B_t}c_t\dd\Pi_0\le \frac p{p-1}||c_t||_{L^{p,\infty}(\Pi_0)}\Pi_0(B_t)^{\frac{p-1}{p}}
\end{equation}
and the Lorentz space structure assumed in \eqref{eqn:lorentz-asp}, estimate \eqref{eqn:proof-idea-Bad-set-numerator} follows by Markov's inequality and the first Borel--Cantelli lemma.
This assumption yields a significant weakening of the central assumptions of Syring and Martin \cite{syring2023} or Gr\"unwald and Mehta \cite{grunwald2020}.

We can show that if the sequence of sets $B_t$ has sufficiently small prior mass, then the second term always vanishes and hence the behaviour of the posterior measure is entirely controlled via the good sets $G_t$.
\begin{theorem}\label{thm:contraction_bad-sets}
Let $B_t\in\mathcal{T}$ be a deterministic sequence of sets satisfying 
\begin{equation}
    \Pi_{0}(B_t)\lesssim e^{-s_t\gamma}
\end{equation}
for some $\gamma>q(\mathfrak{h}^\star+\kappa)$. Assume also that one of assumptions \labelcref{eqn:asp-consistency}, \cref{asp:averaged-g_t}, or \cref{asp:rates-endpoint-ver} holds, then it follows $P_0$--\as that
\begin{equation}
    \Pi_t(B_t|X_1^t)\to 0.
\end{equation}
Moreover, $P_0$--\as for all $t$ sufficiently large, the left hand side is $\mathcal{O}(e^{-s_t\epsilon})$ for all $\epsilon\in(0,\gamma/q-\kappa-\mathfrak{h}^\star)$.
\end{theorem}
\begin{remark}
This theorem needs no rate of convergence. Everything is controlled via the size of the sets $B_t$ vanishing sufficiently quickly which makes this theorem likewise applicable to consistency and contraction rates. Readers familiar with nonparametric Bayesian asymptotics will find this to be the natural analogue of the Gibbs measure to Bayes (see for instance Theorems 6.39 and 6.40 in \cite{ghosal2017}).
\end{remark}
Therefore, under \cref{asp:sieve-gamma}, the posterior on the sieve complement, the bad sets, $\Pi_t(B_t|X_1^t)\to0$ almost surely, under any of the denominator estimates. 

\cref{thm:contraction_bad-sets} immediately yields that prior mass control is necessary for contraction.
\begin{corollary}\label{thm:necessary-prior-mass}
Assume that \cref{asp:lorentz} and \eqref{eqn:asp-consistency} are satisfied. If 
\begin{equation}
\Pi_0(N_{M\epsilon_t})\lesssim e^{-\gamma s_t}
\end{equation}
for some $\gamma>q(\mathfrak h^\star+\kappa)$ and $M>0$, then $P_0$--\as
\begin{equation}
\Pi_t(N_{M\epsilon_t}|X_1^t)\to 0.
\end{equation}
\end{corollary}
\begin{remark}
The case 
\begin{equation}
    e^{-\gamma s_t}\lesssim \Pi_0(N_{M\epsilon_t})\lesssim e^{-bs_t\epsilon_t}
\end{equation}
is an open question.%
\end{remark}

\subsubsection{The good sets: a shell decomposition}

\begin{figure}
\centering
\begin{tikzpicture}[scale=1, every node/.style = {scale=0.75, inner sep=7pt}, line width=0.5pt, line join=round]
\def \domain{plot[closed hobby] coordinates {
      (0,5) (1.2, 5)  (2.2,4) (2.1,2.6) (2.6,1.7) (4,1.4)
      (5.2,0.8) (5.6,-0.1) (4.9,-1.4) (3.5,-1.8) (2.7,-2.2) 
      (2,-3.1) (0.3,-3.4) (-1.4,-2.8) (-3.2,-2.6) (-4.1,-2.2) 
      (-4.4,-0.9) (-4.1,0.7) (-3.8,2) (-4.2,2.8) (-3.9,3.7) (-2.5,3.7)
      (-1,4) 
    }}
\draw[thick] \domain;

\def \setI{plot[closed hobby] coordinates {
  (3.5,0.8) (4,0.5)  (4.4,0.4) (4.5,-0.2) (4.9,-0.5) (4.8,-1)
  (4.2,-0.9) (3.3,-1.1) (2.7,-0.9) (2.8,-0.2) (3,0.5)
}}

\def\Nupper{plot[closed hobby] coordinates {
  (-2,2.4) (-1.2,2.4)  (-1.4,2) (-1.2,1.5) (-1.6,1.3) (-2,1.5)
  (-2.4,1.4) (-2.6,1.7) (-2.6,2.1) (-2.8,2.4)
}}

\def\Nlower{plot[closed hobby] coordinates {
  (0,0.4) (0.5,0.5)  (0.8,0) (0.6,-0.6) (0,-0.3) (-0.7,-0.5)
  (-0.4,0) (-0.4,0.4) 
}}

\def\SOneSeparate{plot[closed hobby] coordinates {
  (-2.5,-0.5) (-2,-0.8)  (-1.8,-1.2) (-2.3,-1.3) (-2.8,-1.4) (-3.2,-1)
  (-2.8,-0.8) 
}}

\def\SOnelower{plot[closed hobby] coordinates {
  (0,0.8) (0.8,0.9)  (1.2,0) (0.8,-0.8) (0,-0.6) (-1,-1)
  (-0.8,0) (-1,0.7) (-0.6,1)
}}

\def\SOneupper{plot[closed hobby] coordinates {
  (-2,2.7) (-1.2,2.8)  (-0.8,2.3) (-1,1.7) (-1.2,1.2) (-1.7,1.1)
  (-2.3,1.3)  (-2.8,1.2) (-3,2) (-2.8,2.8)
}}

\def\STwoseparate{plot[closed hobby] coordinates {
  (-2.5,-0.2) (-2,-0.4)  (-1.7,-1) (-1.9,-1.5) (-2.5,-1.5) (-3.3,-1.7)
  (-3.8,-1.4)  (-3.5,-0.7) (-3,-0.7) 
}}

\def\STwominimisers{plot[closed hobby] coordinates {
  (-2,3.1) (-1,3)  (-0.5,2.4) (-0.6,1.4) (0.3,1.2) (1.3,1)
  (1.3,-0.7)  (1,-1) (-0.7,-1.6) (-1.3,-1) (-1.1,0) (-1.3,0.7)
  (-2.5,0.5) (-3.3,1.3) (-3.2,2) (-3.6,2.6) (-3,3.2)
}}

\def\SThree{plot[closed hobby] coordinates {
  (-2,3.4) (-0.4,3.3)  (0.3,2.4) (0.3,1.5) (1.4,1.5) (2,0.4)
  (1.7,-1)  (1.5,-2) (0.5,-2.2) (0,-2.7) (-1.4,-2) (-3,-2.3)
  (-4,-1.8) (-4.2,-1) (-3.8,-0.4) (-3,-0.2) (-3,0.3) (-3.6,0.6)
  (-3.5,1.7) (-3.7,2.3)  (-3.5,3.5)
}}

\def\Sdagger{plot[closed hobby] coordinates {
  (-1.5,4.5) (-0.5,3.8)  (0.5,3.4) (1,2) (2,1.5) (2.5,0.5)
  (2,-1)  (2.3,-2) (1.4,-2.7) (0,-3) (-0.8,-3.6) (-2,-3.2) (-4,-3) (-4.5,-2) (-4.5,-0.8) (-4,-0.3) (-3.2,0.1) (-3.6,0.4) (-4.1,0.8)
  (-4.3,1.8) (-4.5,3.5) (-3.5,4.3) (-2,5)
}}
\begin{scope}[transparency group, fill opacity=0.35]
\fill[ibm2] \domain;
\begin{scope}
  \clip \domain;
  \fill[ibm1] \Sdagger ;
\end{scope}

\fill[ibm0] \SThree;
\fill[ito3] \STwoseparate; \fill[ito3] \STwominimisers;
\fill[ito2] \SOneSeparate; \fill[ito2] \SOnelower; \fill[ito2] \SOneupper;
\fill[ito1] \Nlower; \fill[ito1] \Nupper;
\fill[ibm3]  \setI;
\end{scope}

\begin{scope}
  \clip \domain;
  \draw[line width=0.7pt] \Sdagger;
\end{scope}
\draw[line width=0.7pt] \SThree;
\draw[line width=0.7pt] \STwoseparate;  \draw[line width=0.7pt] \STwominimisers;
\draw[line width=0.7pt] \SOneSeparate;  \draw[line width=0.7pt] \SOnelower;  \draw[line width=0.7pt] \SOneupper;
\draw[line width=0.7pt] \Nlower;     \draw[line width=0.7pt] \Nupper;
\draw[line width=0.7pt, dash pattern=on 12pt off 3pt] \setI;
\draw[thick] \domain;

\node[align=center] at (-2,2) {$N_{\epsilon_t}$}; \node[align=center] at (0,0) {$N_{\epsilon_t}$};

\node[align=center] at (-2.5,-1) {$S_t^{(1)}$}; \node[align=center] at (-0.6,0.6) {$S_t^{(1)}$}; \node[align=center] at (-1.1,2) {$S_t^{(1)}$};

\node[align=center] at (-3.5,-1.2) {$S_t^{(2)}$}; \node[align=center] at (-2.4,0.9) {$S_t^{(2)}$};

\node[align=center] at (1.2,-1.6) {$S_t^{(3)}$};

\node[align=center] at (1,4) {$S_t^{(j^\dagger(t))}$}; \node[align=center] at (-3.8,0.05) {$S_t^{(j^\dagger(t))}$};

\node[align=center] at (5,-1.7) {$\Theta$};

\node[align=center] at (3.6,-0.2) {$I$};

\end{tikzpicture}
\caption{Illustration of an example shell decomposition of \cref{asp:lorentz-sieve-shell-interaction} partitioning the hypothesis space $\Theta$ into the shells $S_t^{(j)}$ and $N_{\epsilon_t}$ of $h:\Theta\to[0,\infty]$ according to $\epsilon_t$ and $M_0$. This forms a topographical map according to the different wells of $h$. Note, $I\subseteq S_t^{(j^\dagger(t))}$ because of the truncation at $M_0>\mathfrak{h}^\star$.}
\label{fig:shell-decomposition}
\end{figure}
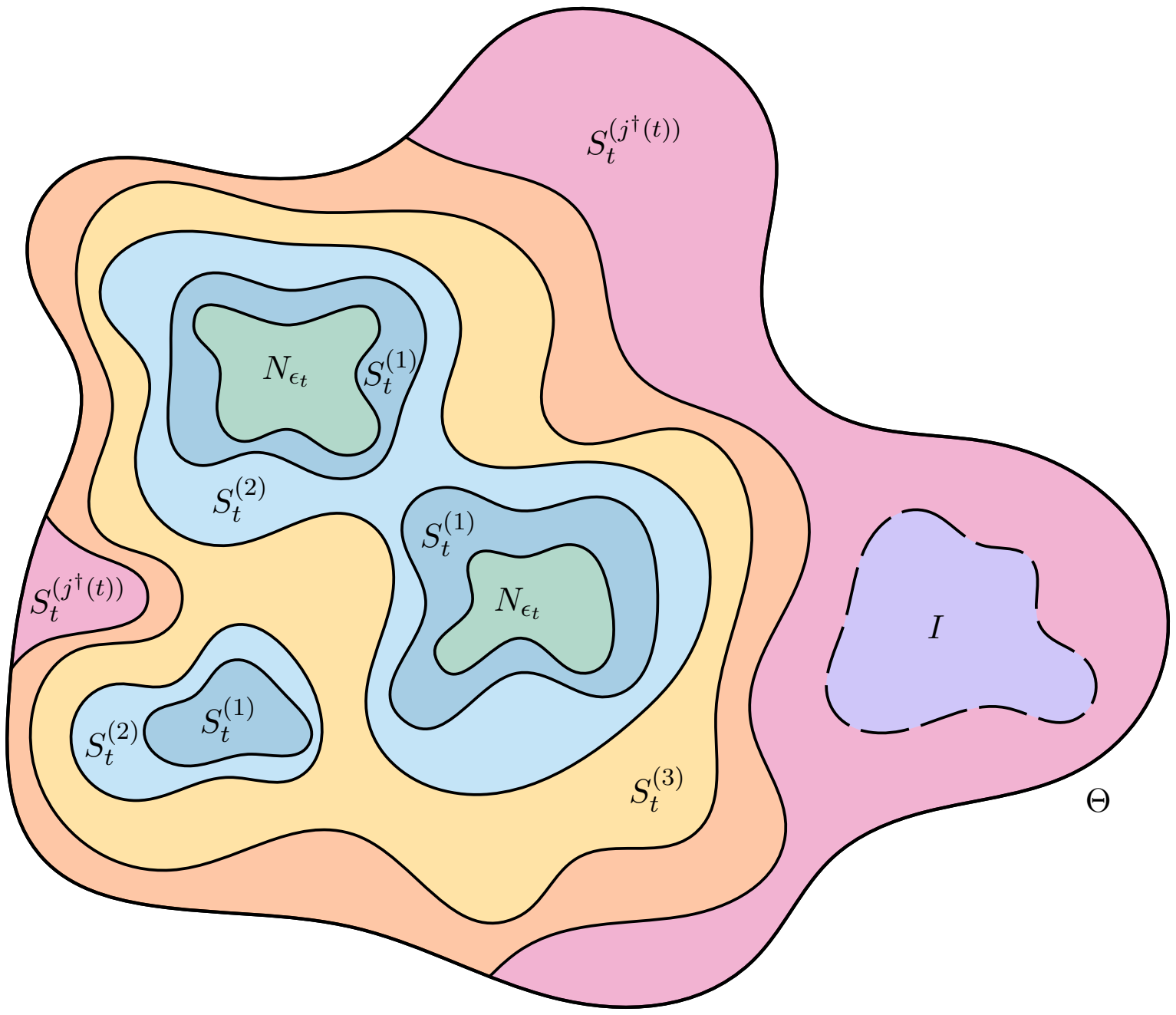

On the sieve, we likewise have to establish control of the numerator, $Z_t(N_{M\epsilon_t}^c\cap G_t)$. We may partition the sieve via the shells \eqref{eqn:shell-definition-intro} $$S_t^{(j)}:=\{P:j\epsilon_t\le h(P)\wedge M_0-\mathfrak{h}^\star<(j+1)\epsilon_t\},\quad j=1,2,\dots,$$ similar to those of \cite{ghosal2017}; see \cref{fig:shell-decomposition} for an illustration. These partition $\Theta$, \cref{lem:shells-partition}, and are empty for $$j>j^\dagger(t):=\lfloor(M_0-\mathfrak{h}^\star)\epsilon_t^{-1}\rfloor\vee 1.$$ 
Hence, at each time $t$ we have a finite partition with $j^\dagger(t)\sim\epsilon_t^{-1}\to\infty$.
\begin{remark}
We stress the importance of $M_0>\mathfrak{h}^\star$ as a useful tool in the proof, irrespective that the conclusions are independent of the specific choice. Particularly, for all $t$, $$I\cap G_t\subseteq S_t^{(j^\dagger(t))}\cap G_t$$ which may generally be non--empty and have positive prior measure, but all hypotheses in $I$ cause $\Delta_t$ to misbehave, and hence $M_0$ allows control on the behaviour of the Gibbs posterior on $I$.

On the remainder of the sieve, $S_t^{(j)}$, for $j\in[1;j^\dagger(t)-1]$, the truncation $M_0$ is invisible.
\end{remark}

Now, \cref{asp:lorentz-sieve-shell-interaction} yields the required control, where we allow variation across the different shells. 
Defining, the random one sided deviation of $\Delta_t$ from $h\wedge M_0$ by
\begin{equation}
U_t(P):=h(P)\wedge M_0- \Delta_t(X_1^t;P)
\end{equation}
allows us to write 
\begin{equation}
R_t=e^{s_tU_t}e^{-s_t(h\wedge M_0)}.
\end{equation}
We establish control on the integral of $R_t$ over $S_t^{(j)}\cap G_t$ by bounding the deterministic part $e^{-s_t(h\wedge M_0)}$ and the random part $e^{s_tU_t}$ separately. The deterministic one is easily bounded by the essential supremum since $e^{-s_t(h\wedge M_0)}\in L^\infty$.

For the random part, we utilise ideas from Calder\'on--Zygmund theory, 
\begin{equation}
\Psi_t^{(j)}:=e^{s_tU_t}\mathds{1}_{S_t^{(j)}\cap G_t}
\end{equation}
where we proceed analogously to the proof of the Marcinkiewicz interpolation theorem (see Calder\'on \cite{calderon1966}, or Theorem 1.3.2 in Grafakos \cite{grafakos2014}). 

Take $\lambda_t^{(j)}=e^{s_tc_1j\epsilon_t}$, 
\begin{equation}
\Psi_t^{(j)}=\Psi_t^{(j)}\mathds{1}_{\{\Psi_t^{(j)}\le \lambda_t^{(j)}\}}+\Psi_t^{(j)}\mathds{1}_{\{\Psi_t^{(j)}> \lambda_t^{(j)}\}}.
\end{equation}
Here, the first part is in $L^\infty$, and the latter in $L^{p_1,\infty}\subset L^1$. More precisely, after a layer--cake argument and utilising the control afforded by \cref{asp:lorentz-sieve-shell-interaction}, integrating the latter over $S_t^{(j)}\cap G_t$ equals integrating
\begin{equation}
\Pi_0\left(\Psi_t^{(j)}>\lambda\right)\le \alpha_1\lambda^{-p_1},\qquad \lambda\ge \lambda_t^{(j)}
\end{equation}
which is integrable over $\lambda\ge \lambda_t^{(j)}$ precisely because $p_1>1$. 

The shells, forming an analogue to the dyadic cubes of Calder\'on--Zygmund decompositions, have bounds geometric in $j$, hence summing $Z_t(S_t^{(j)}\cap G_t)$ establishes
the conclusion of the numerator bound on the good sets, \cref{prop:refined-sieve-estimate}, 
\begin{equation}\label{eqn:proof-discussion-numerator-on-sieve-estimate}
Z_t(G_t\cap\{h\ge \mathfrak{h}^\star+j\epsilon_t\})\lesssim e^{-s_t\mathfrak{h}^\star}e^{-s_t\epsilon_tj(1-c_1)}\quad \forall j\in[1;j^\dagger(t)]
\end{equation}
(eventually \as).

\subsection{Denominator estimates}
In the statement of the contraction rates we have posited two separate hypotheses under which we may demonstrate that contraction rates hold. %
In particular, these serve to prove the intermediate hypothesis that almost surely
\begin{equation}\label{eqn:proof-explanation-eta_t-epsilon_t-bound}
    \limsup_{t\to+\infty}\frac{\eta_t}{\epsilon_t}\le \mathfrak{L}<\infty,
\end{equation}
for
\begin{equation}\label{eqn:eta_t-definition-sec-1}
    \eta_t:=-\frac{1}{s_t}\log\left(e^{s_t\mathfrak{h}^\star}Z_t(\Theta)\right)=-\mathfrak{h}^\star-\Lambda_t,
\end{equation}
where $\Lambda_t=s_t^{-1}\log Z_t(\Theta)$,
which turns out to be a convenient object of study. %

This establishes further 
$$
\limsup_{t\to\infty}\eta_t\le 0\quad\mathrm{a.s.}
$$
which allows the application of \cref{thm:contraction_bad-sets} where we have not made assumption \labelcref{eqn:asp-consistency}. 
Taken together with the numerator estimate \labelcref{eqn:proof-discussion-numerator-on-sieve-estimate}, this may be used to establish control on $\Pi_t(N_{j_0(c')\epsilon_t}^c\cap G_t|X_1^t)$ through careful measure theoretic arguments in \cref{thm:rates-on-sieve-under-eta/epsilon-asp}.

Now to establish \eqref{eqn:proof-explanation-eta_t-epsilon_t-bound}, \cref{asp:averaged-g_t} has two parts, a bound on $\overline{\Phi}_{g_t}^{\mathfrak{h}^\star}(a\epsilon_t)$, as well as summability of the fluctuations of the averaged difference over the set $E_t$ between $\Delta_t$ and $g_t$. 

The first is indeed related to the prior mass control of classical Bayesian contraction rates, in particular this relationship is made precise through \cref{thm:optimal-choice-E_t,lem:bathtub-sandwich}, converting the averaged and optimised hypothesis to the sublevel sets. This is especially inspired by the ``Bathtub principle'' in analysis (see Chapter 1.14 in \cite{lieb1997}) as well as interpolation and decreasing rearrangements theory (cf. \cite{calderon1966}, \cite{bennett1988}) applied here to the increasing rearrangements of $g_t$ and the associated minimal Hardy average. In particular \cref{lem:bathtub-sandwich} demonstrates that
\begin{equation}
-\log \Pi_0(g_t\le\mathfrak{h}^\star+(2a+\bar{b})\epsilon_t)-\log2\le \overline{\Phi}_{g_t}^{\mathfrak{h}^\star}(a\epsilon_t)
\end{equation}
implying
$$
\Pi_0(E_t)\ge\Pi_0(g_t\le\mathfrak{h}^\star+(2a+\bar{b})\epsilon_t)\ge \frac12 e^{-bs_t\epsilon_t}
$$
using also the definition of $E_t$. The condition implied by \cref{lem:bathtub-sandwich} is itself inspired by and a generalised version of the small ball deviations of van der Vaart and van Zanten \cite{vanderVaart2008}. This principally ensures that the sets retain sufficient mass to not vanish as stipulated for sufficiently small sets in \cref{thm:contraction_bad-sets}. In fact, we establish in \cref{thm:optimal-choice-E_t} that the choice of $E_t$ is optimal and optimising this set does not improve the rate of convergence. 

Before we may utilise the second part of the hypothesis \cref{asp:averaged-g_t}, we first need to derive a lower bound on $\log Z_t(\Theta)$ pathwise for each $\omega\in\Omega$, relating this to the average integral over $E_t$ of $\Delta_t(\omega;\cdot)$. In fact, we have 
$$
\log Z_t(\Theta)(\omega)\ge \log \Pi_0(E_t)-s_t\fint_{E_t}\Delta_t(\omega;\cdot)\dd\Pi_0
$$
by a pathwise application of Jensen's inequality, being careful about extended real values, and which we may relate with the previous estimate, to the second condition of \cref{asp:averaged-g_t}. 
Some cumbersome manipulations of this, and application of the first Borel--Cantelli lemma, yields the desired eventually \as statement. This argument is made precise in \cref{thm:pathwise-denominator-contraction}.

Through similar methods, we may make the following observation relating to sharp rates. Denote
\begin{equation}
    A_t(P):=-\frac1{s_t}\log c_t(P),
\end{equation}
and
\begin{equation}
\mathcal{E}_{\Pi_0}[A_t](s_t):=-\log ||e^{-A_t}||_{L^{s_t}(\Pi_0)}.
\end{equation}
We are further led to a natural lower bound on the rate at which the denominator increases.
\begin{theorem}\label{thm:converse-to-rates}
Let $\delta_t\ge0$ satisfy $\sum_te^{-s_t\delta_t}<\infty$. Then, 
\begin{equation}
    P_0\left(\Lambda_t>-\mathcal E_{\Pi_0}[A_t](s_t)+\delta_t\right)\le e^{-s_t\delta_t}
\end{equation}
which is summable. Moreover, eventually $P_0$--\as
\begin{equation}\label{eqn:thm-rates-converse-eta_t-lower-bound}
\eta_t\ge \mathcal{E}_{\Pi_0}[A_t](s_t)-\mathfrak{h}^\star-\delta_t\ge \essinf_{P\in\Theta} A_t(P)-\mathfrak{h}^\star+\frac12\mathcal{J}_{A_t}(s_t)-\frac{\log2}{s_t}-\delta_t,
\end{equation}
where 
\begin{equation}
 \mathcal{J}_{A_t}(s_t):=\inf_{u>0}\left\{u-\frac{\log\Pi_0(A_t\le\mathfrak{a}_t^\star+u)}{s_t}\right\}
\end{equation}
if $\essinf_{P\in\Theta} A_t(P)=\mathfrak a_t^\star>-\infty$.
\end{theorem}

\subsection{On the proofs of \texorpdfstring{\cref{thm:rates-endpoint-control,thm:consistency}}{Theorems \ref*{thm:rates-endpoint-control} and  \ref*{thm:consistency}}}
We also discuss some of the elements in the proofs of the two other main results, \cref{thm:rates-endpoint-control,thm:consistency}.%

\subsubsection{On the use of \texorpdfstring{\cref{asp:rates-endpoint-ver}}{(H\ref*{asp:rates-endpoint-ver})}}
The proof of \cref{thm:rates-endpoint-control} proceeds quite similar to that of \cref{thm:rates}, except for the use of \cref{asp:rates-endpoint-ver} in establishing estimate \labelcref{eqn:proof-explanation-eta_t-epsilon_t-bound}, which we establish in \cref{thm:denominator-via-lorentz-bound}. The main idea of the assumption is that we exhibit a witness on which a subset of the posterior increases sufficiently quickly so that this may serve as a lower bound on $Z_t(\Theta)$ and serves as an eventual \as upper bound on $\eta_t/\epsilon_t$.

The proof relies on first recognising the relationship of $Z_t(\Theta)$ and the Lorentz space $L^{1,\infty}(\Pi_0)$, of which $R_t$ is \as a member as $$R_t\in L^1(\Pi_0)\subset L^{1,\infty}(\Pi_0),$$ via the estimate 
\begin{equation}
    Z_t(\Theta)\ge e^{-s_t(\mathfrak{h}^\star+u)}\Pi_0(\Delta_t<\mathfrak{h}^\star+u)
\end{equation}
where $u>0$ is arbitrary, which we establish in \cref{lem:L^1-infty-bound-on-Z_Theta}. This behaves like an endpoint of the interpolation scale $L^{p,\infty}$ of the sieve bounds, considering instead $p=1$.

Making use of this result, taking $u=a\epsilon_t>0$, and \cref{asp:rates-endpoint-ver} allows us to conclude the estimate \labelcref{eqn:proof-explanation-eta_t-epsilon_t-bound} after a Borel--Cantelli argument. 

\subsubsection{On the proof of \texorpdfstring{\cref{thm:consistency}}{Theorem \ref*{thm:consistency}}}
The proof of the consistency results proceeds similarly to that of the contraction rates, but allows for simplification as we need not consider the rate at which the convergence takes place. 
Likewise to the proof of \cref{thm:rates}, for $A\in\mathcal T$ with $\mathfrak h(A)>\mathfrak h^\star$
\begin{equation}
    \Pi_t(A|X_1^t)\le \Pi_t(A\cap G_t|X_1^t)+ \Pi_t(B_t|X_1^t)
\end{equation}
where the second term vanishes \as as $t\to\infty$ by \cref{thm:contraction_bad-sets} as before. The first term now does not need to deal with the intricacies associated with the prior mass of the sets $N_{M\epsilon_t}$ vanishing, and in fact we have, by \cref{lem:h(A)>h*=>Pi_0(A)<1}, that $\mathfrak{h}(A)>\mathfrak{h}^\star\implies\Pi_0(A)<1$. That is the prior has positive mass on the set we aim to concentrate on.

As before $U_t(P):=(h(P)\wedge M_0)-\Delta_t(X_1^t;P)$, allows us to formulate $$R_t(P)=e^{s_tU_t(P)-s_t(h(P)\wedge M_0)}.$$
The interpolation of the random function $e^{s_tU_t}$ as on the shell before, yields the control on the numerator on the good sets. 

On the denominator, we place the conceptually weak assumption \labelcref{eqn:asp-consistency}, satisfied in most settings of interest, for instance as soon as an ergodic theorem applies, or a generalised asymptotic equipartition property \cite{algoet1988}.
This result took shape mainly as part of our effort in investigating the results of \cite{shalizi2009}, where one of the central problems is the assignment of a rate to the exceptional sets provided by Egorov's theorem. We make the requirements and dependence of the Egorov sets on $\omega$ and the random onset times precise in the proof of \cref{thm:qualitative-denominator-floor}, establishing
\begin{equation}
    \liminf_{t\to+\infty}\Lambda_t\ge -\mathfrak h^\star\qquad P_0-\as
\end{equation}

That 
\begin{equation}
\Pi_t(A\cap G_t|X_1^t)=\frac{Z_t(A\cap G_t)}{Z_t(\Theta)}\to0
\end{equation}
(\as) may then be established by utilisation of these bounds.

\section{Discussion and open questions}\label{sec:discussion}

We have introduced novel contraction rate theorems on Gibbs posteriors in broad generality, including arbitrary hypothesis spaces, nonadditive loss functions, dependent observables, unattained minimisers, possibly infinite limiting loss differences, and make no relation to the Donsker-Varadhan variational formulation which may generally be inappropriate for certain settings. The proof techniques introduced provide new tools for nonparametric posterior asymptotics through interpolation theory and Lorentz spaces. Our results extend to consistency, and weak$^\star$ convergence to the Dirac delta at a unique minimiser, as well as contraction in pre--metrics. Our results provide almost sure convergence of the Bayes posterior in the same generality, and may replace the results of Shalizi \cite{shalizi2009} for which we have identified counterexamples.

While we have phrased our results in terms of observables and an underlying data generating process $P_0$, this is not a requirement of our theory, and nothing assumes that the loss functions make use of the observables. And we may therefore work with a fully deterministic sequence of differences $\Delta_t(P)$ instead, similar in spirit to the results of Miller \cite{miller2021}. Neither do we require the elements of $\Theta$ to be probability measures, we have merely phrased the results in this notation to highlight the novel theoretical contributions to nonparametric Gibbs posteriors but these may for instance be functions in some function space that may not be real valued.

The main open problems left for Gibbs posteriors are nonparametric Bernstein--von Mises results, which even in the Bayes setting are famously hard and impossible in strong topologies \cite{freedman1999}, however Nickl \cite{nickl2020jems} and Giordano and Kekkonen \cite{giordano2020} give some BvM results for inverse problems. In fact, Bayesian inverse problems are Gibbs posteriors on separable Hilbert spaces with Gaussian priors \cite{stuart2010}. While we think that some of the methods presented in the paper can aid in this development, the generality supposed throughout, specifically arbitrary hypothesis spaces through both assuming no topology and not specifying a prior measure, generally does not allow for BvM results and we would have to significantly restrict the setting to attempt these.

Furthermore, we require deterministic learning rates $\beta_t$, as necessitated by our proof approach, however it is natural to consider learning rates informed by previous observables, for instance by performing Bayesian inference on $\beta_t$ \cite{lee2026}.

Another natural extension to the present work is the convergence behaviour of the posterior measure under the non--unique minimisers which, as demonstrated by \cites{berk1966,coates2026}, is difficult even under multiple but finite minimisers. 

\subsection{Outline of the remainder of the paper}

The remainder of the paper is dedicated to proving the main results, which we structure into several parts. %
In \cref{part:Gibbs-measure-properties} we establish tools and properties of the numerator on and off the sieve in \cref{sec:numerator-on-sieves,sec:numerator-off-sieves}. \cref{part:denominator-analysis} provides tools for proving sufficiently fast denominator estimates in \cref{sec:denominator-estimates}. In \cref{part:contraction-analysis} we utilise the results of \cref{part:Gibbs-measure-properties,part:denominator-analysis} to provide the proofs of the main results, in addition to verifying the assumptions. \cref{sec:consistency-proof} establishes first the consistency results, after which we provide the proof of the main theorems on contraction rates in \cref{sec:rates-proof}. The relationship of Gibbs and Bayes posteriors is proved in \cref{sec:Bayes-proof}. \cref{apx:measure-theory-background} presents some measure theoretic definitions and properties used throughout, and in \cref{apx:on-shalizi} we discuss the results of \cite{shalizi2009} and provide a simple counterexample demonstrating that the stated contraction rate theorem may fail.

\part{Properties of the Gibbs posterior measure}\label{part:Gibbs-measure-properties}

\section{Some elementary properties}
This section introduces some elementary properties used throughout the proofs. 
\subsection{On the finiteness of \texorpdfstring{$\mathfrak{h}^\star$}{h*}}
We make the following useful observations.
\begin{lemma}\label{lem:h(A)>h*=>Pi_0(A)<1}
Let $A\in\mathcal{T}$ satisfy $\mathfrak{h}(A)>\mathfrak{h}^\star$, then $\Pi_0(A)<1$. Moreover, if $\Pi_0(A)=0$ then $\mathfrak{h}(A)=+\infty$.
\end{lemma}
\begin{proof}
Suppose first that $\mathfrak{h}(A)<\infty$ and denote $\delta:=(\mathfrak{h}(A)-\mathfrak{h}^\star)/2>0$ by assumption.
Then, using the definition of the essential infimum,
\begin{equation}\label{eqn:essinf-A-def}
\mathfrak{h}(A):=\sup\{b\in\R\cup\{+\infty\}:\Pi_0(A\cap\{h<b\})=0\}
\end{equation}
from which it follows that $h(P)\ge\mathfrak{h}(A)$ for $\aev$ $P\in A$.
Furthermore, since $\mathfrak{h}^\star=\mathfrak{h}(\Theta)$, by definition, $\forall\epsilon>0$
\begin{equation*}
    \Pi_0(h<\mathfrak{h}^\star+\epsilon)>0.
\end{equation*}
Therefore, taking $\epsilon=\delta$, we have that $\Pi_0(h<\mathfrak{h}^\star+\delta)>0$ and hence $$\{h<\mathfrak{h}^\star+\delta\}\subset (\Theta\setminus A) \cup N,$$ where $N$ is a $\Pi_0$--null set. It follows that
\begin{align*}
1=\Pi_0(\Theta)&=\Pi_0(\Theta\setminus A)+ \Pi_0(A)\\
&\ge \Pi_0(h<\mathfrak{h}^\star+\delta)+\Pi_0(A)>\Pi_0(A)
\end{align*}
yielding the first result for $\mathfrak{h}(A)<\infty$. 
Consider now $\mathfrak{h}(A)=+\infty$ and pick any $\delta>0$, then the same argument as before establishes this case.
The moreover part is an immediate consequence of definition \labelcref{eqn:essinf-A-def}, as the bound is satisfied for all $b\in\R$ since $A$ is null. 
\end{proof}

\begin{lemma}\label{lem:h*<infty=>Pi_0(I)<1}
$\mathfrak{h}^\star<+\infty$ if and only if $\Pi_0(I)<\Pi_0(\Theta)$.
\end{lemma}
\begin{proof} 
Suppose first that $\Pi_0(I)<1$, then there exists a positive measure set on which $h<\infty$. Thus, there must exist a constant $b<+\infty$ such that $\{h\ge b\}$ does not have full $\Pi_0$ measure, and hence the essential infimum of $h$ over $\Theta$ is finite.  

For the forward, suppose for a contrapositive that $\Pi_0(I)=\Pi_0(\Theta)$, then by definition of the essential infimum, for all $a<+\infty$ it holds that $h\ge a$ $\Pi_0$--\aev, and thus $\mathfrak{h}^\star=\mathfrak{h}(\Theta)=\infty$. The conclusion follows.
\end{proof}
\subsection{A quantifier reversal lemma}
The following is essentially Lemma 1 in \cite{shalizi2009}.
\begin{lemma}\label{lem:quantifier-reversal}
Let $(X,\mathcal{A})$ and $(Y,\mathcal{B})$ be two arbitrary measurable spaces. 
Suppose that $Q\subseteq X\times Y$ is jointly measurable, and define $Q_x=\{y:(x,y)\in Q\}$ and $Q_y=\{x: (x,y)\in Q\}$. 

Then, if for some probability measures $\mu$ and  $\nu$ on $(X,\mathcal{A})$ and $(Y,\mathcal{B})$ respectively, it holds for ($\nu$--almost) all $y\in Y$ that $\mu(Q_y)=1$, then
\begin{equation}
    \mu\{x:\nu(Q_x)=1\}=1.
\end{equation}
\end{lemma}
\begin{proof}
Tonelli on $\mathds{1}_Q(x,y)$ defined on $(X\times Y,\mathcal{A}\otimes\mathcal{B},\mu\otimes\nu)$ yields
$$
\mu\otimes\nu\,(Q)=\int_X\nu(Q_x)\,\mu(dx)=\int_Y\mu(Q_y)\,\nu(dy)
$$
since $Q\in\mathcal{A}\otimes\mathcal{B}$. The conclusion follows as $\mu(Q_y)=1$ on a subset of $Y$ of full $\nu$ measure, and noting $x\mapsto \nu(Q_x)$ is measurable.
\end{proof}

In particular, it allows us to convert statements that hold for fixed hypotheses $P_0$--almost surely into statements that hold $P_0$--almost surely simultaneously for $\Pi_0$--almost all $P\in\Theta$ \cite{shalizi2009}.

\section{Sieve decomposition I: the good sets}\label{sec:numerator-on-sieves}

This section develops the numerator on the sets of the hypothesis space with favourable behaviour. The sieve $G_t$ is fixed throughout.
\subsection{Elementary properties}
We make the following observations, which it turns out make the analysis of the Gibbs measure more straightforward.
We suppose that $M_0>\mathfrak{h}^\star$ is a fixed constant throughout. 

We recall the following elementary lemma.
\begin{lemma}
Let $(X,\mathcal{B},\mu)$ be a probability space and suppose $f:X\to[-\infty,+\infty]$ is $\mathcal{B}$--measurable, then $\forall A\in\mathcal{B}$ and $c>0$
\begin{equation}\label{eqn:lemma-L1-L-infty-bound}
    \int_Ae^{c f}\dd\mu\le e^{c\esssup_A f}\mu(A),\qquad\int_Ae^{-c f}\dd\mu\le e^{-c\essinf_A f}\mu(A)
\end{equation}
\end{lemma}
This lemma can be proved in many different ways, one of which we give below for completeness.
\begin{proof}
Since $(X,\mu)$ is a finite measure space we have the continuous embedding $L^\infty(X)\hookrightarrow L^p(X)$ for all $p\in[1,+\infty]$ which is naturally proved via H\"older's inequality
\begin{equation*}
\int_A1\cdot\left|e^{\pm c f}\right|d\mu\le ||1||_{L^1(A,\mu)}||e^{\pm cf}||_{L^\infty(A,\mu)}=\mu(A)\cdot\esssup_{x\in A}|e^{\pm c f(x)}|\,.
\end{equation*}
By monotonicity of $x\mapsto e^x$, and recalling $\esssup_A-cf=-c\essinf_A f$, the result follows.
\end{proof}

Since we use the layer cake identity throughout, we give an elementary proof of this below \cite{lieb1997}. 
\begin{lemma}[Layer Cake]
Let $\phi:[0,\infty)\to[0,\infty)$ be $C^1([0,\infty))$, increasing, and $\phi(0)=0$. Suppose $(X,\mu)$ is $\sigma$--finite and $f\in\mathcal L(X)$. Then \begin{equation}
\int_X \phi(|f|)\dd\mu=\int_0^\infty\phi'(t)\mu(|f|>t)\dd t\,.
\end{equation}
\end{lemma}
\begin{proof}
$\phi(|f|)=\phi(|f|)-\phi(0)$, so by the fundamental theorem of calculus,
\begin{gather*}
\int_X\phi(|f|)\dd\mu=\int_X\int_0^{|f|}\phi'(t)\dd t\dd\mu=\int_X\int_0^\infty\phi'(t)\mathds{1}_{|f|>t}\dd t\dd\mu\\
=\int_0^\infty\phi'(t)\int_X\mathds{1}_{|f|>t}\dd\mu\dd t=\int_0^\infty\phi'(t)\mu(|f|>t)\dd t\,.
\end{gather*}
where we have used Tonelli, as the monotonicity of $\phi$ implies $\phi'(t)\ge 0$.
\end{proof}

\subsection{A first sieve estimate}
Recall that $-s_t^{-1}\log R_t=\Delta_t$, and define the one--sided random deviation from $h\wedge M_0$ at time $t$ by
\begin{equation}
U_t(P)= (h(P)\wedge M_0)-\Delta_t(X_1^t;P).
\end{equation}
It is hence immediate
\begin{equation*}
R_t=e^{s_tU_t-s_t(h\wedge M_0)},
\end{equation*}
and
\begin{equation}
    \left\{P:h\wedge M_0>u +\Delta_t(X_1^t;P)\right\}=\{P:U_t(P)>u\}.
\end{equation}

\begin{proposition}\label{prop:first-sieve-estimate}
Let $G_t\in\mathcal{T}$ and $M_0>\mathfrak{h}^\star$ be given, and suppose that there exists $p_1\in(1,+\infty]$, $\alpha_1\in(0,\infty)$, and some $\lambda_t\ge 0$ such that $P_0$--\as 
\begin{equation}\label{eqn:asp-weak-uniformity}
    \Pi_0(G_t\cap\{U_t>u\})\le \alpha_1e^{-p_1s_t u}\,,\qquad \forall u\ge \lambda_t
\end{equation}
then for all $A\in\mathcal{T}$ it holds $P_0$--almost surely that
\begin{equation}
    Z_t(A\cap G_t)\le e^{-s_t(\mathfrak{h}(A)\wedge M_0)}\left(e^{s_t\lambda_t}+\frac{\alpha_1}{p_1-1}\right).
\end{equation}
\end{proposition}

\begin{proof}
Since
\begin{align*}
Z_t(A\cap G_t)&=\int_{A\cap G_t}R_t\dd\Pi_0\\
&=\int_{A\cap G_t}e^{s_tU_t}e^{-s_t(h\wedge M_0)}\dd\Pi_0
\end{align*}
Noting $e^{-s_t(h\wedge M_0)}\in L^{\infty}$ and $\mathfrak{h}(A\cap G_t)\ge \mathfrak{h}(A)$, by H\"older's inequality \eqref{eqn:lemma-L1-L-infty-bound}
\begin{equation}\label{eqn:inequality-A-G_t}
\le e^{-s_t(\mathfrak{h}(A)\wedge M_0)}\int_{G_t}e^{s_t U_t}\dd\Pi_0
\end{equation}
The integral can now be bounded by the layer cake representation
\begin{align*}
    \int_{G_t}e^{s_tU_t}\dd\Pi_0 &=\int_0^\infty\Pi_0(G_t\cap\{e^{s_tU_t}>\lambda\})\dd\lambda\\
    &=\int_0^{e^{s_t\lambda_t}}\Pi_0(G_t\cap\{e^{s_tU_t}>\lambda\})\dd\lambda+\int_{e^{s_t\lambda_t}}^\infty\Pi_0(G_t\cap\{e^{s_tU_t}>\lambda\})\dd\lambda
\end{align*}
The first term is controlled via 
\begin{equation*}
    \int_0^{e^{s_t\lambda_t}}\Pi_0(G_t\cap\{e^{s_tU_t}>\lambda\})\dd\lambda\le \int_0^{e^{s_t\lambda_t}}1\dd\lambda=e^{s_t\lambda_t}
\end{equation*}
and the second via assumption \labelcref{eqn:asp-weak-uniformity}, which is zero for $p_1=+\infty$, yielding for $p_1\in(1,\infty)$
\begin{align*}
\int_{e^{s_t\lambda_t}}^\infty\Pi_0(G_t\cap\{e^{s_tU_t}>\lambda\})\dd\lambda &\le \int_{e^{s_t\lambda_t}}^\infty\alpha_1e^{-p_1\log\lambda}\dd\lambda\\
&=\int_{e^{s_t\lambda_t}}^\infty\alpha_1\lambda^{-p_1}\dd\lambda\\
&=\frac{\alpha_1}{p_1-1}e^{-(p_1-1)s_t\lambda_t}
\end{align*}
where $p_1>1$ by assumption yields the integrability. Combining both bounds leads to
\begin{align*}
\int_{G_t}e^{s_tU_t}\dd\Pi_0&\le e^{s_t\lambda_t}+e^{s_t\lambda_t}\frac{\alpha_1}{p_1-1}e^{-p_1s_t\lambda_t}\\
&=e^{s_t\lambda_t}+\frac{\alpha_1}{p_1-1}e^{-(p_1-1)s_t\lambda_t}\le e^{s_t\lambda_t}+\frac{\alpha_1}{p_1-1}
\end{align*}
having used that $(p_1-1)>0$ and $s_t\lambda_t\ge 0$. Substituting this estimate into inequality \labelcref{eqn:inequality-A-G_t} completes the proof.
\end{proof}
\subsection{A shell decomposition}
Suppose we are provided a sequence $\epsilon_t>0$ for which $\epsilon_t\to0$ as $t\to\infty$ satisfying also that $s_t\epsilon_t\to\infty$, then we may decompose the hypothesis space into shells of the ``neighbourhoods'' of $h(P)$. In particular, we may define for each $t$ and integer $j\in\N$ the shells
\begin{equation}\label{eqn:shells-def}
S_t^{(j)}:=\{P\in\Theta:j\epsilon_t+\mathfrak{h}^\star\le h(P)\wedge M_0< (j+1)\epsilon_t+\mathfrak{h}^\star\}\,.
\end{equation}
Since $\epsilon_t\to0$ we have that $\epsilon_t\le M_0-\mathfrak{h}^\star$ for sufficiently large $t$ so that these shells are not vacuously empty, and we may define the largest integer $j\in\N$ for which the shells are non--empty as 
\begin{equation}\label{eqn:j_max-def}
    j^\dagger(t):=\left\lfloor\frac{M_0-\mathfrak{h}^\star}{\epsilon_t}\right\rfloor \vee 1.
\end{equation}
It is immediate that $j^\dagger(t)\to\infty$ as $t\to\infty$. 
\begin{lemma}\label{lem:shells-partition}
At each time $t$ sufficiently large, the shells and $$N_{\epsilon_t}:=\{P\in\Theta: h(P)<\mathfrak{h}^\star+\epsilon_t\}$$ yield a disjoint partition of the hypothesis space, that is 
\begin{equation}
    \Theta=N_{\epsilon_t}\;\dot{\cup}\;\dot{\bigcup}_{j=1}^{j^\dagger(t)}\;S_t^{(j)}.
\end{equation}
\end{lemma}
\begin{proof}
For any $j\ge j^\dagger(t)+1$, the shell $S_t^{(j)}=\emptyset$. Indeed, $$h\wedge M_0\le M_0<(j^\dagger(t)+1)\epsilon_t+\mathfrak{h}^\star\le j\epsilon_t+\mathfrak{h}^\star$$ for any such $j$. Further, the set $\{h\ge M_0\}\subseteq S_t^{j^\dagger(t)}$ since $h\wedge M_0=M_0$ on this set and $j^\dagger(t)\epsilon_t+\mathfrak{h}^\star\le M_0< (j^\dagger(t)+1)\epsilon_t+\mathfrak{h}^\star$, with the equality below attained on $\{h\ge M_0\}$ if and only if $(M_0-\mathfrak{h}^\star)\epsilon_t^{-1}\in\N$. That the individual shells are disjoint follows by definition, and hence 
\begin{equation*}
    \Theta\setminus \bigcup_{j=1}^{j^\dagger(t)}S_t^{(j)}=\{P\in\Theta:h(P)< \mathfrak{h}^\star+\epsilon_t\}=:N_{\epsilon_t}
\end{equation*}
as desired.
\end{proof}
\begin{remark}
Note that $S_t^{(j^\dagger(t))}\supseteq \{P:h(P)=\infty\}=:I$, and hence controlling the behaviour of the Gibbs posterior on this shell controls the behaviour on the hypotheses on which $h$ is ill--behaved, which allows for $\Pi_0(I)\in[0,1)$.
\end{remark}

The truncation that we have introduced and make use of throughout is merely a technical tool for controlling the behaviour of $R_t$ whenever this diverges to $+\infty$, but is invisible to the conclusion.

\begin{lemma}\label{lem:shell-truncation-invisibility}
For any $j\in[1;j^\dagger(t)]$ and $t$ sufficiently large, it holds 
\begin{equation*}
    \{h\wedge M_0\ge j\epsilon_t+\mathfrak{h}^\star\}=\{h\ge j\epsilon_t+\mathfrak{h}^\star\}\;\mathrm{and}\;\{h\wedge M_0<j\epsilon_t+\mathfrak{h}^\star\}=\{h<j\epsilon_t+\mathfrak{h}^\star\}.
\end{equation*}
\end{lemma}
\begin{proof}
For the first equality, suppose $h\ge j\epsilon_t+\mathfrak{h}^\star$ and it therefore holds that $h\wedge M_0\ge (j\epsilon_t+\mathfrak{h}^\star)\wedge M_0= j\epsilon_t+\mathfrak{h}^\star$ since $M_0\ge j^\dagger(t)\epsilon_t+\mathfrak{h}^\star$.

For the second, likewise by definition $h\wedge M_0<j\epsilon_t+\mathfrak{h}^\star\le j^\dagger(t)\epsilon_t+\mathfrak{h}^\star\le M_0$ and hence $h\wedge M_0=h<M_0$.

The reverse are immediate.
\end{proof}
\subsection{A refined sieve estimate}
We revisit the analysis of \cref{prop:first-sieve-estimate} accounting now for deviations of $U_t:=(h\wedge M_0)-\Delta_t$ over the shells on the sieve. %
The proposed condition is a Lorentz space generalisation of this where the behaviour is controlled on each shell individually and is analogous to those of Ghosal and van der Vaart \cite{ghosal2017}. The proof proceeds similar to the Marcinkiewicz interpolation theorem (see \cite{calderon1966} and \cite{grafakos2014}).
\begin{proposition}\label{prop:refined-sieve-estimate}
Let $(G_t)_{t}$, $\epsilon_t\to0$, and $M_0>\mathfrak{h}^\star$ be given, and assume that there exists constants $\alpha_1\in(0,+\infty)$, $p_1\in(1,+\infty]$, and $c_1\in(0,1)$  independent of $j$ and $t$, such that $P_0$--almost surely for all sufficiently large $t$ and all $j\in[1;j^\dagger(t)]$, it holds that
\begin{equation}\label{eqn:asp-weak-sieve-deviation-bound}
    \Pi_0(S_t^{(j)}\cap G_t\cap \{U_t>u\})\le \alpha_1e^{-p_1s_tu},\quad \forall u\ge c_1j\epsilon_t.
\end{equation}
Then it eventually holds $P_0$--\as, uniformly over $j\in[1;j^\dagger(t)]$, that
\begin{equation}
Z_t(\{h\ge\mathfrak{h}^\star+j\epsilon_t\}\cap G_t)\lesssim e^{-s_t\mathfrak{h}^\star} e^{-s_tj\epsilon_t(1-c_1)}
\end{equation}
where the unspecified constant depends only on $\alpha_1$ and $p_1$.
\end{proposition}
\begin{proof}
As $\epsilon_t\to 0$, we may consider \wolog $t$ sufficiently large such that $$M_0\ge \mathfrak{h}^\star+\epsilon_t.$$

Let $k\in [1;j^\dagger(t)]$ and consider $S_t^{(k)}\cap G_t$. On this set, by \cref{lem:shell-truncation-invisibility}, $h\ge \mathfrak{h}^\star+k\epsilon_t$, and $R_t=e^{-s_t(h\wedge M_0)}e^{s_tU_t}$. 

Set $\lambda_t^{(k)}=c_1k\epsilon_t$, then by layer--cake and \cref{prop:first-sieve-estimate} that
\begin{equation*}
\begin{aligned}
\int_{S_t^{(k)}\cap G_t}e^{s_tU_t}\dd\Pi_0 &= \int_0^\infty\Pi_0(S_t^{(k)}\cap G_t\cap\{e^{s_tU_t}>\lambda\})\dd\lambda\\
&\le e^{s_t\lambda_t^{(k)}}\Pi_0(S_t^{(k)})+\int_{e^{s_t\lambda_t^{(k)}}}^\infty\alpha_1\lambda^{-p_1}\dd\lambda
\end{aligned}
\end{equation*}
For $p_1=+\infty$, the set $S_t^{(k)}\cap G_t\cap\{U_t>u\}$ is $\Pi_0$--null on $u\ge\lambda_t^{(k)}$ , that is $$S_t^{(k)}\cap G_t\cap\{e^{s_tU_t}>\lambda\}$$ is null on $\lambda\ge e^{s_t\lambda_t^{(k)}}$ making the second integral zero.

For $p_1\in(1,\infty)$ the integrability of the second integral follows by $p_1>1$, hence
\begin{equation*}
\int_{e^{s_t\lambda_t^{(k)}}}^\infty\alpha_1\lambda^{-p_1}\dd\lambda=\frac{\alpha_1}{p_1-1}e^{-(p_1-1)s_t\lambda_t^{(k)}}\le \frac{\alpha_1}{p_1-1}e^{s_t\lambda_t^{(k)}}\,.
\end{equation*}
We understand the above constant as zero for $p_1=\infty$.
Substituting this into the previous yields the estimate
\begin{equation*}
\int_{S_t^{(k)}\cap G_t}e^{s_tU_t}\dd\Pi_0 \le \left(1+\frac{\alpha_1}{p_1-1}\right)e^{s_t\lambda_t^{(k)}}.
\end{equation*}
Taking $C=(1+\alpha_1/(p_1-1))$, we thus get that $P_0$--\as
\begin{align*}
Z_t(S_t^{(k)}\cap G_t) &\le C e^{-s_t(\mathfrak{h}(S_t^{(k)})\wedge M_0)}e^{s_tc_1k\epsilon_t}\\
&\le Ce^{-s_t(\mathfrak{h}^\star+k\epsilon_t)}e^{s_tc_1k\epsilon_t}\\
&=Ce^{-s_t\mathfrak{h}^\star}e^{-s_tk\epsilon_t(1-c_1)}
\end{align*}
by definition of $S_t^{(k)}$. Now, as $\forall j_0\in[1;j^\dagger(t)]$
\begin{equation*}
\bigcup_{k=j_0}^{j^\dagger(t)}(S_t^{(k)}\cap G_t)=\{h\ge \mathfrak{h}^\star+j_0\epsilon_t\}\cap G_t,
\end{equation*}
using the disjointness of the shells, \cref{lem:shells-partition}, and the previous estimate, it follows $P_0$--\as that
\begin{align*}
Z_t\left(\{h\ge \mathfrak{h}^\star+j_0\epsilon_t\}\cap G_t\right)&=\sum_{k=j_0}^{j^\dagger(t)}Z_t\left(S_t^{(k)}\cap G_t\right)\\
&\le \sum_{k=j_0}^{j^\dagger(t)}Ce^{-s_t\mathfrak{h}^\star}e^{-s_tk\epsilon_t(1-c_1)}\\
&=Ce^{-s_t\mathfrak{h}^\star}\sum_{k=j_0}^{j^\dagger(t)}e^{-s_tk\epsilon_t(1-c_1)}
\end{align*}
Noting that the sum is geometric as $e^{-s_t\epsilon_t(1-c_1)}< 1$, by elementary analysis
\begin{align*}
&=Ce^{-s_t\mathfrak{h}^\star}\left(e^{-s_tj_0\epsilon_t(1-c_1)}\left(\frac{1-\left(e^{-s_t\epsilon_t(1-c_1)}\right)^{j^\dagger(t)-j_0+1}}{1-e^{-s_t\epsilon_t(1-c_1)}}\right)\right)\\
&\le\frac{C}{1-e^{-s_t\epsilon_t(1-c_1)}}e^{-s_t\mathfrak{h}^\star}e^{-s_tj_0\epsilon_t(1-c_1)}.
\end{align*}
Further, as $s_t\epsilon_t\to\infty$ we have that $(1-e^{-s_t\epsilon_t(1-c_1)})\to 1$, and hence 
$$C(1-e^{-s_t\epsilon_t(1-c_1)})^{-1}\to C$$
as $t\to\infty$. Noting that $s_t\epsilon_t$ is not assumed to be monotone, but by the existence of the limit, for all $t$ sufficiently large, it must hold that $s_t\epsilon_t(1-c_1)\ge \log 2$, and hence $C(1-e^{-s_t\epsilon_t(1-c_1)})^{-1}\le2C$ eventually. We thus conclude.
\end{proof}

\section{Sieve decomposition II: the bad sets}\label{sec:numerator-off-sieves}

The aim of this section is establishing sufficient control of the numerator of the Gibbs posterior off the sieve where the interaction of the observables and the hypotheses through the connection $\Delta_t$ may be arbitrarily misbehaved by considering sets of exponentially small prior mass. 
\subsection{A weak type Hardy--Littlewood inequality}
We first record the following weak type Hardy--Littlewood inequality proved via decreasing rearrangements which will be instrumental in the proof of \cref{prop:numerator-bad-sets}. And in particular, this proves the continuous embedding for $p>1$
\begin{equation}
    L^{p,\infty}\hookrightarrow L^1
\end{equation}
on a finite measure space (see also Bennett and Sharpley \cite{bennett1988}).
\begin{proposition}\label{prop:hardy-littlewood}
Let $(X,\mu)$ be a probability space and $p\in(1,\infty]$. Then for any measurable set $B\subseteq X$ with $\mu (B)>0$ and measurable $f :X\to(-\infty,\infty]$ it holds
\begin{equation}\label{eqn:hardy-littlewood}
    \int_B|f|\dd\mu\le \frac{p}{p-1}||f||_{L^{p,\infty}(\mu)}\mu(B)^{\frac{p-1}{p}}\,.
\end{equation}
\end{proposition}
\begin{proof}
Consider first the case where $p<\infty$. We may assume without loss of generality that $f$ is nonnegative as otherwise we may take $|f|$ and use $\int f\le\int |f|$. Applying the layer cake representation to the left hand side yields
\begin{align*}
\int_B f\dd\mu &= \int_0^\infty\mu(\{x\in B: f(x)>\lambda\})\dd\lambda\\
&= \int_0^\infty \mu(B\cap\{f>\lambda\})\dd\lambda\\
&\le \int_0^\infty \left(\mu(B)\wedge\mu(f>\lambda)\right)\dd\lambda
\end{align*}
We may now relate this to the Lorentz quasi--norm of $f$.

Denote by $\mu_f(\lambda):=\mu(f>\lambda)$ the distribution function of $f$, and by $f^\star(t):=\inf\{\lambda\ge0:\mu_{f}(\lambda)\le t\}$ the decreasing rearrangements of $f$. It holds that
\begin{equation}\label{eqn:shalizinewproofclaim}
    \int_0^\infty \left(\mu(B)\wedge\mu(f>\lambda)\right)\dd\lambda= \int_0^{\mu(B)}f^\star(s)\,ds\,.
\end{equation}

Indeed, one may obtain this by first considering
\begin{equation}\label{eqn:proofnewshaliziiff}
\mu_{f}(\lambda)> t \iff \lambda<f^\star(t).
\end{equation}
For this, it will be useful to recall that $f^\star(\mu_{f}(\lambda))\le \lambda$ provided $\mu_{f}(\lambda)<\infty$ and that $\mu_{f}(f^\star(t))\le t$ provided $f^\star(t)<\infty$ (see Proposition 1.7 in Chapter 2 of \cite{bennett1988}).

Suppose first that $\mu_{f}(\lambda)> t$. If $f^\star(t)=\infty$ the result holds trivially, so suppose $f^\star(t)<\infty$ and for a contrapositive that $f^\star(t)\le \lambda$. Now as $\mu_{f}$ is decreasing, $\mu_{f}(\lambda)\le\mu_{f}(f^\star(t))\le t$. The forward direction of \eqref{eqn:proofnewshaliziiff} hence holds.

Suppose now that $f^\star(t)> \lambda$, then suppose for a contrapositive that $\mu_{f}(\lambda)\le t$. Then as $f^\star$ is decreasing, $f^\star(t)\le f^\star(\mu_{f}(\lambda))\le \lambda$. Hence, the reverse direction and thus \eqref{eqn:proofnewshaliziiff} holds.

Now, for $\mu(B)=b\in[0,1]$, consider $E\subseteq[0,\infty)^2$,
\begin{align*}
E:&=\{(s,\lambda):s<\min\{b,\mu_{f}(\lambda)\}\}\\
&=\{(s,\lambda):s<b\,,\, s <\mu_{f}(\lambda)\}\\
&=\{(s,\lambda):s<b\,,\, \lambda <f^\star(s)\}
\end{align*}
by equivalence \labelcref{eqn:proofnewshaliziiff}. It is straightforward to see that this is measurable as the countable intersection of measurable sets.

Integrating $E$ with respect to $s$ yields $\min\{b,\mu_{f}(\lambda)\}$ as a function of $\lambda$, and integrating again with respect to $\lambda$ yields the right hand side of \labelcref{eqn:shalizinewproofclaim}. Conversely, integrating first with respect to $\lambda$, one gets $f^\star(s)\mathds{1}_{s<b}(s)$.

Hence, as $\mathds{1}_E\ge 0$, Tonelli yields
\begin{align*}
\int_0^\infty \left(\mu(B)\wedge\mu(f>\lambda)\right)\dd\lambda
=&\int_0^\infty\int_0^\infty\mathds{1}_E(s,\lambda)\,ds\,d\lambda\\
=&\int_0^\infty\int_0^\infty\mathds{1}_E(s,\lambda)\,d\lambda\,ds\\
=&\int_0^\infty f^\star(s)\mathds{1}_{s<b}(s)\,ds=\int_0^{b}f^\star(s)\,ds
\end{align*}
establishing \eqref{eqn:shalizinewproofclaim}.

To conclude, by definition of the Lorentz quasi--norm, for all $\lambda>0$
\begin{equation*}
    \mu_{f}(\lambda)\le \left(\frac{||f||_{L^{p,\infty}}}{\lambda}\right)^p\,.
\end{equation*}
Taking $\lambda'=||f||_{L^{p,\infty}}s^{-1/p}$, it follows $\mu_{f}(\lambda')\le s \implies f^\star(s)\le \lambda'$ by definition of $s\mapsto f^\star(s)$. Now, let $q=p/(p-1)$ be conjugate to $p$, then
\begin{equation*}
    \int_0^{b}f^\star(s)\,ds\le ||f||_{L^{p,\infty}}\int_0^{b}s^{-1/p}ds=||f||_{L^{p,\infty}}\cfrac{b^{1-\frac1p}}{1-\frac1p}=q\cdot||f||_{L^{p,\infty}}\cdot b^{1/q}
\end{equation*}
as desired. 

The case where $p=+\infty$ holds trivially since $\mu$ is a probability measure, so that we may bound the integral by the essential supremum over $B$ times the size of the set, and we take $q:=\frac{p}{p-1}=1$, completing the proof.
\end{proof}
In the statement of the proposition, we have not assumed that $f\in L^{p,\infty}$ which is not strictly necessary but otherwise makes the right hand side vacuously $+\infty$ whenever $f\notin L^{p,\infty}$ for some $p>1$.

\subsection{Controlling behaviour on exponentially thin sets}
\begin{proposition}\label{prop:numerator-bad-sets}
Under only \cref{asp:lorentz} pick deterministic $B_t\in\mathcal{T}$ for which $\Pi_0(B_t)\le\alpha e^{-s_t\gamma}$ for some $\alpha,\gamma>0$. Then $\forall\epsilon>0$, it holds $P_0$--\as that
\begin{equation}
   Z_t(B_t):=\int_{B_t}R_t\dd\Pi_0 \le\exp(-s_t(\gamma/q-\kappa-\epsilon) )
\end{equation}
for all but finitely many $t$.
\end{proposition}
\begin{proof}
By standing assumption, the map $(X_1^t,P)\mapsto \exp(-s_t\Delta_t(X_1^t; P))\mathds{1}_{B_t}$ is $\sigma(X_1^t)\otimes\mathcal{T}$ measurable, and further nonnegative, taking values in the extended reals. Hence, applying Tonelli yields  
\begin{align*}
    \E_{P_0^{(t)}}Z_t(B_t)&= \int_{\Xi^{(t)}}\int_{B_t}R_t(x_1^t;P)\,\Pi_0(dP)\,P_0^{(t)}(dx_1^t)\\
    &=\int_{B_t}\int_{\Xi^{(t)}}R_t(x_1^t; P)\,P_0^{(t)}(dx_1^t)\,\Pi_0(dP)\\
    &=\int_{B_t} c_t(P)\,\Pi_0(dP)\,.
\end{align*}
The weak type Hardy--Littlewood inequality \labelcref{eqn:hardy-littlewood} now yields that
\begin{equation}
\E_{P_0^{(t)}}Z_t(B_t)\le q||c_t||_{L^{p,\infty}}\,\Pi_0(B_t)^{1/q}\le q||c_t||_{L^{p,\infty}}\,\alpha^{1/q} e^{-s_t\gamma/q}\,.
\end{equation}

Moreover, by \cref{asp:lorentz}, for any $\epsilon>0$ and all $t$ sufficiently large, we have that $||c_t||_{L^{p,\infty}}\le \exp((\epsilon+\kappa) s_t)$, and thereby eventually
\begin{equation}
\E_{P_0^{(t)}}Z_t(B_t)\le q\alpha^{1/q}\exp(-s_t(\gamma/q-\epsilon-\kappa))\,.
\end{equation}

Applying Markov's inequality,
\begin{align*}
P_0^{(t)}\big\{x_1^t :\; &Z_t(B_t)>\exp(-s_t(\gamma/q-\kappa-2\epsilon))\big\} \\
&\le \exp(s_t(\gamma/q-\kappa-2\epsilon))\E_{P_0^{(t)}}[Z_t(B_t)]\\
&\le \exp(s_t(\gamma/q-\kappa-2\epsilon))\cdot q\alpha^{1/q}\exp(-s_t(\gamma/q-\epsilon-\kappa))\\
&=q\alpha^{1/q}\exp(-s_t\epsilon)
\end{align*}
for all $\epsilon>0$ and all $t$ sufficiently large. Finally, as the sequence $\exp(-\delta s_t)$ is assumed summable for all $\delta>0$, the first Borel--Cantelli lemma implies that the above bound may occur only finitely often, and the desired conclusion thus holds.
\end{proof}

From the proposition, it is immediate that we get the following inequality.
\begin{corollary}
It holds $P_0$ almost surely that
\begin{equation}
    \limsup_{t\to\infty} s_t^{-1}\log Z_t(B_t)\le \kappa-\gamma/q\,.
\end{equation}
\end{corollary}

\part{Analysis of the Gibbs posterior normalisation}\label{part:denominator-analysis}

\section{Denominator estimates}\label{sec:denominator-estimates}

This section is dedicated to providing more precise control of the denominator of the Gibbs posterior. Whereas the numerator has required upper bounds on the unnormalised Gibbs posterior, the denominator requires lower bound control, so that we may ensure that $Z_t(\Theta)$ increases sufficiently quickly. We will in particular 
study the denominator via the random quantity
\begin{equation}
    \eta_t:=-\frac1{s_t}\log \left(e^{^{s_t\mathfrak{h}^\star}}{Z_t(\Theta)}\right)=-\mathfrak{h}^\star-\Lambda_t,\quad\mathrm{where}\quad\Lambda_t=\frac1{s_t}\log Z_t(\Theta).
\end{equation}
\subsection{The denominator is well defined}
An immediate, but useful observation is that the denominator of the Gibbs posterior is almost surely finite for all but finitely many $t$.
\begin{lemma}\label{lem:finite-Z_t-Theta}
By \cref{asp:lorentz}, $Z_t(\Theta)<\infty$ \as for all but finitely many $t$.
\end{lemma}
\begin{proof}
Applying the Hardy--Littlewood inequality \labelcref{eqn:hardy-littlewood} with $B=\Theta$, and Tonelli yields that
\begin{align*}
q\cdot||c_t||_{L^{p,\infty}}\ge \int_\Theta c_t\dd\Pi_0=\int_\Theta\int_{\Xi^{(t)}}R_t\dd P_0^{(t)}\dd\Pi_0=\int_{\Xi^{(t)}}Z_t(\Theta)\dd P_0^{(t)}.
\end{align*}
By \cref{asp:lorentz}, the left hand side of the above is finite for all but finitely many $t$. And therefore, the set of $\omega\in\Omega$ on which $Z_t(\Theta)(\omega)=\infty$ has $P_0^{(t)}$ measure zero. Thus, $Z_t<+\infty$ $P_0$ almost surely.
\end{proof}
Thus, and recalling that $Z_t(\Theta)>0$ $P_0$--almost surely by standing assumption, taking logarithms yields the following result which ensures that we may study the behaviour of $\eta_t$. 
\begin{corollary}
It holds that $\eta_t:=-\mathfrak{h}^\star-\Lambda_t>-\infty$ \as for all but finitely many $t$. Moreover by standing assumption, $\eta_t<+\infty$ \as for all $t$. And thus, for all $t$ sufficiently large, $P_0$--almost surely
\begin{equation}
    \eta_t\in(-\infty,+\infty).
\end{equation}
\end{corollary}

An immediate estimate on the behaviour of $\Lambda_t$ is the following.
\begin{lemma}
It holds for all $\epsilon>0$ and $t$ sufficiently large that
\begin{equation}
    P_0(\Lambda_t>\kappa +2\epsilon)\le q e^{-s_t\epsilon}
\end{equation}
and therefore
\begin{equation}
    \sum_tP_0(\Lambda_t>\kappa +2\epsilon)<+\infty.
\end{equation}
Moreover, as a corollary, it holds $P_0$--almost surely that
\begin{equation}
    \limsup_{t\to\infty}\Lambda_t\le \kappa
\end{equation}
\end{lemma}
\begin{proof}
By the previous argument in \cref{lem:finite-Z_t-Theta}, we have that 
\begin{equation*}
    \int_{\Xi^{(t)}}Z_t(\Theta)\dd P_0^{(t)}\le q||c_t||_{L^{p,\infty}}
\end{equation*}
from which \cref{asp:lorentz} yields for all $\epsilon>0$ that eventually
\begin{equation*}
    \int_{\Xi^{(t)}}Z_t(\Theta)\dd P_0^{(t)}\le qe^{s_t(\kappa+\epsilon)}
\end{equation*}

Thus, by Markov's inequality, we obtain that
\begin{align*}
    P_0\left(Z_t(\Theta)>e^{s_t(\kappa+2\epsilon)}\right)
    &\le e^{-s_t(\kappa+2\epsilon)}\int_{\Xi^{(t)}}Z_t(\Theta)\dd P_0^{(t)}\\
    &\le e^{-s_t(\kappa+2\epsilon)}qe^{s_t(\kappa+\epsilon)}\\
    &=qe^{-s_t\epsilon}
\end{align*}
which is summable by assumption on $s_t$, so the first Borel--Cantelli lemma yields that this may occur only finitely often. This therefore immediately yields that
\begin{equation}\label{eqn:limsup_Lambda_t_bound}
    \limsup_{t\to\infty}\Lambda_t\le \kappa + 2\epsilon.
\end{equation}

Taking any countable sequence of $\epsilon_m\downarrow0$, and considering the full $P_0$ measure sets $\Omega_m\subseteq\Omega$ on which \labelcref{eqn:limsup_Lambda_t_bound} holds for $\epsilon_m$, we may consider $$\Omega_\infty:=\bigcap_{m\in\N}\Omega_m.$$ 

To conclude, this set is a countable intersection of full measure sets and hence has full $P_0$--measure, which completes the proof.
\end{proof}
This yields the following direct corollary.
\begin{corollary}
$P_0$--almost surely
\begin{equation*}
\liminf_{t\to+\infty}\eta_t\ge -\mathfrak{h}^\star-\kappa
\end{equation*}
\end{corollary}

\subsection{Qualitative denominator estimates}
The following theorem relates a rather weak assumption on the relationship between the observables and the difference function to a qualitative denominator estimate in lieu of establishing rates at which the conclusion is achieved. In particular, this result will be instrumental in the proofs of \cref{thm:consistency,thm:contraction_bad-sets} which do not require rates. 
\begin{theorem}\label{thm:qualitative-denominator-floor}
Suppose that $\Pi_0$ almost everywhere it holds that 
\begin{equation}\label{eqn:limsup-convergence-Delta_t}
    P_0\left(\limsup_{t\to\infty}\Delta_t(X_1^t(\cdot);P)\le h(P) \right)=1
\end{equation}
then it holds $P_0$--almost surely that 
\begin{equation}
    \liminf_{t\to\infty}\Lambda_t\ge -\mathfrak{h}^\star.
\end{equation}
\end{theorem}
\begin{proof}
The set $B=\{(\omega,P):\limsup_t\Delta_t>h\}$ is jointly measurable, and by assumption has $P_0$--measure zero $\Pi_0$--almost everywhere. By the quantifier reversal \cref{lem:quantifier-reversal}, we may take the $P_0$ full measure set $C\subset \Omega$ on which assumption \labelcref{eqn:limsup-convergence-Delta_t} does hold, and for any such $\omega\in C$ we notice that the section $B_\omega$ of $B$ has $\Pi_0(B_\omega)=0$. 

Fix such an $\omega\in C$ and let $\delta>0$. We may then set $N_\delta:=\{h < \mathfrak{h}^\star+\delta\}$ with $m_\delta:=\Pi_0(N_\delta)>0$ by definition of the essential infimum. 

Define also the functions
\begin{equation*}
    f_t:=1\wedge \left(\left(\sup_{s\ge t}\Delta_s(\omega;\cdot\,)-(\mathfrak{h}^\star+2\delta)\right)\vee 0\right), 
\end{equation*}
which is measurable, taking values in $[0,1]$, non--increasing, and $f_t\downarrow0$ \aev on $N_\delta$. Now, by Egorov's theorem $\exists E_\delta^\omega\subseteq N_\delta$ with $\Pi_0(E_\delta^\omega)>m_\delta/2$ and $f_t\downarrow0$ uniformly on $E_\delta^\omega$. Therefore, there exists a $T_0=T_0(\omega,\delta)$ such that for all $t\ge T_0$ $\sup_{E_\delta^\omega}\Delta_t\le\mathfrak{h}^\star+3\delta$.

This leads us to conclude
\begin{equation}
    Z_t(\Theta)\ge Z_t(E_\delta^\omega)\ge \frac{m_\delta}2 e^{-s_t(\mathfrak{h}^\star+3\delta)}\qquad\forall t\ge T_0(\omega,\delta).
\end{equation}
Thereby, it holds $P_0$--almost surely that
\begin{equation}\label{eqn:qualitative-denom-proof-liminf}
\liminf_{t\to+\infty}\Lambda_t\ge -\mathfrak{h}^\star-3\delta.
\end{equation}
Taking any countable sequence $\delta_k\downarrow0$, and considering the full $P_0$ measure sets $\Omega_k\subseteq\Omega$ on which \labelcref{eqn:qualitative-denom-proof-liminf} holds for $\delta_k$, we may consider $$\Omega_\infty:=\bigcap_{k\in\N}\Omega_k.$$ 

To conclude, this set is a countable intersection of full measure sets and hence has full $P_0$--measure, which completes the proof.
\end{proof}

\begin{remark}
This theorem makes the dependence of the random Egorov sets on $\omega\in\Omega$ explicit, but we notice that this cannot give us a contraction rate, which is one of the intricacies unaccounted for in \cite{shalizi2009}, here in particular Lemma 10.
\end{remark}

It is worth noting that assumption \labelcref{eqn:limsup-convergence-Delta_t} is quite general and satisfied under general conditions. 

\begin{corollary}
If for $\Pi_0$--almost every $P\in\Theta$ it holds $P_0$--almost surely that
\begin{equation}\label{eqn:a-s-pointwise-convergence}
    s_t^{-1}\log R_t(P)\to-h(P),
\end{equation}
then it holds with $P_0$ probability 1 that
$$\liminf_t\Lambda_t\ge -\mathfrak{h}^\star.$$
\end{corollary}
\begin{remark}
    Recall that the left hand side of \labelcref{eqn:a-s-pointwise-convergence} equals $-\Delta_t(X_1^t;P)$ and that the right hand side is defined as the limit in $t$ of the expectation of the right hand side, and therefore this is satisfied if some ergodic theorem or asymptotic equipartition property holds.
\end{remark}

\subsection{Small ball deviations}\label{subsec:small-ball-dev}
Before relating the denominator to the given rate $\epsilon_t$ in the next subsection, we introduce here some tools that will be instrumental in the proof of \cref{thm:rates}. In particular, the following ideas are inspired by the work of van der Vaart and van Zanten \cite{vanderVaart2008} in order to investigate the small deviations of the prior on sets close to $\mathfrak{h}^\star$. The set up is also inspired by exponential Tauberian theorems (see e.g. Kasahara and Kosugi \cite{kasahara2002}) as well as interpolation theory in particular the Peetre $K$--functional (see e.g. Calder\'on \cite{calderon1966} or Bennett and Sharpley \cite{bennett1988}).

\subsubsection*{Small deviation}
We let $g:\Theta\to(-\infty,+\infty]$ be a deterministic function or defined for a fixed $\omega\in\Omega$, where we denote $\mathfrak{g}^\star=\essinf_\Theta g>-\infty$. Define further,  for finite $c\ge \mathfrak g^\star$
\begin{equation}\label{eqn:def-Phi_g^c}
    \Phi_g^c(u):=-\log \Pi_0(g \le c+u)\qquad\forall u>0,
\end{equation}
which we may observe to be non--increasing, finite, and nonnegative for all $u>0$. If we do not specify a $c\in\R$ then we implicitly take $c=\mathfrak{g}^\star$.
Define also the free energy
\begin{equation}
    \mathcal{E}_{\Pi_0}[g](s):=-\frac1s\log\int_\Theta e^{-sg}\dd\Pi_0=-\log ||\exp(-g)||_{L^s(\Pi_0)}
\end{equation}
and  further the concentration function
\begin{equation}
    \mathcal{J}^c_g(s):=\inf_{u>0}\left\{u+\frac{\Phi^c_g(u)}{s}\right\}.
\end{equation}
\begin{lemma}
$s\mapsto\mathcal{E}_{\Pi_0}[g](s)$ is non--increasing, $\mathcal{E}_{\Pi_0}[g](s)\downarrow \mathfrak{g}^\star$, and satisfies
\begin{equation}
\mathfrak{g}^\star\le \mathcal{E}_{\Pi_0}[g](s)\le \mathfrak{g}^\star+u+s^{-1}\Phi_g(u) \qquad \forall u>0.
\end{equation}
\end{lemma}
\begin{proof}
The non--increasing part is immediate by the continuous embedding $L^{s'}\hookrightarrow L^s$ for all $\infty\ge s'\ge s> 0$ by H\"older's inequality, and in particular the (quasi--)norm relation $||e^{-g}||_{L^{s'}}\ge ||e^{-g}||_{L^{s}}$ on $(\Theta,\Pi_0)$. 
For the monotone convergence we make use of the norm relation and using \labelcref{eqn:lemma-L1-L-infty-bound} notice that $||e^{-g}||_{L^\infty}=e^{-\mathfrak{g}^\star}$.

Now, by definition of essential infimum $g\ge \mathfrak{g}^\star$ $\Pi_0$--\aev, and take $u>0$ then, for the upper bound, we notice
\begin{align*}
    \int_\Theta e^{-sg}\dd\Pi_0 &\ge \int_{\{g\le\mathfrak{g}^\star+u\}} e^{-sg}\dd\Pi_0\\
    &\ge\int_{\{g\le\mathfrak{g}^\star+u\}} e^{-s(\mathfrak{g}^\star+u)}\dd\Pi_0\\
    &=\Pi_0(g\le\mathfrak{g}^\star+u)e^{-s(\mathfrak{g}^\star+u)}.
\end{align*}
Thus, after taking logarithms, the result follows.
\end{proof}
Optimising over the choice of $u>0$ leads us to the natural connection to $\mathcal{J}_g(s)$.
\begin{corollary}
$\mathcal{E}_{\Pi_0}[g](s)\le \mathfrak{g}^\star+\mathcal{J}_g(s)$.
\end{corollary}

We may also make the following observation which makes the relationship to the learning rate $\beta_t$ more explicit.
\begin{corollary}
For all $s>0$ we have $\mathcal{E}_{\Pi_0}[g](s)\le \mathcal{E}_{\Pi_0}[g](\beta s)$ $\forall\beta\in(0,1]$.
\end{corollary}
\begin{proof}
    This follows by the monotonicity of $s\mapsto \mathcal{E}_{\Pi_0}[g](s)$, and this being non--increasing.
\end{proof}

For the next results we further define the following.
\subsubsection*{Rearrangement identities} 
Recall the increasing rearrangements of $g$ 
\begin{equation}
    g_\star(\lambda)=\inf\{v\in[-\infty,+\infty]: \Pi_0(g\le v)\ge \lambda\}\qquad\lambda\in(0,1],
\end{equation}
which is non--decreasing. Likewise, denote also its Hardy average function by 
\begin{equation}
    g_{\star\star}(v):=\frac1v\int_0^vg_\star\dd\,\mathrm{Leb}
\end{equation}
and define also the following quantity
\begin{equation}
    \overline{\Phi}^{c}_g(a):=-\log\sup\{v\in(0,1]:g_{\star\star}(v)\le c+a\}.
\end{equation}
We record some elementary observations.
\begin{lemma}\label{lem:elementary-props-increasing-rearrangements}
The following identities hold.
\begin{gather}
 g_\star(\lambda)\le y\iff\Pi_0(g\le y)\ge \lambda,\\
 \{\lambda:g_\star(\lambda)\le y\}= (0,\Pi_0(g\le y)],\\
 \forall \phi\ge 0\qquad \int_\Theta\phi\,\circ g \dd\Pi_0=\int_0^1\phi\,\circ g_\star\dd\,\mathrm{Leb},\\
 g\ge c_0\; \Pi_0\mathrm{-a.e.} \iff g_\star\ge c_0\;\mathrm{on}\;(0,1]\\
 \forall t \in (0,1] \qquad g_{\star\star}(t)\le g_\star(t)
\end{gather}
\end{lemma}
\begin{proof}
The first is immediate as $g_\star$ is non--decreasing and left continuous, and hence directly implies the second. The third is an application of the layer cake representation. The fourth follows likewise by definition, using $\Pi_0(g<c_0)=0$. The last follows likewise by the non--decreasing property of $g_\star$.
\end{proof}
\begin{lemma}\label{lem:bathtub-onesided}
Suppose $g\ge c_0:=c-b_0$ for finite $b_0\ge 0$ $\Pi_0$--\aev Then for all $E\in\mathcal{T}$ with $\Pi_0(E)>0$, 
\begin{equation}
    \fint_Eg\dd\Pi_0:=\frac1{\Pi_0(E)}\int_Eg\dd\Pi_0\ge g_{\star\star}(\Pi_0(E)).
\end{equation}
\end{lemma}
\begin{proof}
Denote $\lambda:=\Pi_0(E)$ and $\tau=g_\star(\lambda)$.
Suppose first that $\tau<\infty$, and hence the following is well defined, 
$$
g=\tau-(\tau-g)\ge \tau-(\tau-g)\vee 0,
$$
also using the almost everywhere lower boundedness. Thus, 
\begin{align*}
    \int_E g\dd\Pi_0\ge\int_E\tau-(\tau-g)\vee 0\dd\Pi_0&=\tau\lambda -\int_E(\tau-g)\vee0\dd\Pi_0\\
    &\ge\tau\lambda-\int_\Theta(\tau-g)\vee0\dd\Pi_0.
\end{align*}
Since $g\ge c_0$ \aev, we find $(\tau-g)\vee0\le (\tau-c_0)\vee 0=\tau-c_0$ \aev, which shows that the right hand side is finitely integrable. Now, by the previous lemma, we have 
\begin{equation*}
    \int_\Theta(\tau-g)\vee0\dd\Pi_0=\int_0^1(\tau-g_\star)_+\dd\,\mathrm{Leb}.
\end{equation*}
By the non--decreasing property and the definition, we have that $\forall u\ge \lambda$, $g_\star(u)\ge g_\star(\lambda)=\tau$, and hence $(\tau-g_\star(u))\vee0=0$ for such $u$. We therefore arrive at
\begin{equation*}
    \int_0^1(\tau-g_\star)_+\dd\,\mathrm{Leb}=\int_0^\lambda(\tau-g_\star)_+\dd\,\mathrm{Leb}=\lambda\tau-\lambda g_{\star\star}(\lambda),
\end{equation*}
which completes the case where $\tau<\infty$.

Suppose now that $\tau=+\infty$, then by definition, $\Pi_0(g\le x)<\lambda$ for all $x<+\infty$. 
For $n\in\N$ with $n\ge c_0\vee 1$ let $g_n:=g\wedge n\uparrow g$, then as $g_n-c_0\ge0$ \aev and we are on a finite measure space, the monotone convergence theorem implies that
$$
\fint_Eg_n\dd\Pi_0\to\fint_E g\dd\Pi_0.
$$
Now, we make use of $(g\wedge n)_\star=g_\star\wedge n$. Indeed, let $\mu\in(0,1]$ be arbitrary and suppose that for some $v<n$ we have $\Pi_0(g\le v)=\Pi_0(g_n\le v)\ge \mu$, then $(g\wedge n)_\star(\mu)\le v<n$, and if no such $v<n$ exists, then for all $v\ge n$, $\Pi_0(g_n\le v)=\Pi_0(g_n\le n)=1\ge \mu$, in which case $(g\wedge n)_\star(\mu)=n$. Hence $(g\wedge n)_\star=g_\star\wedge n$ and for $\lambda=\Pi_0(E)$ we may apply the previous result in the finite $\tau$ setting to $g_n$, so that for all $n\in\N$
\begin{equation}\label{eqn:increasing-rearrangement-bound-finite-n}
\fint_Eg_n\dd\Pi_0\ge \lambda^{-1}\int_0^\lambda (g_\star\wedge n)\dd\,\mathrm{Leb}.
\end{equation}
Moreover, making use of $g\ge c_0$ \aev iff $g_\star\ge c_0$ \aev on $(0,1]$ from the previous lemma, $g_\star\wedge n-c_0\ge 0$ and noting that $(g_\star\wedge n-c_0)\uparrow g_\star-c_0$, the monotone convergence theorem yields again that
$$
\int_0^\lambda g_\star\wedge n\dd\,\mathrm{Leb}\to \int_0^\lambda g_\star\dd\,\mathrm{Leb}.
$$
Therefore, as the limits exist, these commute with the inequality \labelcref{eqn:increasing-rearrangement-bound-finite-n}, so that we may conclude the desired result. The proof is complete.
\end{proof}

We also state the following lemma.
\begin{lemma}\label{lem:bathtub-sandwich}
Suppose $\exists b_0\ge0$ with $g\ge c-b_0>-\infty$ $\Pi_0$--\aev Then $\forall a>0$
\begin{equation}\label{eqn:bathtub-sandwich}
\Phi_g^c(2a+b_0)-\log 2\le \overline{\Phi}_g^c(a)\le\Phi_g^c(a).
\end{equation}
\end{lemma}
\begin{proof}
Turn first to the upper bound $\overline{\Phi}_g^c(a)\le\Phi_g^c(a)$. 
If $\Phi_g^c(a)=+\infty$ then the upper bound is trivial. We therefore suppose $\Phi_g^c(a)<\infty$, from which it follows that $v_0:=\Pi_0(g\le c+a)>0$. This in turn implies that $g_\star(w)\le c+a$ for all $w\in(0,v_0]$ by \cref{lem:elementary-props-increasing-rearrangements}, and also that $g_{\star\star}(v_0)\le c+a$. We therefore have that $$v_0\in\{v\in(0,1]:g_{\star\star}(v)\le a+c\}$$ from which we may conclude that
\begin{equation*}
\overline{\Phi}_g^c(a):=-\log\sup\{v\in(0,1]:g_{\star\star}(v)\le a+c\}\le -\log v_0,
\end{equation*}
proving the upper bound.

Turn now to the lower bound $\Phi_g^c(2a+b_0)-\log 2\le \overline{\Phi}_g^c(a)$. The case where $\overline{\Phi}_g^c(a)=+\infty$ is again trivial, so we may suppose $\overline{\Phi}_g^c(a)$ is finite. 

Let $v\in\{v\in(0,1]:g_{\star\star}(v)\le a+c\}$ be arbitrary, and which exists by the finiteness assumption. Then, for this $v$ it holds that $g_{\star\star}(v)\le a+c$, and hence
\begin{equation*}
    v\,g_{\star\star}(v)=\int_0^{v/2}g_\star\dd\,\mathrm{Leb}+\int_{v/2}^vg_\star\dd\,\mathrm{Leb}\le v(a+c).
\end{equation*}
Making use of the equivalence $g\ge c_0$ \aev $\iff g_\star\ge c_0$ on $(0,1]$,
we have
\begin{equation*}
    \int_0^{v/2}g_\star\dd\,\mathrm{Leb}\ge \frac v2 c_0,
\end{equation*}
and, since $g_\star$ is non--decreasing
\begin{equation*}
    \int_{v/2}^v g_\star \dd\,\mathrm{Leb}\ge \frac v2 g_\star(v/2).
\end{equation*}
Combining both estimates yields
\begin{align*}
    v(a+c)\ge v\,g_{\star\star}(v)&\ge \frac v2 c_0 + \frac v2 g_\star(v/2)\\
    \implies\quad 2a+c+b_0&\ge g_\star(v/2)
\end{align*}
from which \cref{lem:elementary-props-increasing-rearrangements} again yields that
\begin{equation*}
    \Pi_0(g\le 2a+c+b_0)\ge \frac v2.
\end{equation*}
Therefore, we have that 
\begin{equation*}
    \Phi_g^c(2a+b_0)=-\log\Pi_0(g\le c+(2a+b_0))\le -\log v+\log 2.
\end{equation*}
Recalling that $v$ was arbitrary, this yields
\begin{equation}
    \Phi_g^c(2a+b_0)-\log 2\le -\log v\qquad\forall v\in\{v\in(0,1]:g_{\star\star}(v)\le a+c\}.
\end{equation}
Taking an increasing sequence $v\uparrow$ yields the supremum and we may conclude.
\end{proof}

\subsection{Proof of \texorpdfstring{\cref{thm:optimal-choice-E_t}}{Theorem \ref*{thm:optimal-choice-E_t}} completed}

We are now in a position to prove \cref{thm:optimal-choice-E_t}, which illuminates the choice of $E_t$.
\begin{proof}[Proof of \cref{thm:optimal-choice-E_t} completed]
That  \labelcref{eqn:condition-A-thm-optimality-g_t-E_t} implies \labelcref{eqn:condition-B-thm-optimality-g_t-E_t} follows by definition of $\Phi_g^c(a)$ and the average integral.

The intermediate implication, $\mathrm{\labelcref{eqn:condition-B-thm-optimality-g_t-E_t}}{\implies}\mathrm{\labelcref{eqn:condition-C-thm-optimality-g_t-E_t}}$, follows by \cref{lem:bathtub-onesided} applied in combination with the definition of $a\mapsto\overline{\Phi}_g^c(a)$.

The final implication
$$\mathrm{\labelcref{eqn:condition-C-thm-optimality-g_t-E_t}}\overset{((2a+\bar{b})\epsilon_t,bs_t\epsilon_t+\log2)}{\implies}\mathrm{\labelcref{eqn:condition-A-thm-optimality-g_t-E_t}}$$
is proved by \cref{lem:bathtub-sandwich}, where we note the changes in the constants. 
\end{proof}

\subsection{Quantitative denominator estimates}
In this section, we provide two different routes to the desired intermediary conclusion. These are fundamentally concerned with analysing the following object, which we are led to by applying the methods established in \cref{subsec:small-ball-dev}
\begin{equation}
    -\Lambda_t(\omega)=\mathcal{E}_{\Pi_0}[\Delta_t(\omega;\cdot)](s_t), %
\end{equation}
recalling that $\eta_t:=-\mathfrak{h}^\star-\Lambda_t$.
And in particular, we establish that eventually almost surely, 
\begin{equation}
    \limsup_{t\to+\infty}\frac{\eta_t}{\epsilon_t}\le \mathfrak{L}<\infty
\end{equation}
where $\mathfrak{L}$ is independent of $\omega\in\Omega$.

\subsubsection{Endpoint control}
We first phrase a lower bound on $Z_t(\Theta)$ via the $L^{1,\infty}$ norm of $R_t$.
\begin{lemma}\label{lem:L^1-infty-bound-on-Z_Theta}
    It holds for all $u>0$ that
    \begin{equation}
        Z_t(\Theta)=||R_t||_{L^1(\Pi_0)}\ge ||R_t||_{L^{1,\infty}(\Pi_0)}\ge e^{-s_t(\mathfrak{h}^\star+u)}\Pi_0(\Delta_t<\mathfrak{h}^\star+u).
    \end{equation}
\end{lemma}
\begin{proof}
The first inequality is Markov's inequality, whereas the latter follows by taking for all $u>0$, in the definition of the $L^{1,\infty}(\Pi_0)$ quasi--norm, $\lambda_u=e^{-s_t(\mathfrak{h}^\star+u)}$ and noticing that 
\begin{equation*}
\left\{R_t>e^{-s_t(\mathfrak{h}^\star+u)}\right\}=\left\{\Delta_t< \mathfrak{h}^\star+u\right\}
\end{equation*}
yields the desired conclusion.
\end{proof}

Using this, we may now observe the following result, which will be used to establish \cref{thm:rates-endpoint-control}.
\begin{theorem}\label{thm:denominator-via-lorentz-bound}
Suppose there exists $a,\alpha_2\in(0,+\infty)$ and $p_2\in[0,\infty)$ with
\begin{equation}
    \sum_tP_0\left(D_t:=\left\{\omega:\Pi_0(\Delta_t(\omega;\cdot\,)<\mathfrak{h}^\star+a\epsilon_t)<\alpha_2e^{-p_2s_ta\epsilon_t}\right\}\right)<+\infty.
\end{equation}
Then, on $D_t^c$ it holds that
\begin{equation}
    \eta_t\le (1+p_2)a\epsilon_t-s_t^{-1}\log\alpha_2,
\end{equation}
and therefore, for all $L>\mathfrak{L}:=(1+p_2)a$ eventually $P_0$--\as that
\begin{equation}
    \eta_t\le L\epsilon_t.
\end{equation}
Moreover, 
\begin{equation}
    \limsup_{t\to+\infty}\frac{\eta_t}{\epsilon_t}\le\mathfrak{L}\qquad P_0-\as
\end{equation}
\end{theorem}
\begin{proof}
Notice first that by the first Borel--Cantelli lemma, the bound in $D_t$ may occur at most finitely often. Further, on $D_t^c$, which contains $P_0$--almost every $\omega\in\Omega$ eventually, and using \cref{lem:L^1-infty-bound-on-Z_Theta} with $u=a\epsilon_t>0$, it follows
\begin{align*}
    Z_t(\Theta)&\ge e^{-s_t(\mathfrak{h}^\star+a\epsilon_t)}\Pi_0(\Delta_t<\mathfrak{h}^\star+a\epsilon_t)\\
    &\ge e^{-s_t(\mathfrak{h}^\star+a\epsilon_t)}\alpha_2e^{-p_2s_ta\epsilon_t}\\
    &=\alpha_2\exp(-s_t(\mathfrak{h}^\star+a\epsilon_t+p_2a\epsilon_t))\\
    &=\alpha_2\exp(-s_t(\mathfrak{h}^\star+(1+p_2)a\epsilon_t))
\end{align*}
Therefore, on $D_t^c$, \ie eventually \as,
\begin{align*}
s_t^{-1}\log Z_t(\Theta)&=\Lambda_t\ge \frac{\log\alpha_2}{s_t}-\mathfrak{h}^\star-((1+p_2)a\epsilon_t)\\\iff \eta_t&\le (1+p_2)a\epsilon_t-\frac{\log\alpha_2}{s_t}
\end{align*}
We may then conclude that $\forall L>\mathfrak{L}:=(1+p_2)a$, eventually \as, 
$$\eta_t\le L\epsilon_t$$
and in particular, 
$$\limsup_{t\to+\infty}\frac{\eta_t}{\epsilon_t}\le \mathfrak{L},$$
holds $P_0$--\as
\end{proof}

\subsubsection{Pathwise averages} In addition to the previous result, we also phrase the following result that yields the same conclusion, albeit with different constants, both of which will be of use later.

\begin{theorem}\label{thm:pathwise-denominator-contraction}
Let $g_t$, a sequence of deterministic $\mathcal{T}$--measurable functions, be given. Suppose $\exists a,b,\bar{b}\ge0$ and $a'>0$ such that eventually it holds $g_t\ge \mathfrak{h}^\star-\bar{b}\epsilon_t$ $\Pi_0$--\aev, that
\begin{equation}\label{eqn:asp-overline-Phi_g_t^h^star}
    \overline{\Phi}_{g_t}^{\mathfrak{h}^\star}(a\epsilon_t)\le bs_t\epsilon_t
\end{equation}
and also that, for $E_t:=\{g_t\le \mathfrak{h}^\star+(2a+\bar{b})\epsilon_t\}$,
\begin{equation}\label{eqn:asp-summability-F_t(E_t)-for-g_t-def}
    \sum_t P_0(F_t(E_t)>a'\epsilon_t)<\infty,\quad F_t(E_t):=\fint_{E_t}(\Delta_t-g_t)\dd\Pi_0.
\end{equation}
Then, for all $L>2a+a'+b+\bar{b}=:\mathfrak{L}$ eventually $P_0$--almost surely
\begin{equation}
    \eta_t\le L\epsilon_t
\end{equation}
and moreover, it holds $P_0$--almost surely that
\begin{equation}
    \limsup_{t\to\infty}\frac{\eta_t}{\epsilon_t}\le \mathfrak{L}.
\end{equation}
\end{theorem}

\begin{proof}
We apply \cref{lem:bathtub-sandwich} to assumption \labelcref{eqn:asp-overline-Phi_g_t^h^star} where the lower bound of inequality \labelcref{eqn:bathtub-sandwich} yields that
\begin{equation*}
\Phi_{g_t}^{\mathfrak{h}^\star}((2a+\bar{b})\epsilon_t)\le bs_t\epsilon_t+\log 2
\end{equation*}
from which it follows, utilising the definition of $E_t$ and also \labelcref{eqn:def-Phi_g^c}, that
\begin{equation}\label{eqn:Pi_0-E_t-initial-bound}
\Pi_0(E_t)\ge\frac12e^{-bs_t\epsilon_t}.
\end{equation}
Next, we claim that for each $t$ sufficiently large and for each $\omega\in\{Z_t(\Theta)(\omega)<\infty\}$, which has full measure eventually by \cref{lem:finite-Z_t-Theta}, it holds that
\begin{equation}\label{eqn:log-Z_t_Theta-bound}
    \log Z_t(\Theta)(\omega)\ge \log\Pi_0(E_t)-s_t\fint_{E_t}\Delta_t(\omega;\cdot)\dd\Pi_0.
\end{equation}
\begin{proof}[Proof of \labelcref{eqn:log-Z_t_Theta-bound}]
Since \labelcref{eqn:Pi_0-E_t-initial-bound} holds eventually $ \Pi_0(E_t)\in(0,1]$, and the first term is unproblematic. In the case $Z_t(\Theta)(\omega)=0$ implies $\Delta_t(\omega;\cdot)=+\infty$ $\Pi_0$--\aev so the second, by definition of $E_t$, is trivially infinite, and the estimate holds for $Z_t(\Theta)(\omega)=0$. The case $Z_t(\Theta)(\omega)=+\infty$ is excluded. 

Assuming now that $Z_t(\Theta)(\omega)\in(0,\infty)$, restricting the domain of integration of $Z_t(\Theta)$ to the controlled $Z_t(E_t)$ and applying Jensen's inequality, we obtain
\begin{align*}
\log Z_t(\Theta)&=\log\int_\Theta R_t\dd\Pi_0\\
&\ge \log\int_{E_t}R_t\dd\Pi_0\\
&=\log\Pi_0(E_t)\fint_{E_t}R_t\dd\Pi_0\\
&={\log\Pi_0(E_t)}+\log\int_{E_t}R_t\dd\Pi_{0}|_{E_t}\\
&\ge {\log\Pi_0(E_t)}+\int_{E_t}\log R_t\dd\Pi_{0}|_{E_t}
\end{align*}
where we have used that we may condition the measure $\Pi_0(\cdot)$ on the set $E_t$, to yield the probability measure $\Pi_{0}|_{E_t}(\cdot)=\frac{\Pi_0(\cdot\,\cap E_t)}{\Pi_0(E_t)}$ to which we may apply Jensen's inequality. Thus, estimate \labelcref{eqn:log-Z_t_Theta-bound} is established.
\end{proof}
In particular, making use of the estimate, we have for each $t$ the inclusion
\begin{equation}\label{eqn:set-inclusion-eta_t-bound}
    \left\{\eta_t>(2a+a'+b+\bar{b})\epsilon_t +\frac{\log 2}{s_t}\right\}\subseteq\{F_t(E_t)>a'\epsilon_t\}\cup N_t\subseteq\Omega.
\end{equation}
where $N_t$ is a $P_0$ null set.

Indeed by inequality \labelcref{eqn:Pi_0-E_t-initial-bound},
$$\frac{\log \Pi_0(E_t)}{s_t}\ge -b\epsilon_t-\frac{\log 2}{s_t}$$
and therefore, by the estimate \labelcref{eqn:log-Z_t_Theta-bound} we also have, as subsets of $\Omega$, that
\begin{equation}\label{eqn:bound-Lambda_t+prior-condition_E_t-bounding-neg-avg-Delta_t}
\begin{gathered}
\left\{\frac{1}{s_t}\log Z_t(\Theta)\ge\frac{\log\Pi_0(E_t)}{s_t}-\fint_{E_t}\Delta_t\dd\Pi_0\right\}\\
\subseteq\left\{\frac{1}{s_t}\log Z_t(\Theta)+b\epsilon_t+\frac{\log 2}{s_t}\ge -\fint_{E_t}\Delta_t\dd\Pi_0\right\}.
\end{gathered}
\end{equation}
up to inclusion of $P_0$--measure zero sets.

Now, by definition of $E_t$ we have that 
\begin{equation*}
    \fint_{E_t}g_t\dd\Pi_0\le \mathfrak{h}^\star+(2a+\bar{b})\epsilon_t,
\end{equation*}
and
\begin{equation*}
    \fint_{E_t}\Delta_t\dd\Pi_0\le \mathfrak{h}^\star+(2a+\bar{b})\epsilon_t+F_t(E_t).
\end{equation*}
It therefore follows that, modulo a $P_0$--null set, we have the set inclusion
\begin{align*}
 \left\{F_t(E_t)>a'\epsilon_t\right\}&=\left\{-\fint_{E_t}g_t\dd\Pi_0-a'\epsilon_t>-\fint_{E_t}\Delta_t\dd\Pi_0\right\}\\
 &\supseteq\left\{-\mathfrak{h}^\star-(2a+\bar{b})\epsilon_t-a'\epsilon_t>-\fint_{E_t}\Delta_t\dd\Pi_0\right\}\\
\intertext{and, by the inclusion \labelcref{eqn:bound-Lambda_t+prior-condition_E_t-bounding-neg-avg-Delta_t}, possibly up to another $P_0$ null set, that}
 &\supseteq\left\{-\mathfrak{h}^\star-(2a+\bar{b})\epsilon_t-a'\epsilon_t>\frac{1}{s_t}\log Z_t(\Theta)+b\epsilon_t+\frac{\log 2}{s_t}\right\}\\
 &=\left\{-\mathfrak{h}^\star-\frac1{s_t}\log Z_t(\Theta)>\epsilon_t\left(2a+a'+b+\bar{b}+\frac{\log 2}{s_t\epsilon_t}\right)\right\}\\
 &=\left\{\eta_t>\epsilon_t\left(2a+a'+b+\bar{b}+\frac{\log 2}{s_t\epsilon_t}\right)\right\}
\end{align*}
where the final step has used the definition of $\eta_t$, establishing the inclusion \labelcref{eqn:set-inclusion-eta_t-bound}. The sets $N_t$ are a countable collection of  $P_0$--measure zero sets, and hence their union is $P_0$--null and thus does not contribute to the summation.

Finally, as the sequence $P_0(F_t(E_t)>a'\epsilon_t)$ is summable, it follows that 
$$\sum_tP_0\left(\eta_t>\left(2a+a'+b+\bar{b}+\frac{\log 2}{s_t\epsilon_t}\right)\epsilon_t\right)<+\infty$$ 
and therefore, by the first Borel--Cantelli lemma the bound may occur only finitely many times.
That is, eventually $P_0$--almost surely
\begin{equation}
    \eta_t\le \epsilon_t(2a+a'+b+\bar{b})+\frac{\log2}{s_t}
\end{equation}
and hence, for any $L>2a+a'+b+\bar{b}$, it holds for all $t$ sufficiently large $P_0$ almost surely that
\begin{equation*}
    \eta_t\le L\epsilon_t
\end{equation*}
and hence that $P_0$--almost surely
\begin{equation}
    \limsup_{t\to+\infty}\frac{\eta_t}{\epsilon_t}\le 2a+a'+b+\bar{b}.
\end{equation}
The proof of \cref{thm:pathwise-denominator-contraction} is hence complete.
\end{proof}
In particular, we may optimise the constants $(a,a',b,\bar{b})$ as to minimise the limit superior bound $\mathfrak{L}:=2a+a'+b+\bar{b}$.

We note that a canonical choice for $g_t$ will be a variation of $K_t:=\E_{P_0^{(t)}}\Delta_t$, as $\E_{P_0^{(t)}}F_t(E_t)=0$, but we may modify $g_t$ off the set $E_t$, as defined in the theorem by for instance truncating $K_t$ to enforce the essential infimum of $g_t$ to be bounded below or to restrict the variance. In particular, the following lemma licenses this.
\begin{lemma}\label{lem:exchanging-integration-for-K_t-and-Delta_t}
It holds for all but finitely many $t$ that
\begin{equation}
    \int_{\Xi^{(t)}}\fint_{E_t}\Delta_t\dd\Pi_0\dd P_0^{(t)}=\fint_{E_t} K_t\dd\Pi_0
\end{equation}
\end{lemma}

\begin{proof}
The condition is not trivial as Fubini--Tonelli does not immediately apply. Write
\begin{equation}
    \E_{P_0^{(t)}}\overline{\Delta}_t(E)=\frac{1}{s_t\Pi_0(E)}\int_{\Xi^{(t)}}\int_{E}-\log R_t\dd\Pi_0\dd P_0^{(t)}.
\end{equation}
We may then split $(-\log R_t)$ into its positive and negative parts. On the positive part, we may apply Tonelli without problem. The negative part is slightly more involved. First notice that the negative part equals $(\log R_t)\vee 0$, and hence using $\log^+ x\le x$ on $x\ge0$, we have that
\begin{equation}
    \E_{P_0^{(t)}}[(\log R_t)\vee 0]\le \E_{P_0^{(t)}}R_t=:c_t
\end{equation}
since $R_t\ge 0$ by definition, and we may increase the domain of integration.

From this it follows that for any $p>1$
\begin{equation*}
\int_{E}\int_{\Xi^{(t)}}(\log R_t)\vee 0\dd P_0^{(t)}\dd\Pi_0\le\int_Ec_t\dd\Pi_0\le {||c_t||_{L^1(\Pi_0)}}\lesssim||c_t||_{L^{p,\infty}(\Pi_0)},
\end{equation*}
where we have used the natural embedding of $L^{p,\infty}\hookrightarrow L^1$ on a probability space. The exact bound can be verified by the Hardy--Littlewood inequality \labelcref{eqn:hardy-littlewood} and thus is $q=p/(p-1)$, but the exact constant is unimportant here.
By \cref{asp:lorentz}, it thus follows that there exists a $p>1$ such that the above is finite for all but finitely many $t$. The exchange of the order of integration for the negative part is hence licensed by Fubini for all but finitely many $t$, and thus the desired conclusion follows.
\end{proof}

\subsection{A converse denominator estimate: proof of \texorpdfstring{\cref{thm:converse-to-rates}}{Theorem \ref*{thm:converse-to-rates}} completed}
We complete the section by proving \cref{thm:converse-to-rates}.

Recall the following deterministic function
\begin{equation}
    A_t(P):=-\frac{1}{s_t}\log c_t(P).
\end{equation}
and denote by 
\begin{equation*}
    \mathfrak{a}_t^\star:=\essinf_{P\in\Theta} A_t(P).
\end{equation*}
In particular, contraction rates through our and virtually every other approach require that almost surely $\limsup_t\eta_t/\epsilon_t\le \mathfrak{L}$, if we have a lower bound on $\eta_t$ that increases faster than $\epsilon_t$, then contraction at the prescribed rate is impossible.\footnote{Although this assumption is not strictly necessary, not imposing this would require control of how the numerator and denominator interact.}

The proof of \cref{thm:converse-to-rates} follows by combining the following two results.

\begin{lemma}
It holds that
\begin{equation}\label{eqn:Expectation-of-Z_Theta-relation-to-free-energy-A_t}
\E_{P_0^{(t)}}[Z_t(\Theta)]=\exp\left(-s_t\mathcal{E}_{\Pi_0}[A_t](s_t)\right)
\end{equation}
and further that pointwise
\begin{equation}
    \limsup_{t\to\infty}A_t(P)\le h(P).
\end{equation}
\end{lemma}
\begin{proof}
An immediate application of Tonelli yields that
\begin{align*}
\mathcal{E}_{\Pi_0}[A_t](s_t)&=-\frac{1}{s_t}\log\int_\Theta e^{-s_tA_t}\dd\Pi_0\\&=-\frac{1}{s_t}\log\int_\Theta\int_{\Xi^{(t)}} R_t\dd P_0^{(t)}\dd\Pi_0\\
&=-\frac{1}{s_t}\log\int_{\Xi^{(t)}} Z_t(\Theta)\dd P_0^{(t)}=-\frac1{s_t}\log\E_{P_0^{(t)}}[Z_t(\Theta)]
\end{align*}
which yields the first. For the second, an immediate application of Jensen's inequality yields that 
\begin{equation*}
    A_t(P)\le K_t(P):=\frac{1}{s_t}\E_{P_0^{(t)}}[-\log R_t]
\end{equation*}
where the right hand side is assumed to converge, pointwise in $P$, to $h(P)$, by assumption \labelcref{asp:1} which completes the proof.
\end{proof}

\begin{theorem}\label{prop:nonasymptotic-exp-tauberian-thm}
    For $g$ as before, it holds for all $s>0$ that
    \begin{equation}\label{eqn:nonasymptotic-Tauberian-slack}
        \mathfrak{g}^\star+\frac12\mathcal{J}_g(s)-s^{-1}\log2\le \mathcal{E}_{\Pi_0}[g](s)\le \mathfrak{g}^\star+\mathcal{J}_g(s).
    \end{equation}
\end{theorem}
\begin{proof}
The upper bound is immediate.
For the lower bound, consider first
\begin{equation*}
e^{s\mathfrak{g}^\star}\int_\Theta e^{-sg}\dd\Pi_0=\int_\Theta e^{-s(g-\mathfrak{g}^\star)}\dd\Pi_0
\end{equation*}
Then, using the layer--cake identity for $x\mapsto e^{-sx}$ on $x\in[0,+\infty]$ and a change of variables, we are led to
\begin{equation*}
    \int_\Theta e^{-s(g-\mathfrak{g}^\star)}\dd\Pi_0=\int_0^\infty se^{-sv}\,\Pi_0(g-\mathfrak{g}^\star\le v)\dd v
\end{equation*}
where we may forego the absolute values as $g-\mathfrak{g}^\star$ is almost everywhere nonnegative. Then taking $\Phi_g(u)=-\log\Pi_0(g\le\mathfrak{g}^\star+u)$, we get
\begin{equation*}
    =\int_0^\infty se^{-sv-\Phi_g(v)}\dd v,
\end{equation*}
since $u\mapsto \Phi_g(u)$ is right continuous. 

Now
\begin{equation*}
    \frac12(sv+\Phi_g(v))+\frac12sv\ge\frac{s}{2} \mathcal{J}_g(s)+\frac{sv}2.
\end{equation*}
as $\mathcal{J}_g(s)$ is an infimum over all $u>0$ independent of $v$.

And as $\Phi_g$ is nonnegative, it follows
\begin{equation*}
-sv-\Phi_g(v)\le -sv-\frac12\Phi_g(v)
\end{equation*}
which then implies
\begin{align*}
e^{-sv-\Phi_g(v)}&\le e^{ -sv-\frac12\Phi_g(v)}\\
&=e^{-\frac12(sv+\Phi_g(v))}e^{-\frac12sv}\\
&\le e^{-\frac s2\mathcal{J}_g(s)}e^{-\frac{sv}2}.
\end{align*}

It thereby follows that
\begin{align*}
\int_0^\infty se^{-sv-\Phi_g(v)}\dd v &\le \int_0^\infty se^{-\frac s2\mathcal{J}_g(s)}e^{-\frac{sv}2}\dd v\\
&= \exp\left({-\frac s2\mathcal{J}_g(s)}\right) \int_0^\infty se^{-\frac{sv}2}\dd v=2\exp\left({-\frac s2\mathcal{J}_g(s)}\right).
\end{align*}
Taking logarithms and rearranging yields the desired lower bound.
\end{proof}

\begin{proof}[Proof of \cref{thm:converse-to-rates} completed]
    The first part is a standard application of Markov's inequality with equality \labelcref{eqn:Expectation-of-Z_Theta-relation-to-free-energy-A_t}, where the summability follows by standing assumption on $s_t$. 
    
    The second follows by the first Borel--Cantelli lemma applied to the first part which yields the eventual almost sure conclusion, and the estimate \labelcref{eqn:thm-rates-converse-eta_t-lower-bound} therefore follows by direct application of inequality \labelcref{eqn:nonasymptotic-Tauberian-slack} of \cref{prop:nonasymptotic-exp-tauberian-thm}.
\end{proof}

\part{Analysis of the contraction}\label{part:contraction-analysis}

\section{Verification of assumptions}\label{sec:verification-of-assumptions}

This section provides the proof of the example, \cref{thm:DP-W_1-example,lem:DP-W_1-small-ball-estimates}, as well as the simplified assumptions \cref{lem:pointwise-a-s-convergence-to-h,lem:uniform-convergence-on-sieve}.
\subsection{On simplified assumptions}
\begin{proof}[Proof of \cref{lem:pointwise-a-s-convergence-to-h}]
As the limit exists, it must equal the limsup. Hence, for \aev $P\in\Theta$
$$
P_0\left(\lim_{t\to\infty}\frac1{s_t}\log R_t(\omega,P)=-\limsup_{t\to\infty}\Delta_t(X_1^t(\omega);P)=-h(P)\right)=1
$$
as claimed.
\end{proof}
\begin{proof}[Proof of \cref{lem:uniform-convergence-on-sieve}]
From \eqref{eqn:uniform-convergence-on-sieve}, eventually \as
$$
\esssup_{G_t}U_t\le c_1\lambda_t\iff\Pi_0(G_t\cap\{U_t>u\})=0\quad\forall u\ge c_1\lambda_t
$$
yielding \eqref{eqn:asp-concentration-sieve}. Further, if $\lambda_t\le\epsilon_t$, then $\forall u\ge c_1j\epsilon_t\ge c_1\epsilon_t\ge c_1\lambda_t$
$$
0\le\Pi_0\rest (G_t\cap S_t^{(j)}) (\{U_t>u\})\le \Pi_0\rest G_t (\{U_t>u\})=0
$$
and \cref{asp:lorentz-sieve-shell-interaction} follows.
\end{proof}
\subsection{Proof of \texorpdfstring{\cref{lem:DP-W_1-small-ball-estimates}}{Lemma \ref*{lem:DP-W_1-small-ball-estimates}}}
Steps 2 and 3 of the intermediate argument have been adapted from Nguyen \cite{nguyen2013}.
\begin{proof}[Proof of \cref{lem:DP-W_1-small-ball-estimates}]
We divide the argument into steps.
\begin{step}An estimate of the packing number.
Fix $\nu\in\P(K)$ and $\epsilon>0$.
Let $\{x_1,\dots,x_N\}\subset K$ be a maximal $\epsilon$--separated subset of $K$, where $N=N_\mathrm{pack}(K,\rho,\epsilon)$ is the $\epsilon$--packing number. As $K$ is compact, this is finite, and as $G_0$ is AD--regular, $N\sim\epsilon^{-d}$. 

Indeed, as $(x_i)_{i=1}^N$ is maximal $\epsilon$--separated, the open balls $(B(x_i,\epsilon/2))_{i=1}^N$ are disjoint
\begin{gather*}
1\ge \sum_{i=1}^{N}G_0(B(x_i,\epsilon/2)\cap K)\ge NC^{-1}\left(\frac\epsilon2\right)^d\implies N\lesssim\epsilon^{-d},\\
\intertext{and $(B(x_i,\epsilon))_{i=1}^N$ covers $(K,\rho)$ with open $\epsilon$--balls,}
1=G_0(K)\le \sum_{i=1}^{N}G_0(B(x_i, \epsilon)\cap K)\le N C\epsilon^{d}\implies N\gtrsim \epsilon^{-d}.
\end{gather*}
Since $K$ admits an AD--regular probability measure, $d\ge0$ is the Hausdorff dimension of $K$ and equals the Assouad dimension (see Theorem 6.4.1 in \cite{fraser2021}). 
\end{step}
\begin{step}A partition of $K$.
We construct a partition of $K$ as follows, let $S_i=B(x_i,\epsilon/2)$, which are disjoint for all $i\in[1;N]$, by assumption of being $\epsilon$--separated. Assign the remaining elements of $K$ to the set $S_i'$ with the nearest point $x_i$ in the packing $(x_i)_{i=1}^N$, where ties are broken by index precedence, partitioning $K$. It is immediate that $S_i\subseteq S_i'\subset B(x_i,\epsilon)$ and $\mathrm{diam}(S_i')\le 2\epsilon$.

Denote by $g_i=G_0(S_i')$, $m_i=P(S_i')$, and $\hat p_i=\nu(S_i')$. By definition of the Dirichlet process, since $P\sim\Pi_0^{(\mathrm{DP})}$, it follows 
$$
(m_1,\dots,m_N)\sim\mathrm{Dir}(\alpha g_1,\dots,\alpha g_N).
$$
We may take $\hat p_N=\max_i\hat p_i\ge 1/N$ \wolog after relabelling.
\end{step}

\begin{step}A first upper bound on $W_1$--balls under the DP prior.
Next, we establish an upper bound on the $1$--Wasserstein metric by construction of a coupling of $P$ and $\nu$. In particular, take $m_i\wedge \hat p_i$ as the probability mass of $\nu$ and $P$ on $S_i'$ for each $i\in[1;N]$. The remaining mass, $\frac12||m-\hat p||_1$, may be associated arbitrarily. Then, taking the expectation of this yields
\begin{equation}\label{eqn:proof-dp-prior-upper-bound-W_1}
    W_1(P,\nu)\le 2\epsilon+{\textstyle \frac12}||m-\hat p||_1\mathrm{diam}(K).
\end{equation}
Since $G_0$ is AD--regular, it is atomless, and we may take $\epsilon>0$ sufficiently small so that $\alpha g_i\le 1$ for all $i\in[1;N]$. Then, letting $\delta=\epsilon/\mathrm{diam}(K)$, we find from estimate \labelcref{eqn:proof-dp-prior-upper-bound-W_1}
\begin{align}
\Pi&_0^{(\mathrm{DP})}(P:W_1(P,\nu)\le 3\epsilon)
\ge \Pi_0^{(\mathrm{DP})}\left(||m-\hat p||_1\le \delta\right).\label{eqn:proof-dp-prior-ineq-1}\\
\intertext{Letting $\eta:=1/(N(N-1))\wedge \delta/(2N)$ and as $||m-\hat p||_1\le2\sum_{i=1}^{N-1}|m_i-\hat p_i|$}
&\ge \Pi_0^{(\mathrm{DP})}(|m_i-\hat p_i|\le \eta, \forall i\in[1;N-1])\label{eqn:proof-dp-prior-ineq-2}.\\
\intertext{Integrating over the $N-1$ dimensional simplex, since $\eta\le 1/(N(N-1))$ implies $m\in\Delta_{N-1}$ for $\{|m_i-\hat p_i|\le \eta, \forall i\in[1;N-1]\}$, yields}
&=\cfrac{\Gamma(\alpha)}{\prod_{i=1}^{N}\Gamma(\alpha g_i)}\int_{\Delta_{N-1}\cap\left\{|m_i-\hat p_i|\le \eta\right\}}\prod_{i=1}^{N-1}m_i^{\alpha g_i-1}\left(1-\sum_{i=1}^{N-1}m_i\right)^{\alpha g_N-1}d m_i\\
\intertext{by definition of a DP. Now, utilising $(1-\sum_{i=1}^{N-1}m_i)^{\alpha g_N-1}=m_N^{
\alpha g_N-1}\ge1$ as $\alpha g_N\le 1$}
&\ge \frac{\Gamma(\alpha)}{\prod_{i=1}^{N}\Gamma(\alpha g_i)}\prod_{i=1}^{N-1}\int_{(\hat p_i-\eta)\vee 0}^{(\hat p_i+\eta)\wedge 1}m_i^{\alpha g_i-1}\dd m_i\\
\intertext{and finally, since $\Gamma(a)\le 1/a$ for $a\in(0,1]$ and elementary integration,}
&\ge \Gamma(\alpha)\eta^{N-1}\prod_{i=1}^N\alpha g_i.
\end{align}
\end{step}

\begin{step}A refined estimate utilising AD--regularity.
Making use of the AD--regularity of $G_0$, the previous yields
\begin{align}
\Pi_0^{(\mathrm{DP})}(P:W_1(P,\nu)\le 3\epsilon)&\ge {\Gamma(\alpha)\alpha^N}\eta^{N-1}\prod_{i=1}^NG_0(S_i')\notag\\
&\ge {\Gamma(\alpha)\alpha^N}\eta^{N-1}\prod_{i=1}^NG_0(B(x_i,\epsilon/2)\cap K)\notag\\
&\ge \frac{\Gamma(\alpha)\alpha^N}{N^{N-1}}\left(\frac{\epsilon}{2\mathrm{diam}(K)}\wedge\frac{1}{N-1}\right)^{N-1} C^{-N}\left(\frac\epsilon2\right)^{Nd}\notag
\end{align}
Taking logarithms and rearranging yields the estimate
\begin{multline}\label{eqn:proof-dp-prior-final-bound}
-\log \Pi_0^{(\mathrm{DP})}(P:W_1(P,\nu)\le 3\epsilon)\le (N-1)\log N-\log\Gamma(\alpha)-N\log \alpha\\+(N-1)\log\left(\frac{2\mathrm{diam}(K)}{\epsilon}\vee ({N-1})\right)+N\log C-Nd\log \frac\epsilon2
\end{multline}
Now, as $N\sim\epsilon^{-d}$, if $\epsilon<e^{-1}$, all terms are of order $\lesssim\epsilon^{-d}\log(1/\epsilon)$. Hence,
$$
\epsilon_0=\min\left(\frac{\mathrm{diam}(K)}{2},\frac{(\alpha C)^{-1/d}}{2},\frac1e\right)
$$
suffices throughout. 
\end{step}
And therefore, estimate \labelcref{eqn:proof-dp-prior-final-bound} establishes the result after redefining $\epsilon$ as $\epsilon/3$ where the $1/3$ scaling exclusively affects the constant $C_0$.
\end{proof}

\subsection{The DP prior with ergodic observables: proof of \texorpdfstring{\cref{thm:DP-W_1-example}}{Theorem \ref*{thm:DP-W_1-example}}}
We assume throughout this subsection the setting of \cref{thm:DP-W_1-example}.
\begin{lemma}\label{lem:DP-W_1-limiting-ptoperties}
The following hold true $$h(P)=W_1(P,\mu),$$ $$\tilde r_t:=W_1(\mu,\mathbb P_t)\to0\qquad\mathrm{a.s.,}$$ $$|\Delta_t(X_1^t;P)-h(P)|\le 2\tilde r_t,\quad\mathrm{and}$$ $$\Theta_\star=\{\mu\}.$$
\end{lemma}
\begin{proof}
Let $(f_j)_{j=1}^\infty\subset C_b(K)$ be a dense, countable sequence. Then by Birkhoff's pointwise ergodic theorem (cf. \cite{kerr2016}), it is immediate, as $\nu_0$ is ergodic, that for all $j\in\N$ $\nu_0$--\as
\begin{equation*}
\frac1t\sum_{i=0}^{t-1}f(T^iX_1)\overset{t\to\infty}{\longrightarrow}\int_K f\dd\mu. 
\end{equation*}
And hence, $\mathbb P_t\overset{w^\star}{\to}\mu$, in the weak$^\star$ sense.
Then, as the $1$--Wasserstein metric metrizes weak$^\star$ convergence of probability measures on compact metric spaces (see Theorems 6.9 and 6.18 in \cite{villani2009}), $$W_1(P,\mathbb P_t)\to W_1(P,\mu)$$ $\nu_0$--\as As $W_1(P_1,P_2)\le\mathrm{diam}(K)$ for all $P_1,P_2\in\mathcal P(K)$, the results follow by the bounded convergence theorem.
\end{proof} 
\begin{proposition}
For all $\epsilon>0$, 
$$ \Pi_t^{(\mathrm{DP})}(\{P\in\Theta: W_1(P,\mu)\ge \epsilon\}|X_1^t)\to0$$
$\nu_0$--almost surely. Moreover, it is $\nu_0$--\as that
\begin{equation*}
\Pi_t^{(\mathrm{DP})}(dP|X_1^t)\to \delta_{\mu}(dP).
\end{equation*}
\end{proposition}
\begin{proof}
As $\Pi_0^{(\mathrm{DP})}(\{P:h(P)\le\epsilon\})>0$ for all $\epsilon>0$ by \cref{lem:DP-W_1-small-ball-estimates} and \cref{lem:DP-W_1-limiting-ptoperties}, the assumptions of \cref{thm:consistency} are satisfied. Moreover, $(\Theta,W_1)$ is a metric space, $\Theta_\star=\{\mu\}$, and $h(P)>0$ for all $P\ne \mu$, \cref{thm:weak-star-convergence-unique-min} yields the final conclusion.
\end{proof}

\begin{proof}[Proof of \cref{thm:DP-W_1-example}]
We demonstrate the assumptions of \cref{thm:rates} are satisfied. The standing measurability hypothesis holds trivially, in fact $P\mapsto W_1(P,\nu)$ is continuous for all $\nu\in\mathcal P(K)$, and $\beta_t>0$ is chosen appropriately.
The condition on $h(P):=W_1(P,\mu)$ holds from the previous lemma with $I=\emptyset$ and $\mathfrak{h}^\star=0$.

For \cref{asp:lorentz}, we have by the triangle inequality
\begin{gather*}
|\Delta_t(X_1^t;P)|=|W_1(P,\mathbb P_t)-W_1(\mu, \mathbb P_t)|\le W_1(P,\mu)=h(P)\le\mathrm{diam}(K)\\\intertext{and thus}
R_t(P)=e^{-s_t\Delta_t(X_1^t;P)}\le e^{s_th(P)}\le e^{s_t\mathrm{diam}(K)}.
\end{gather*}
Hence $c_t(P)\le e^{s_t\mathrm{diam}(K)}$ and we may take $p=\infty$, and $\kappa=\mathrm{diam}(K)$.

For \cref{asp:sieve-gamma} we take $G_t\equiv \Theta$, so this is trivially satisfied.

In \cref{asp:lorentz-sieve-shell-interaction}, notice
\begin{equation*}
    \Delta_t(X_1^t;P)=W_1(P,\mathbb P_t)-W_1(\mu, \mathbb P_t)\ge h(P)-2W_1(\mu,\mathbb P_t)=h(P)-2r_t
\end{equation*}
and hence
\begin{equation*}
    \esssup_{P\in\Theta}\left(\frac1{s_t}\log R_t(P)+h(P)\right)\le \esssup_{P\in\Theta}\left(h(P)-(h(P)-2r_t)\right)\le 2r_t.
\end{equation*}
So we may apply \cref{lem:uniform-convergence-on-sieve}, with $p_1=\infty$, for all $c_1\in(0,1)$ satisfying $2r_t\le c_1 \epsilon_t$.

For \cref{asp:averaged-g_t}, pick $g_t=h$, so $\bar b=0$. Let $a>0$, by \cref{lem:DP-W_1-small-ball-estimates}, and as $\epsilon_t\to0$, eventually $a\epsilon_t<\epsilon_0$ and hence
$$
-\log\Pi_0^{(\mathrm{DP})}(\{P:h(P)\le a\epsilon_t\})\lesssim (a\epsilon_t)^{-d}\log\left(\frac1{a\epsilon_t}\right).
$$
We now claim, for $\epsilon_t$ as specified, that eventually
\begin{equation}\label{eqn:DP-W_1-thm-proof-claim}
    \epsilon_t^{-d}\log\left(\frac1{\epsilon_t}\right)\lesssim s_t\epsilon_t,\qquad\mathrm{and}\qquad\epsilon_ts_t\to\infty.
\end{equation}
\begin{proof}[Proof of estimate \labelcref{eqn:DP-W_1-thm-proof-claim}]
Let 
$$
b_t:=\left(\frac{\log s_t}{s_t}\right)^{1/(d+1)},
$$
then $s_t\epsilon_t\ge As_tb_t=As_t^{d/(d+1)}(\log s_t)^{1/(d+1)}\to\infty$. Further, $x\mapsto x^{-d}\log(1/x)$ is decreasing for $x$ sufficiently small, and satisfied for all $\epsilon_t$ with $t$ sufficiently large since $\epsilon_t\to0$. Hence,
\begin{align*}
\epsilon_t^{-d}\log\left(\frac1{\epsilon_t}\right)\le (Ab_t)^{-d}\log\left(\frac1{Ab_t}\right)\le A^{-d}b_t^{-d}\log\frac1{b_t}.
\end{align*}
as $A\ge 1$. Plugging in the definition of $b_t$,
\begin{align}
\log\frac1{b_t}&=\frac1{d+1}\left(\log s_t-\log\log s_t\right)\le \frac{\log s_t}{d+1}\notag\\
\intertext{implying}
b_t^{-d}\log\frac1{b_t}&\le \frac1{d+1}s_t^{d/(d+1)}(\log s_t)^{1/(d+1)}\notag\\
&=\frac{1}{d+1}\frac{s_t}{s_t^{1/(d+1)}}(\log s_t)^{1/(d+1)}=\frac1{d+1}s_tb_t,\notag\\
\intertext{and hence}
\epsilon_t^{-d}&\log\left(\frac1{\epsilon_t}\right)\le \frac{A^{-d}}{d+1}s_tb_t\le \frac{A^{-(d+1)}}{d+1}s_t\epsilon_t
\end{align}
establishing estimate \labelcref{eqn:DP-W_1-thm-proof-claim}.
\end{proof}
Utilising the estimate, eventually
$$
-\log\Pi_0^{(\mathrm{DP})}(\{P:h(P)\le a\epsilon_t\})\lesssim \epsilon_ts_t,
$$
and by \cref{thm:optimal-choice-E_t}, 
$$
\overline\Phi_h^{\mathfrak{h}^\star}(a\epsilon_t)\lesssim s_t\epsilon_t.
$$
For the second part of \cref{asp:averaged-g_t}, notice that $|\Delta_t(X_1^t;P)|\le h(P)$ for all $x\in K$, and hence $\Delta_t-h\le 0$, and hence 
$$
\nu_0\left(\fint_{E_t}(\Delta_t-h)\dd\Pi_0^{(\mathrm{DP})}>a'\epsilon_t\right)=0
$$
for all $t$. We have established all assumptions of \cref{thm:rates} hold. The proof of \cref{thm:DP-W_1-example} is complete.
\end{proof}

\section{Gibbs posterior consistency}\label{sec:consistency-proof}
\subsection{Proof of \texorpdfstring{\cref{thm:contraction_bad-sets}}{Theorem \ref*{thm:contraction_bad-sets}} completed}
We first prove the following theorem which we then use to conclude the proof.
\begin{theorem}\label{prop:posterior-consistency-on-exponentially-thin-sets}
Assume the setting of \cref{prop:numerator-bad-sets} where $\gamma/q>\mathfrak{h}^\star+\kappa$, and suppose also that almost surely
\begin{equation}\label{eqn:asp-liminf-bound-Z_t(Theta)}
\liminf_{t\to+\infty} s_t^{-1}\log Z_t(\Theta)\ge -\mathfrak{h}^\star,
\end{equation}
then eventually $P_0$--almost surely for all $\epsilon\in(0,\gamma/q-\kappa-\mathfrak{h}^\star)$
\begin{equation}
    \Pi_t(B_t|X_1^t)\lesssim e^{-\epsilon s_t}.
\end{equation}
\end{theorem}
\begin{proof}
Writing
\begin{equation*}
\Pi_t(B_t|X_1^t)=\frac{Z_t(B_t)}{Z_t(\Theta)}
\end{equation*}
by \cref{prop:numerator-bad-sets} and assumption \labelcref{eqn:asp-liminf-bound-Z_t(Theta)}, for all $\epsilon>0$ it holds eventually \as that
\begin{equation*}
\begin{aligned}
\frac{Z_t(B_t)}{Z_t(\Theta)}&\le\frac{\exp(-s_t(\gamma/q-\kappa-\epsilon/2))}{\exp(-s_t(\mathfrak{h}^\star+\epsilon/2))}\\
&=\exp(s_t(\mathfrak{h}^\star-\gamma/q+\kappa+\epsilon)).
\end{aligned}
\end{equation*}
Now, as $\gamma/q>\mathfrak{h}^\star+\kappa$ we may choose $\epsilon>0$ sufficiently small so that the convergence holds $P_0$--\as, \ie $\epsilon\in(0,\gamma/q-\kappa-\mathfrak{h}^\star)$. Taking $$\epsilon'=\gamma/q-\kappa-\mathfrak{h}^\star-\epsilon\in(0,\gamma/q-\kappa-\mathfrak{h}^\star)$$ thus yields the desired bound.
\end{proof}
\begin{proof}[Proof of \cref{thm:contraction_bad-sets} completed]
We establish first that assumption \labelcref{eqn:asp-liminf-bound-Z_t(Theta)} holds. Indeed,
under \labelcref{eqn:asp-consistency}, we have that assumption \labelcref{eqn:limsup-convergence-Delta_t} is satisfied and hence \cref{thm:qualitative-denominator-floor} yields $$\limsup_{t\to+\infty}\eta_t\le 0.$$ Under \cref{asp:averaged-g_t} we may apply \cref{thm:pathwise-denominator-contraction} and under \cref{asp:rates-endpoint-ver} we make use of \cref{thm:denominator-via-lorentz-bound}, which both yield some $L>0$ for which it eventually \as holds that
$$
\eta_t\le L\epsilon_t\quad\implies\quad\limsup_{t\to+\infty}\eta_t\le 0.
$$
By definition $\eta_t=-\mathfrak{h}^\star-\Lambda_t$ and hence eventually \as
\begin{gather*}
    0\ge\limsup_{t\to+\infty}(-\mathfrak{h}^\star-\Lambda_t)=-\mathfrak{h}^\star+\limsup_{t\to+\infty}(-\Lambda_t)=-\mathfrak{h}^\star-\liminf_{t\to+\infty}\Lambda_t\\
    \implies\liminf_{t\to+\infty}\frac1{s_t}\log Z_t(\Theta)\ge -\mathfrak{h}^\star
\end{gather*}
By the setting of the theorem we are now in a position to apply \cref{prop:posterior-consistency-on-exponentially-thin-sets} which completes the proof.
\end{proof}

\subsection{Proof of \texorpdfstring{\cref{thm:consistency}}{Theorem \ref*{thm:consistency}} completed}
\begin{bareproof}
Recall that
\begin{align*}
    \Pi_t(A|X_1^t)&=\Pi_t(A\cap G_t|X_1^t)+\Pi_t(A\cap G_t^c|X_1^t)\\
    &\le \Pi_t(A\cap G_t|X_1^t) + \Pi_t(G_t^c|X_1^t).
\end{align*}
By \cref{thm:contraction_bad-sets}, $\Pi_t(G_t^c|X_1^t)\to 0$ \as so it remains to prove that for such sets $A$, $\Pi_t(A\cap G_t|X_1^t)\to0$ \as 
Further
\begin{equation}\label{eqn:proof-consistency-Pi_t-A_t-G_t-equality}
    \Pi_t(A\cap G_t|X_1^t)=\frac{Z_t(A\cap G_t)}{Z_t(\Theta)}.
\end{equation}
Note first that if $\Pi_0(A)=0$, then the numerator of the above is trivially zero for all $t$, and by standing assumption the denominator is \as positive, thus the result follows for this case.

Assume now $\Pi_0(A)>0$, and note that the assumption that $\mathfrak{h}(A)>\mathfrak{h}^\star$ yields, invoking \cref{lem:h(A)>h*=>Pi_0(A)<1}, that $\Pi_0(A)<1$. 
Now, since $\Pi_0(A)\in(0,1)$, by \cref{prop:first-sieve-estimate} we have that for all such sets $A$, eventually $P_0$--almost surely
\begin{equation}\label{eqn:proof-consistency-numerator-on-sieve-bound}
Z_t(A\cap G_t)\le e^{-s_t(
\mathfrak{h}(A)\wedge M_0)}\left(e^{s_t\lambda_t}+\frac{\alpha_1}{p_1-1}\right)
\end{equation}
and further, by invoking \cref{thm:qualitative-denominator-floor}, it holds $P_0$--\as that 
\begin{equation}\label{eqn:proof-consistency-denominator-liminf-Lambda_t}
    \liminf_{t\to+\infty}\Lambda_t\ge-\mathfrak{h}^\star.
\end{equation}
Hence, for all $\epsilon>0$ eventually $P_0$--\as
\begin{equation*}
\Lambda_t\ge -\mathfrak{h}^\star-\epsilon
\end{equation*}
and therefore, using $\Lambda_t=s_t^{-1}\log Z_t(\Theta)$, 
\begin{equation}\label{eqn:proof-consistency-Z_t-Theta-bound}
    Z_t(\Theta)\ge e^{-s_t(\mathfrak{h}^\star+\epsilon)}
\end{equation}
\as for all $t$ sufficiently large. Combining estimate \labelcref{eqn:proof-consistency-numerator-on-sieve-bound} with inequality \labelcref{eqn:proof-consistency-Z_t-Theta-bound} we obtain eventually \as the following upper bound on equality \labelcref{eqn:proof-consistency-Pi_t-A_t-G_t-equality}
\begin{align}
\Pi_t(A\cap G_t|X_1^t) \le& \frac{\exp(-s_t(\mathfrak{h}(A)\wedge M_0))}{\exp(-s_t(\mathfrak{h}^\star+\epsilon))}\left(e^{s_t\lambda_t}+\frac{\alpha_1}{p_1-1}\right)\notag\\
=&\exp(s_t(\mathfrak{h}^\star+\epsilon-\mathfrak{h}(A)\wedge M_0))\left(e^{s_t\lambda_t}+\frac{\alpha_1}{p_1-1}\right),\notag\\
\intertext{which we may split as}
=& \exp(s_t(\mathfrak{h}^\star+\lambda_t+\epsilon-\mathfrak{h}(A)\wedge M_0))\label{eqn:proof-consistency-posterior-on-sieve-with-lambda_t}\\
&+\frac{\alpha_1}{p_1-1}\exp(s_t(\mathfrak{h}^\star+\epsilon-\mathfrak{h}(A)\wedge M_0)).\label{eqn:proof-consistency-posterior-on-sieve-constant-multiple}
\end{align}
Now, as $\mathfrak{h}(A)\wedge M_0>\mathfrak{h}^\star$, we may choose $\epsilon$ sufficiently small such that \labelcref{eqn:proof-consistency-posterior-on-sieve-constant-multiple} converges to 0 \as noting that $\alpha_1/(p_1-1)$ is constant. For part \labelcref{eqn:proof-consistency-posterior-on-sieve-with-lambda_t}, we notice that $\lambda_t\to 0$ (\as) and hence eventually for all $\eta>0$ we have $\lambda_t\le\eta$, from which we find that \as for all $t$ sufficiently large
\begin{equation}\label{eqn:proof-consistency-final-estimate}
\exp(s_t(\mathfrak{h}^\star+\lambda_t+\epsilon-\mathfrak{h}(A)\wedge M_0))\le \exp(s_t(\mathfrak{h}^\star+\eta+\epsilon-\mathfrak{h}(A)\wedge M_0)).
\end{equation}
Thus, we may choose $\eta>0$ sufficiently small, depending on the sufficiently small $\epsilon>0$ chosen for bound \labelcref{eqn:proof-consistency-posterior-on-sieve-constant-multiple}, such that $\mathfrak{h}(A)\wedge M_0-\mathfrak{h}^\star-\epsilon>\eta$, so that estimate \labelcref{eqn:proof-consistency-final-estimate} converges $P_0$--\as to 0. 
The proof of \cref{thm:consistency} is hence complete.
\end{bareproof}

\subsection{Proof of \texorpdfstring{\cref{thm:weak-star-convergence-unique-min}}{Proposition \ref*{thm:weak-star-convergence-unique-min}} completed}
\begin{bareproof}
By the hypothesis that $\Theta$ is first countable at $\nu^\star$, we have that there exists a local neighbourhood basis around $\nu^\star$, that is $\{U_1,U_2,\dots\}$ with $\nu^\star\in U_m$ for all $m\in\N$. We may assume \wolog that the basis is non--increasing, as otherwise we may define the non--increasing sequence of sets
\begin{gather*}
    \tilde{U}_m:=\bigcap_{k=1}^mU_k,\\
    \intertext{satisfying}
    \mathfrak{h}(\Theta\setminus\tilde U_m)=\min_{k\in[1;m]}\mathfrak{h}(\Theta\setminus U_k)>\mathfrak{h}^\star,
\end{gather*}
which forms a basis at $\nu^\star$.
Now, since $\forall m\in\N\;\mathfrak{h}(\Theta\setminus U_m)>\mathfrak{h}^\star$ and applying \cref{thm:consistency} yields for all $m\in\N$ $P_0$--\as
$$
\Pi_t(\Theta\setminus U_m|X_1^t)\to 0.
$$
Letting $N_m\subset \Omega$ be the $P_0$ null sets where the above does not hold for $U_m$, then 
$$
P_0\left(N_\infty:=\bigcup_{m=1}^\infty N_m\right)=0.
$$
Likewise, we may discard the $P_0$ null sets on which $Z_t(\Theta)\in\{0,\infty\}$ for each $t$. Let $C\subseteq\Omega$ be the full measure remainder on which the convergence holds for all $m\in\N$ and $Z_t(\Theta)\in(0,\infty)$ for all $t$. On this set, let $f:\Theta\to\R$ be a bounded, $\mathcal{T}$--measurable function continuous at $\nu^\star$, and let $\delta>0$. There hence exists a $V\subseteq\Theta$ on which $|f-f(\nu^\star)|\le \delta$, and since $(U_m)_m$ is a basis, $\exists m\in\N$ for which $U_m\subseteq V$. It follows, for all $t$ sufficiently large
\begin{align*}
\left|\int_\Theta f\dd\Pi_t-f(\nu^\star)\right|&\le \int_\Theta|f-f(\nu^\star)|\dd\Pi_t\\
&=\int_{U_m}|f-f(\nu^\star)|\dd\Pi_t+\int_{\Theta\setminus U_m}|f-f(\nu^\star)|\dd\Pi_t\\
&\le \delta \Pi_t(U_m|X_1^t)+||f-f(\nu^\star)||_{\infty}\Pi_t(\Theta\setminus U_m|X_1^t)\\
&\le \delta + 2||f||_\infty\Pi_t(\Theta\setminus U_m|X_1^t)
\end{align*}
and 
\begin{equation*}
\delta + 2||f||_\infty\Pi_t(\Theta\setminus U_m|X_1^t)\to\delta
\end{equation*}
since the convergence of $\Pi_t(\Theta\setminus U_m|X_1^t)\to0$ holds on $C$ and $f$ being bounded. Taking $\delta\downarrow0$ yields
\begin{equation*}
\int_\Theta f\dd\Pi_t\to f(\nu^\star).
\end{equation*}
Now, the conditions hold simultaneously for all $f\in C_b(\Theta)$ and as $C$ has full measure, the conclusion follows.
\end{bareproof}

\section{Gibbs posterior contraction rates}\label{sec:rates-proof}

This section completes the proof of our main result, \cref{thm:rates}. The argument proceeds mainly by revisiting the results obtained in \cref{part:Gibbs-measure-properties,part:denominator-analysis}, taking into account now the results obtained in \cref{sec:consistency-proof}.

\subsection{A first Gibbs contraction result}
We turn our attention first to \cref{prop:refined-sieve-estimate}, and provide some more careful analysis under some further assumptions on the behaviour of $\eta_t$, which we will show next to be satisfied.
\begin{theorem}\label{thm:rates-on-sieve-under-eta/epsilon-asp}
Suppose that condition \labelcref{eqn:asp-weak-sieve-deviation-bound} holds for some $\alpha_1$, $p_1$, and $c_1$. Suppose also that $\exists\mathfrak{L}\in\R$ deterministic such that $P_0$--almost surely
\begin{equation}\label{eqn:asp-eta_t-bound}
\limsup_{t\to+\infty}\frac{\eta_t}{\epsilon_t}\le \mathfrak{L}<+\infty.
\end{equation}
Then for all $c'\in(c_1,1)$ there exists a constant $j_0(c')$ such that eventually $P_0$ almost surely, uniformly over $j\in[j_0(c');j^\dagger(t)]$, it holds
\begin{equation}
    \Pi_t(\{h\ge \mathfrak{h}^\star+j\epsilon_t\}\cap G_t|X_1^t)\le \exp(-s_tj\epsilon_t(1-c')).
\end{equation}
The constant $j_0(c')$ is independent of $t$ depending only on $\mathfrak{L},c',c_1$ and equals
\begin{equation}
    j_0(c')=\left\lfloor \frac{\mathfrak{L}\vee0}{c'-c_1}\right\rfloor+ 1.
\end{equation}
\end{theorem}
\begin{proof}
Let $c'\in(c_1,1)$, and $t$ sufficiently large so that $j_0=j_0(c')\in[1;j^\dagger(t)]$ to be specified later. It holds by standing assumption and \cref{lem:finite-Z_t-Theta} that $\exists \Omega_0\subseteq\Omega$ of full measure and $T_0\ge 1$ such that for $\omega\in \Omega_0$ and all $t\ge T_0$, $Z_t(\Theta)(\omega)\in(0,\infty)$. Now, \cref{prop:refined-sieve-estimate} yields a full measure set $\Omega_1\subseteq\Omega$ and a $T_1\ge 1$ for which its conclusion holds for all $\omega\in \Omega_1$ and for all $t\ge T_1$ uniformly for all $j\in[1;j^\dagger(t)]$. Hence, on $\Omega_2:=\Omega_0\cap \Omega_1$, of full measure, and for all $t\ge T_2:=T_0\vee T_1$, we find for all $j\in[j_0;j^\dagger(t)]$ that
\begin{align}
    \Pi_t(\{h\ge\mathfrak{h}^\star+j\epsilon_t\}\cap G_t) &= \frac{Z_t(\{h\ge\mathfrak{h}^\star+j\epsilon_t\}\cap G_t)}{\exp({-s_t(\mathfrak{h}^\star+\eta_t)})}\notag\\
    &\le\exp\left\{-s_t\left(\frac{-\log 2C}{s_t}+\mathfrak{h}^\star+j\epsilon_t(1-c_1)-\mathfrak{h}^\star-\eta_t\right)\right\}\notag\\
    &=\exp\left\{-s_tj\epsilon_t\left(1-c_1-\frac{\eta_t}{\epsilon_t j}-\frac{\log2C}{s_tj\epsilon_t}\right)\right\}\label{eqn:estimate-first-contraction-rate-proof}
\end{align}
where we are precise about the deterministic constant $C$ from the conclusion of \cref{prop:refined-sieve-estimate}.

Now, by assumption \labelcref{eqn:asp-eta_t-bound}, there exists a full measure $\Omega_3\subseteq\Omega$ so that for each $\omega\in \Omega_3$ there exists a finite $T_3(\omega)\ge 1$ such that $\forall t\ge T_3(\omega)$
$$
\frac{\eta_t(\omega)}{\epsilon_t}\le\mathfrak{L}+\frac\delta2j_0.
$$
Moreover, we have a deterministic $T_4\ge 1$ satisfying for all $t\ge T_4$
$$
\frac{\log2C}{s_t\epsilon_t}\le \frac\delta2j_0.
$$
In particular, we may bound estimate \labelcref{eqn:estimate-first-contraction-rate-proof} by
\begin{equation}\label{eqn:first-sieve-rate}
    \le \exp(-s_tj\epsilon_t(1-c'))
\end{equation}
which holds for $(\omega,t)$ if and only if
\begin{equation}
    1-c_1-\frac{\eta_t(\omega)}{\epsilon_t j}-\frac{\log2C}{s_tj\epsilon_t}\ge 1-c'
\end{equation}
or equivalently,
\begin{equation}\label{eqn:c'-estimate}
    \frac{\eta_t(\omega)}{j\epsilon_t}+\frac{\log 2C}{s_tj\epsilon_t}\le c'-c_1.
\end{equation}
Optimising this bound over $j$, we may determine the smallest $j$ for which the above holds. Taking the limit superior of the left hand side of \labelcref{eqn:c'-estimate} pointwise for $\omega\in \Omega_2\cap \Omega_3$ leads us to 
\begin{equation*}
c'-c_1> \frac{\mathfrak{L}}{j}\iff j >\frac{\mathfrak{L}}{c'-c_1}. 
\end{equation*}
Thus for the estimate \labelcref{eqn:first-sieve-rate} to converge to $0$ we require 
\begin{equation*}
    j^\dagger(t)\ge j\ge j_0(c'):=\left\lfloor \frac{\mathfrak{L}\vee0}{c'-c_1}\right\rfloor+ 1,
\end{equation*}
where this lower bound is independent of $t$. Note, there exists a $T_5\ge1$ for which $j_0(c')\le j^\dagger(t)$ for all $t\ge T_5$.

Finally, we may pick $\delta\in(0,c'-c_1-\mathfrak{L}/j_0(c'))$ which is non--empty by definition, so that for $\omega\in \Omega_3$ and all $t\ge T_3(\omega)\vee T_4\vee T_5$
$$
\frac{\eta_t(\omega)}{\epsilon_t}+\frac{\log 2C}{s_t\epsilon_t}\le \mathfrak{L}+\delta j_0(c')<\mathfrak{L}+(c'-c_1)j_0(c')-\mathfrak{L}
$$
and
$$
\frac{\eta_t(\omega)}{j\epsilon_t}+\frac{\log 2C}{js_t\epsilon_t}<(c'-c_1)\frac{j_0(c')}{j}\le c'-c_1
$$
uniformly for all $j\in[j_0(c');j^\dagger(t)]$.

To conclude, estimate \labelcref{eqn:first-sieve-rate} holds for $\omega\in \Omega_2\cap \Omega_3$, for all $t\ge T_2 \vee T_3(\omega)\vee T_4\vee T_5$, and since $\Omega_2$ and $\Omega_3$ have full $P_0$ measure so does their intersection. As $c'\in(c_1,1)$ was arbitrary, the conclusion follows.
\end{proof}

The theorem yields the following estimate for the control of the divergence behaviour of $R_t$ on the set $I:=\{h=+\infty\}$, and moreover, as a consequence of the truncation via $M_0>\mathfrak{h}^\star$ on $\{h\ge M_0\}$. 
\begin{corollary}
For $I=\{h=+\infty\}$ and $\Pi_0(I)<1$, then eventually, $P_0$--\as
\begin{equation}
    \Pi_t(I\cap G_t|X_1^t)\le \exp(-s_t(1-c')(M_0-\mathfrak{h}^\star -\epsilon_t)).
\end{equation}
\end{corollary}
\begin{proof}
Using $j^\dagger(t):=\lfloor(M_0-\mathfrak{h}^\star)\epsilon_t^{-1}\rfloor$, and hence $j^\dagger(t)\epsilon_t>M_0-\mathfrak{h}^\star-\epsilon_t$, and we conclude since $\{h\ge M_0\}\subseteq\{h\ge \mathfrak{h}^\star+j^\dagger(t)\epsilon_t\}$ eventually.
\end{proof}
\subsection{Proof of \texorpdfstring{\cref{thm:rates}}{Theorem \ref*{thm:rates}} completed}
\begin{bareproof}
Let $M>0$ to be determined later and write 
$$N_{M\epsilon_t}:=\{P\in\Theta:h(P)<\mathfrak{h}^\star+M\epsilon_t\}.$$ 
To demonstrate $\Pi_t(N_{M\epsilon_t}|X_1^t)\to 1$ almost surely, it is equivalent to prove that $\Pi_t(N_{M\epsilon_t}^c|X_1^t)\to 0$ \as Then, we may decompose the Gibbs posterior over $N_{M\epsilon_t}^c$ by partitioning it through $G_t$ into disjoint subsets as 
\begin{equation}
\begin{aligned}
\Pi_t(\{h\ge\mathfrak{h}^\star+M\epsilon_t\}|X_1^t) &= \Pi_t(N_{M\epsilon_t}^c\cap G_t|X_1^t)+\Pi_t(N_{M\epsilon_t}^c\cap G_t^c|X_1^t)\\
&\le\Pi_t(N_{M\epsilon_t}^c\cap G_t|X_1^t)+\Pi_t(G_t^c|X_1^t).
\end{aligned}
\end{equation}
Now, by standing assumption on $G_t^c$, this satisfies the hypothesis of \cref{thm:contraction_bad-sets}, and hence 
\begin{equation}
    \Pi_t(G_t^c|X_1^t)\to 0\qquad P_0\mathrm{-}\as
\end{equation}

It therefore remains to be shown that (\as)
\begin{equation*}
    \Pi_t(N_{M\epsilon_t}^c\cap G_t|X_1^t)\to 0.
\end{equation*}
Indeed, by assumption \cref{asp:averaged-g_t} we are in a position to apply \cref{thm:pathwise-denominator-contraction}. It therefore holds $P_0$--almost surely that
\begin{equation*}
    \limsup_{t\to+\infty}\frac{\eta_t}{\epsilon_t}\le 2a+a'+b+\bar{b}=:\mathfrak{L}.
\end{equation*}
Hence, for any $c'\in(c_1,1)$ \cref{thm:rates-on-sieve-under-eta/epsilon-asp} applies and we are led to the estimate
\begin{equation}\label{eqn:rates-proof-eqn-sieve-bound}
    \Pi_t(G_t\cap\{h\ge \mathfrak{h}^\star+j\epsilon_t\}|X_1^t)\le \exp(-s_t\epsilon_t j_0(c')(1-c'))
\end{equation}
which holds eventually $P_0$--almost surely for all $j\in[j_0(c');j^\dagger(t)]$.
Recognising now that for all $t$
\begin{equation}
     N_{M\epsilon_t}^c\subseteq \{h\ge\mathfrak{h}^\star+j_0(c')\epsilon_t\}
\end{equation}
as soon as $M\ge j_0(c')$, where 
\begin{equation*}
    j_0(c'):=\left\lfloor \frac{\mathfrak{L}\vee0}{c'-c_1}\right\rfloor+ 1
\end{equation*}
is independent of $t$. Therefore, as $s_t\epsilon_t\to\infty$ the estimate \labelcref{eqn:rates-proof-eqn-sieve-bound} converges to zero \as for all such $M\ge j_0(c')$. 
The finite sample contraction rate follows by noting $\lfloor M\rfloor\le j^\dagger(t)$ eventually and 
$$\{h\ge \mathfrak{h}^\star+M\epsilon_t\}\subseteq\{h\ge \mathfrak{h}^\star+\lfloor M\rfloor\epsilon_t\}.$$
 We may therefore apply \cref{thm:rates-on-sieve-under-eta/epsilon-asp} with $j=\lfloor M\rfloor\ge j_0(c')$ and the posterior on $G_t^c$ decays faster or at the same rate by \cref{thm:contraction_bad-sets}.
The proof of \cref{thm:rates} is complete.
\end{bareproof}

\subsection{Proof of \texorpdfstring{\cref{thm:rates-endpoint-control}}{Theorem \ref*{thm:rates-endpoint-control}}}

The conditions for \cref{thm:rates-endpoint-control,thm:rates} detail two distinct methods for verifying that condition \labelcref{eqn:asp-eta_t-bound} holds. Thus, reasoning as in the proof of \cref{thm:rates} we may likewise conclude \cref{thm:rates-endpoint-control}. 

\begin{proof}[Proof of \cref{thm:rates-endpoint-control} completed]
It follows by \cref{asp:rates-endpoint-ver} that we may apply \cref{thm:denominator-via-lorentz-bound}. From this, we have almost surely that
$$\limsup_{t\to\infty}{\eta_t}/{\epsilon_t}\le\mathfrak{L}$$ 
with $\mathfrak{L}=(1+p_2)a$. Arguing as in the proof of \cref{thm:rates}, we derive the desired result.
\end{proof}

\subsection{Pre--metric contraction}
Further, we may now establish \cref{prop:pre-metric-contraction}.
\begin{proof}[Proof of \cref{prop:pre-metric-contraction} completed]
By assumption, we have that
$$
h(P)-\mathfrak{h}^\star\gtrsim d(P,\tilde{\Theta}_\star)^w.
$$
Writing $N_{M\epsilon_t}$ as before with $M$ as in \cref{thm:rates}, we have that
\begin{align*}
N_{M\epsilon_t}:=\{h< \mathfrak{h}^\star+M\epsilon_t\}=\{h-\mathfrak{h}^\star< M\epsilon_t\}.
\end{align*}
Thus, using relation \labelcref{eqn:pre-metric-asp}, we find that
$$
M\epsilon_t>h(P)-\mathfrak{h}^\star\gtrsim d(P,\tilde{\Theta}_\star)^{w}
$$
and therefore,
\begin{equation}
    N_{M\epsilon_t}:=\{P:h(P)-\mathfrak{h}^\star< M\epsilon_t\}\subseteq\left\{P: d(P,\tilde{\Theta}_\star)\lesssim(M\epsilon_t)^{1/w}\right\}.
\end{equation}
We thus conclude by application of \cref{thm:rates,thm:rates-endpoint-control}.
\end{proof}

\subsection{Almost sure contraction rates for Bayes posterior measures}\label{sec:Bayes-proof}

We conclude by returning to the relationship of Gibbs and Bayes posteriors.

The Bayes posterior naturally satisfies a supermartingale property that allows %
us to bound $c_t\le 1$. This inequality is in general strict.
\begin{proposition}\label{prop:Bayes-c_t<=1}
For all $t$, the dominated Bayes posterior satisfies 
\begin{equation}
    c_t(P):=\int_{\Xi^{(t)}}R_t(P)\dd P_0^{(t)}\le 1
\end{equation}
and hence satisfies \cref{asp:lorentz} with $p=\infty$ and $\kappa=0$.
\end{proposition}
\begin{proof}
Let $P^{(t)}\in\Theta$ be any probability measure on $(\Xi^{(t)},\mathcal{X}^{(t)})$. We may obtain its Lebesgue decomposition with respect to $P_0^{(t)}$ as
\begin{equation}
    P^{(t)}=P^{(t),s}+ P^{(t),ac}\quad\mathrm{satisfying}\quad P^{(t),s}\perp P^{(t),ac}
\end{equation}
where the finite measure $P^{(t),ac}$ is absolutely continuous with respect to $P_0^{(t)}$, that is 
\begin{equation}
    P^{(t),ac}\ll P_0^{(t)}
\end{equation}
and further that
\begin{equation}
    P^{(t),s}\perp P_0^{(t)}.
\end{equation}

From this we may rigorously define the likelihood ratio $R_t$ as 
\begin{equation}
    \dd P_0^{(t)}R_t =\dd P^{(t),ac},
\end{equation}
where
\begin{equation}
    R_t=\frac{dP^{(t),ac}}{dP_0^{(t)}}.
\end{equation}
Thus,
\begin{equation}
\begin{aligned}
\int_{\Xi^{(t)}}R_t\dd P_0^{(t)}&= \int_{\Xi^{(t)}}\frac{dP^{(t),ac}}{dP_0^{(t)}}(x_1^t)\dd P_0^{(t)}(x_1^t)\\&=\int_{\Xi^{(t)}}P^{(t),ac}(dx_1^t)= P^{(t),ac}(\Xi^{(t)})=1-P^{(t),s}(\Xi^{(t)})\le 1.
\end{aligned}
\end{equation}
Therefore, $c_t\in L^{\infty}(\Theta,\Pi_0)$ and hence $p=+\infty$, and as $$\log c_t\le\log 1=0,$$ $\kappa=0$.
\end{proof}
One may further establish that $R_t$ is a supermartingale for Bayes (see e.g. \cite{walker2004})
but we shall not do this here.

It thus remains to be verified that the remaining assumptions are satisfied.
\begin{proof}[Proof of \cref{thm:Bayes} completed]
The first part of the theorem is established by \cref{prop:Bayes-c_t<=1}.

The consistency result follows by application of \cref{thm:consistency} by the assumptions of the theorem. 
Finally, the contraction rates for Bayes posteriors follow by appealing to \cref{thm:rates,thm:rates-endpoint-control} depending on which assumption has been posited.
\end{proof}

\appendix
\crefalias{section}{appendix}
\section{Additional background material}\label{apx:measure-theory-background}
We provide some non--standard definitions used throughout the paper implicitly, for the reader's convenience. In this we broadly follow \cites{billingsley1999, rudin1987,grafakos2014}.

\begin{definition}[Essential Infimum]\label{def:essential_infimum}
Take a measurable set $E\subseteq X$, and a measurable function $f:X\to[-\infty,+\infty]$, then we define the essential infimum as
\begin{equation}
    \essinf_{E}f:=\inf\{b\in \R:\mu(\{x\in E:f(x)<b\})>0\},
\end{equation}
where we identify the case where $\mu(E)=0$ with $\essinf_Ef=+\infty$.
\end{definition}
It naturally follows that if $\mu(E)>0$ then for every $a>\essinf_Ef$ there exists a set $S_a\subseteq E$ of positive measure with $f(x)< a$ on $S_a$, and hence that all sublevel sets of $f$ with level greater than the essential infimum have positive measure. 

\begin{definition}[Weak$^\star$ convergence \cite{billingsley1999}]\label{def:weak-star_convergence}
We say a sequence of probability measures $\mu_n$ on $(X,\mathcal B)$ weak$^\star$ converges to $\mu$, written $\mu_n\overset{w^\star}{\to}\mu$, if
\begin{equation}
    \int_X f\dd\mu_n\to\int_Xf\dd\mu\qquad\forall f\in C_b(X)
\end{equation}
where $C_b(X)$ is the space of all bounded continuous functions on $X$.
\end{definition}

\begin{definition}[Lorentz spaces $L^{p,\infty}$]\label{def:lorentz-space}
For the $\sigma$--finite measure space $(X,\mu)$ and $p\in[1,\infty)$, we define the Lorentz space $L^{p,\infty}$ as
\begin{equation}
    L^{p,\infty}(X,\mu):=\left\{f\in\mathcal L(X):||f||_{L^{p,\infty}}<\infty\right\}/\sim
\end{equation}
where 
\begin{equation}\label{eqn:lorentz-norm}
||f||_{L^{p,\infty}}:=\sup_{\lambda>0}\lambda \cdot\left[\mu(\{x:|f(x)|>\lambda\})\right]^{1/p}.
\end{equation}
For $p=+\infty$, we identify $L^{\infty,\infty}=L^\infty$.
\end{definition}

Immediately, by Markov's inequality, for $f\in L^p(\mu)$ it holds $f\in L^{p,\infty}(\mu)$ as
$$
\mu(|f|>\lambda)\le \frac{||f||^p_{L^p}}{\lambda^p}\qquad\forall\lambda>0.
$$
If additionally $(X,\mu)$ is a probability space, it is worth recalling the continuous embedding of the Lorentz spaces. In particular, for any $1\le p< q\le+\infty$ we have 
\begin{equation*}
    L^{q}\hookrightarrow L^{q,\infty}\hookrightarrow L^{p}\hookrightarrow L^{p,\infty}\,.
\end{equation*}

We further define the following space.
\begin{definition}[Zygmund space $L\log L$]\label{def:L-log-L-space}
Suppose $(X,\mu)$ is a $\sigma$--finite measure space, then we define 
\begin{equation}
    L\log L\,(X,\mu):=\left\{f\in\mathcal L:\int_X |f|\log^+|f|\dd\mu<\infty\right\}/\sim
\end{equation}
with associated norm $||\cdot||_{L\log L}$.
\end{definition}
In general, unless of course the measure $\mu$ trivialises the situation, we have 
$$L\log L\subsetneq L^1\subsetneq L^{1,\infty}.$$
In fact, we may observe the following result.

\begin{proposition}
Suppose $(X,\mu)$ is a probability space, $R:X\to[0,\infty]$ and $R\in L^1(\mu)$ with $R\not\equiv 0$. Define the posterior $\hat\nu$ via its Radon--Nikodym derivative as
$$\frac{d\hat \nu}{d\mu}(x):=\frac{R(x)}{||R||_{L^1(\mu)}}$$
then
\begin{equation}
    D_\mathrm{KL}(\hat\nu:\mu)<\infty \iff R\in L\log L(X,\mu).
\end{equation}
\end{proposition}
\begin{proof}
It follows, writing $Z=||R||_{L^1(\mu)}$, that
\begin{align*}
D_\mathrm{KL}(\hat\nu:\mu)&=\E^{\mu}\left[\frac{d\hat\nu}{d\mu}(x)\log\left(\frac{d\hat\nu}{d\mu}(x)\right)\right]\\
&=\int_X \frac1ZR(x)\log\left(\frac1ZR(x)\right)\mu(dx)\\
&=\frac1Z\int_XR(x)\left(\log R(x)-\log Z\right)\mu(dx)\\
&=\frac1Z\int_X R(x)\log R(x)\,\mu(dx)-\log Z
\end{align*}
Now, as $Z\in(0,\infty)$ by assumption, $\log Z\in \R$ and $Z^{-1}\in(0,\infty)$, hence 
$$
D_\mathrm{KL}(\hat\nu:\mu)<\infty\iff \int_X R\log R\dd\mu<\infty\iff R\in L\log L(X,\mu),
$$
as desired.
\end{proof}

\section{Counterexamples to results of Shalizi (2009)}\label{apx:on-shalizi}
Before the proceeding discussion on the results of \cite{shalizi2009}, we preface the following. The ideas presented by Shalizi \cite{shalizi2009} are incredibly useful, the relationship between the replicator dynamics and the Bayes posteriors are fascinating, and moreover that paper has largely motivated even attempting the generality we explore in this work. 

\subsubsection*{On the errors in the proofs}
As part of this paper, we have explored the work of Shalizi \cite{shalizi2009} which considers the asymptotics of Bayes posteriors under similar premises to those of our paper, and came to the conclusion that the proofs stated there are incorrect and do not establish the claimed Theorems 2--6 inclusive, and moreover, we suspect that the proposed assumptions are insufficient to establish the claimed results, Theorems 2--5. These are genuinely insufficient for Theorem 6, the contraction rates result, in \cite{shalizi2009}, as we in fact establish a counterexample in \cref{prop:shalizi-counterexample-rates} below.

We are unaware of any previous work pointing these issues out, and at point of submission of this work, there has been no corrigendum published, nor has an arXiv version been updated. We include a detailed discussion of where the results and proofs fail since the paper has been cited extensively and the claimed theoretical contributions have been widely used.

While we have exclusively focused on asymptotic consistency and contraction rates in this work, Theorems 3 and 6 in \cite{shalizi2009} respectively, we remark that Theorems 2, 4, and 5 are likewise stated incorrectly, but we will not explore this in detail, and that while Theorem 1 holds, the proof requires a slight modification in considering the set $I$, which is precisely why we make use of $M_0>\mathfrak{h}^\star$ throughout our results. For Theorem 5 in \cite{shalizi2009}, we are aware of a counterexample and an immediate fix, we only sketch the argument.  
Provided it is assumed that the posterior vanishes on sets $A$ for which $\mathfrak{h}(A)>\mathfrak{h}^\star$, Theorem 5 of \cite{shalizi2009} may be repaired by considering the limit superior of the Ces\`aro average of the divergences instead of the limit superior of the individual ones, but we have not investigated whether this may be further improved. For a counterexample, one may consider the case of a single hypothesis that errs infinitely often, making the limsup positive, but spaced exponentially far apart, making $\mathfrak{h}^\star=0$.

For the reader's convenience we include here an account of the assumptions made in \cite{shalizi2009}, adopting the notation used throughout this paper, before turning our attention to the results.
\begin{definition} We refer to the following as the assumptions of \cite{shalizi2009}.
\begin{enumerate}
\item $R_t$ is $\sigma(X_1^t)\otimes\mathcal T$--measurable for all $t$.
\item For every $P\in\Theta$, 
\begin{equation*}
    h(P)=\lim_{t\to\infty}t^{-1}\E_{P_0}\left[\log\frac{P_0}{P}(X_1^t)\right]
\end{equation*}
exist and is $\mathcal T$--measurable.
\item For each $P\in\Theta$,
\begin{equation*}
    \lim_{t\to\infty}t^{-1}\log R_t(P)=-h(P)
\end{equation*}
$P_0$--almost surely.
\item For $I:=\{P\in\Theta:h(P)=\infty\}$, $\Pi_0(I)<1$.
\item There exists $G_t\to\Theta$ measurable, for which
\begin{enumerate}[label=(\roman*)]
    \item $\Pi_0(G_t^c)\le \alpha e^{-t\gamma}$ for some $\alpha>0,\gamma>2\mathfrak h^\star$,
    \item The convergence $t^{-1}\log R_t(P)\to -h(P)$ is uniform in $P$ over $G_t\setminus I$, 
    \item $\mathfrak h(G_t)\to\mathfrak h^\star$.
\end{enumerate}
\item Let $\tau(A,\delta)$ be a random stopping time defined for $\delta>0$ and $A\in\mathcal T$ as
\begin{equation*}
    \tau(A,\delta):=\sup\left\{s\in\N:\frac1s\log Z_s(A)>\delta+\limsup_{s'\to\infty}\frac1{s'}\log Z_{s'}(A)\right\}
\end{equation*}
whenever the limit superior is finite. The sets $G_t$ from the previous may be chosen so that, for every $\delta>0$, the stopping time $\tau(G_t,\delta)\le t$ eventually $P_0$--almost surely.
\item Moreover, the sets $G_t$ may be chosen so that for any $A\in\mathcal T$ with positive prior measure, $\mathfrak h(A\cap G_t)\to\mathfrak h(A)$.
\end{enumerate}
\end{definition}
The key assumptions to digest here are 5(ii) and 6. In 5(ii), the behaviour of the convergence of $\Delta_t$ on the set $I$, on which the limiting function $h$ is infinite, is entirely uncontrolled. Furthermore we may read ``the convergence [\dots] is uniform in $P$ over $G_t\setminus I$'' in two ways, the first being for each fixed $n\in\N$ the convergence $t^{-1}\log R_t(P)\to-h(P)$ is uniform over this fixed set $G_n\setminus I$, \ie
\begin{equation*}
\forall n\in\N\quad P_0\left(\omega:\lim_{t\to\infty}\sup_{P\in G_n\setminus I}\left|\frac1t\log R_t(\omega, P)+h(P)\right|=0\right)=1,
\end{equation*}
and the second being a sort of diagonalisation over $t$ in the sense
\begin{equation*}
   P_0\left(\omega: \lim_{t\to\infty}\sup_{P\in G_t\setminus I}\left|\frac{1}{t}\log R_t(\omega,P)+h(P)\right|=0\right)=1.
\end{equation*}
Notice here the difference in the domains of the suprema.
The second is a stronger notion, however Appendix B of \cite{shalizi2009} makes it clear that the former is meant; we nevertheless establish the latter in the counterexample to Theorem 6. But even the latter notion leaves the behaviour on $G_t\cap I$ unaccounted for as we will demonstrate later.

Now, we briefly turn to Assumption 6, in which we have a random stopping time $\tau(A,\delta)(\omega)$, and \cite{shalizi2009} assumes for all $\delta>0$ the diagonalisation that for all $t\ge \tau(G_t,\delta)$ $P_0$--almost surely we have
\begin{equation*}
    \frac1t\log Z_t(G_t)-\limsup_{s\to\infty}\frac1s\log Z_s(G_s)\le \delta,
\end{equation*}
or equivalently the above holds \as for all but for finitely many times $t$.
Notice however that this assumption depends crucially on not only $\omega$, but also on the sets $G_t$ themselves inside the limit superior, where we integrate $R_t$ with respect to the prior $\Pi_0$ over this set $G_t$. This however does not allow us to exchange the sequence $G_t$ with the sequence $G_t\cap A$ even when $\Pi_0(A)>0$, as this fundamentally changes the random stopping times $\tau(A\cap G_t,\delta)(\omega)$.

For consistency, perhaps most significant to our purposes is Lemma 10 of \cite{shalizi2009}, where Fubini is incorrectly applied to a set dependent on $\omega$, as the sets produced via Egorov's theorem will in general depend on $\omega$, and moreover, they are assigned a rate of contraction, which is likewise not permitted, unproved, and unassumed. 
Our \cref{thm:qualitative-denominator-floor} may be used as a replacement and holds under weaker assumptions and is the only part of Lemma 10 used in Theorem 3 \cite{shalizi2009}, but it is only a one--sided result which does not in itself replace the incorrect Lemma 10 of \cite{shalizi2009} entirely. This however only yields the desired denominator bound, but the matching numerator result incorrectly applies Assumption 6 and thus fails.
To see this, the proof of Theorem 3 assumes that we may replace $G_t$ with $G_t\cap A$ and retain the diagonalisation structure proposed in Assumption 6 of \cite{shalizi2009}, but Assumption 7 does not provide this. This is due to the change in the random stopping time $\tau(G_t\cap A,\delta)(\omega)$ from $\tau(G_t,\delta)(\omega)$ of which only the latter is controlled by Assumption 6 as the limit superior itself is modified by the intersection, and uniformity does not provide this as no control is assumed in 5(ii) on the part of $A\cap G_t\cap I$, even if a diagonal uniformity condition had been assumed.

In fact, Assumption 5(iii) and 7 of \cite{shalizi2009} are subsumed in Assumption 5(i) and cannot serve as the replacement for the missing diagonalisation structure or serve to provide sufficient prior mass on $N_{\epsilon_t}$. 

\begin{lemma}\label{lem:asp-shalizi7}
Let $G_t\in\mathcal{T}$ satisfy $\Pi_0(G_t^c)\to0$ as $t\to\infty$, then for any $A\in\mathcal{T}$, $\mathfrak{h}(G_t\cap A)\to \mathfrak{h}(A)$, and in particular this holds for $A=\Theta$, that is $\mathfrak{h}(G_t)\to\mathfrak{h}^\star$.
\end{lemma}
\begin{proof}
Trivially, $\mathfrak{h}(G_t\cap A)\ge \mathfrak{h}(A)$ for all $t$ as the (essential) infimum over some set is necessarily greater or equal to any set containing it. Note also, this establishes the result if $\mathfrak{h}(A)=\infty$ and in particular also the case where $\Pi_0(A)=0$, so we may assume $\mathfrak{h}(A)<\infty$ and $\Pi_0(A)>0$ in the remainder. 

Now, let $\epsilon>0$, and define $N_{\epsilon|A}:=\{P\in A:h(P)\le \mathfrak{h}(A)+\epsilon\}$ and by definition of the essential infimum, this set has positive prior mass, $\Pi_0(N_{\epsilon|A}):=m_\epsilon>0$. Now, since $\Pi_0(G_t^c)\to 0$ by assumption, as soon as $\Pi_0(G_t^c)<m_\epsilon$, the intersection $G_t\cap N_{\epsilon|A}$ has positive prior measure. Thereby, $\mathfrak{h}(G_t\cap A)\le \mathfrak{h}(A)+\epsilon$ for all $t$ sufficiently large. As $\epsilon>0$ was arbitrary, the result follows.
\end{proof}

For the contraction rate theorem, \cite{shalizi2009} further posits the following additional hypothesis.
\begin{definition}[Additional assumption of Theorem 6 in \cite{shalizi2009}]
For $\epsilon_t>0$ for which $\epsilon_t\to0$ and $t\epsilon_t\to\infty$, it holds for all $\delta>0$ that
\begin{equation*}
    \tau(G_t\cap N_{\epsilon_t^c},\delta)\le t
\end{equation*}
eventually $P_0$--almost surely.
\end{definition}

Regarding contraction rates, Theorem 6 in \cite{shalizi2009}, fails in general and the proof and assumptions are genuinely insufficient. In particular, \cite{shalizi2009} posits that for all $\eta_3>0$, eventually \as
\begin{equation*}
\Pi_t(N_{\epsilon_t}^c\cap G_t|X_1^t)\le \exp(t(-\epsilon_t+\eta_3))
\end{equation*}
from which \cite{shalizi2009} concludes that this converges to 0. This is false, and in fact, this diverges. Indeed, recalling that $\epsilon_t\to 0$, $\epsilon_t$ is eventually strictly smaller than any $\eta_3>0$ so that the above diverges to $+\infty$, and hence the posterior does not concentrate. Moreover, taking a decreasing sequence of $\eta_3$ (which is permitted) fails as it has no rate attached, since the ``eventually \as'' part for $\eta_3$  directly yields that each different $\eta_3$ has its own random time $T_0=T(\omega,\eta_3)$ attached after which the bound applies, and thus, attaching the rate $\eta_3=o(\epsilon_t)$ fails precisely because of the random time $T_0$. Moreover, Lemma 4 of \cite{shalizi2009} establishes that sufficient prior mass on sets is necessary for posterior contraction on these sets, however the contraction rate theorem does not account for this.

\subsubsection*{A simple counterexample}
We may consider a simple, parametric counterexample to Theorem 6 of \cite{shalizi2009} where the small ball control fails but all assumptions are satisfied. From this, we may establish that Theorem 6 of \cite{shalizi2009} fails. 
\begin{proposition}\label{prop:shalizi-counterexample-rates}
Suppose $X_i\overset{\iid}{\sim}\mathrm{Bernoulli}(\theta_0)$ with $\theta_0=1/2$. Let $$\Theta=\{\mathrm{Bernoulli(\theta)}:\theta\in[1/4,3/4]\}$$ which is parametrised by $\tilde{\Theta}=[1/4,3/4]$ and suppose that for $m>2$ the prior density with respect to $(\tilde{\Theta},\mathrm{Leb})$ satisfies
$$
\pi_0(\theta)\propto\exp(-|\theta-\theta_0|^{-m}).
$$
which is well defined. 
Then, taking $\epsilon_t=t^{-2/m}$, the assumptions of Theorem 6 in \cite{shalizi2009} are satisfied but 
$$
\Pi_t^{(B)}(\{\theta:h(\theta)\le \mathfrak{h}^\star+\epsilon_t\}|X_1^t)\to 0
$$
$P_0$--almost surely.
\end{proposition}
\begin{proof}
We first show that Assumptions 1--7 of \cite{shalizi2009}, as well as the additional assumption imposed in the statement of Theorem 6 in \cite{shalizi2009}, hold. We restate these when we prove each.

Since all hypotheses are probability densities on $\tilde{\Theta}$ with respect to the Lebesgue measure, and writing $\overline{X_1^t}=t^{-1}\sum_{i=1}^tX_i$,
$$R_t(\theta)=\frac{p_\theta(X_1^t)}{p_{\theta_0}(X_1^t)}=\frac{\theta^{t\overline{X_1^t}}(1-\theta)^{t-t\overline{X_1^t}}}{\theta_0^{t\overline{X_1^t}}(1-\theta_0)^{t-t\overline{X_1^t}}}$$
is well defined on $\tilde{\Theta}$ and measurable with respect to the Lebesgue measure for all $t$. The numerator and denominator are strictly positive for all $t$ since $\theta,\theta_0\in\tilde{\Theta}$ bounded away from $0$ and $1$. Assumption 1, $R_t$ is $\sigma(X_1^t)\otimes\mathcal T$ measurable for all $t$, is thus satisfied.

Moreover, Assumption 2 and 4 are trivially satisfied as
\begin{equation*}
h(\theta)=\lim_{t\to\infty}t^{-1}\E_{P_{\theta_0}^{\otimes t}}[-\log R_t(\theta)]=\lim_{t\to\infty}\E_{P_{\theta_0}}\left[\log\frac{p_\theta(X)}{p_{\theta_0}(X)}\right]=D_\mathrm{KL}(P_{\theta_0}||P_\theta)
\end{equation*}
from which it is straightforward that
$$
h(\theta)=\theta_0\log\frac{\theta_0}\theta+(1-\theta_0)\log\frac{1-\theta_0}{1-\theta}.
$$
Hence, nonnegative $h(\theta)$ exists for all $\theta\in\tilde\Theta$ and is $\mathcal T$ measurable (Assumption 2),
and this is strictly finite on $\tilde{\Theta}$ so $I=\{h=+\infty\}=\emptyset$ (Assumption 4, $\Pi_0(I)<1$). 

For Assumptions 3 and 5 of \cite{shalizi2009}, we may take the sieve $G_t\equiv \tilde{\Theta}$, on which the convergence is uniform. To see this, we note that, using $\theta_0=1/2$
\begin{align*}
||t^{-1}\log R_t(\theta)+h(\theta)||_u&=\left|\left|t^{-1}\log \frac{\theta^{t\overline{X_1^t}}(1-\theta)^{t-t\overline{X_1^t}}}{\theta_0^{t\overline{X_1^t}}(1-\theta_0)^{t-t\overline{X_1^t}}} +h(\theta)\right|\right|_u\\
&=\left|\left|\left(\theta_0-\overline{X_1^t}\right)\log\frac{\theta_0}\theta -\left(\theta_0-\overline{X_1^t}\right)\log\frac{1-\theta_0}{1-\theta}\right|\right|_u\\
&=\left|\left|\left(\theta_0-\overline{X_1^t}\right)\log\frac{\theta_0(1-\theta)}{\theta(1-\theta_0)}\right|\right|_u\\
&=\sup_{\theta\in\tilde{\Theta}}\left|\log\frac{1-\theta}{\theta}\right|\left|\theta_0-\overline{X_1^t}\right|=\log3\left|\theta_0-\overline{X_1^t}\right|.
\end{align*}
We use $||\cdot||_u$ to distinguish the uniform norm from the $L^\infty$ norm.
The strong law of large numbers, hence yields a full measure set $\Omega'\subseteq\Omega$ on which the right hand side converges to 0. Hence, we have established both Assumption 3 and 5, since $t^{-1}\log R_t\to-h$ not only pointwise \as (Assumption 3 of \cite{shalizi2009}) but also uniformly on the sieve $G_t\equiv \tilde\Theta$ which satisfies the mass conditions on this trivially (Assumption 5 of \cite{shalizi2009}). 

Assumption 7 of \cite{shalizi2009}, $\mathfrak{h}(A\cap G_t)\to \mathfrak{h}(A)$, is immediate by the trivial choice of sieve $G_t\equiv \tilde \Theta$, but also implied by Assumption 5 per application of \cref{lem:asp-shalizi7}.

For Assumption 6 and the additional assumptions of Theorem 6, we note first that $\epsilon_t\to 0$ and $t\epsilon_t=t^{1-2/m}\to\infty$ since $m>2$. 
We establish this through the following lemma.
\begin{lemma}\label{lem:shalizi-counterexample}
Suppose there exists an $\Omega'\subseteq\Omega$ of full $P_0$ measure, such that $\forall\omega\in\Omega'$
\begin{equation*}
    \sup_{P\in\Theta}\left|\frac1t\log R_t(P,\omega)+h(P)\right|\to0.
\end{equation*}
Then, for all $\delta>0$ and $\omega\in\Omega'$ there exists a finite $T(\delta,\omega)\ge 1$ for which it holds for any $A\in\mathcal{T}$ with $\Pi_0(A)>0$ that
\begin{equation*}
\tau(A,\delta)(\omega)< t \qquad\forall t\ge T(\delta,\omega),
\end{equation*}
where
$$
\tau(A,\delta):=\sup\left\{s\in\N:\frac1s\log Z_s(A)>\delta +\limsup_{s'\to\infty}\frac1{s'}\log Z_{s'}(A)\right\}.
$$
And in particular this holds for any positive $\epsilon_t\to0$ and $A_t:=\{P:h(P)>\mathfrak{h}^\star+\epsilon_t\}$, that is
$$
\tau(A_t,\delta)(\omega)<t,\qquad \forall t\ge T(\delta,\omega).
$$
\end{lemma}
\begin{proof}[Proof of \cref{lem:shalizi-counterexample}]
Let $A\in\mathcal{T}$ with $\Pi_0(A)>0$. The conclusion is immediate when $\mathfrak{h}(A)=\infty$, so suppose this is finite. 
Now, by assumption, for all $\delta>0$ and $\omega\in\Omega'$ there exists a finite $T(\delta,\omega)\ge1$ such that $\forall t\ge T(\delta,\omega)$ 
$$
\sup_{P\in\Theta}\left|\frac1t\log R_t(P,\omega)+h(P)\right|\le \delta
$$
and hence,
$$
Z_t(A,\omega)\le \int_Ae^{-t(h-\delta)}\dd\Pi_0\le e^{-t(\mathfrak{h}(A)-\delta)}\Pi_0(A).$$
That is
$$
\limsup_{s\to+\infty}\frac1s\log Z_s(A,\omega)\le -\mathfrak{h}(A)+\delta<\infty,
$$
since $s^{-1}\log \Pi_0(A)\to 0$.

Further, we have for any $\eta,\delta'>0$ on $A_\eta:=A\cap\{P:h(P)\le \mathfrak{h}(A)+\eta\}$, which has positive measure by definition of the essential infimum, that for all $t\ge T(\delta',\omega)$
$$
Z_t(A,\omega)\ge Z_t(A_\eta,\omega)\ge e^{-t(\mathfrak{h}(A)+\eta+\delta')}\Pi_0(A_\eta)
$$
from which it follows
$$
\liminf_{s\to+\infty}\frac1s\log Z_s(A,\omega)\ge -\mathfrak{h}(A)
$$
by taking countable sequences of $\eta_m,\delta'_m\to0$. The $\delta'>0$ here is allowed to vary as we only use the limit inferior part of the statement.

Combining these estimates, for all $\omega\in\Omega'$
$$
-\mathfrak{h}(A)\le \liminf_{s\to\infty}\frac1s\log Z_s(A,\omega)\le\limsup_{t\to\infty}\frac1s\log Z_s(A,\omega)\le -\mathfrak{h}(A)+\delta.
$$
Therefore, $\forall t\ge T(\delta,\omega)$, using $\log \Pi_0(A')\le 0$ for all $A'\in\mathcal{T}$,
$$
\frac1t\log Z_t(A,\omega)\le -\mathfrak{h}(A)+\delta\le\delta+\limsup_{s\to\infty}\frac1s\log Z_s(A,\omega)<\infty
$$
and hence 
$$
\tau(A,\delta)(\omega)\le T(\delta,\omega)-1<t,\qquad\forall t\ge T(\delta,\omega).
$$
Since $A$ is arbitrary, the conclusion follows. The second part holds as $\Pi_0(A_t)>0$ eventually and hence the previous applies with $\mathfrak{h}(A_t)\ge \epsilon_t$.
\end{proof}
\cref{lem:shalizi-counterexample} establishes both Assumption 6, \ie for all $\delta>0$ $\tau(G_t,\delta)\le t$  \as for all sufficiently large $t$, as well as the additional assumption of Theorem 6, $\forall\delta>0$ we have $\tau(N_{\epsilon_t}^c\cap G_t,\delta)\le t$ \as for all $t$ sufficiently large.  

Having established that the assumptions of Theorem 6 of \cite{shalizi2009} hold, we now show the conclusion fails. In particular, we have (where we use the definition of \cite{shalizi2009} for $N_{\epsilon_t}$ instead of the one of our results; they are equivalent in this setting with ours yielding a slightly stronger result in general)
\begin{align*}
N_{\epsilon_t}:=\{\theta:h(\theta)\le \mathfrak{h}^\star+\epsilon_t\}\subseteq\{\theta:2(\theta-\theta_0)^2\le \epsilon_t\}=:\tilde{N}_{\epsilon_t}
\end{align*}
which holds by application of Pinsker's inequality for Bernoulli random variables, that is $D_\mathrm{KL}(\theta_0||\theta)\ge 2(\theta-\theta_0)^2$. Since the model is well specified, $\mathfrak{h}^\star=0$ attained uniquely at $\theta=1/2$. Hence, on $\tilde{N}_{\epsilon_t}$
\begin{align*}
|\theta-\theta_0|^2\le \epsilon_t/ 2 &\iff \left(|\theta-\theta_0|^2\right)^{{ m}/2}\le \left(\epsilon_t/2\right)^{ m/2}\\
    &\iff \exp(-|\theta-\theta_0|^{-m})\le \exp\left(-\epsilon_t^{-m/2} 2^{m/2}\right)
\end{align*}
having used that $x\mapsto e^{-1/|x|^a}$ for $a\ge1$ is increasing in $|x|$ and continuous at $0$. Now, it is immediate that
\begin{align*}
\Pi_0(N_{\epsilon_t}) &\le \Pi_0(\tilde{N}_{\epsilon_t})\\
&=\int_{\theta\in \tilde{N}_{\epsilon_t}}\pi_0(\theta)\dd\theta\\
&=\int_{1/4}^{3/4}\mathds{1}_{\exp(-|\theta-\theta_0|^{-m})\le \exp\left(-\epsilon_t^{- m/2}  2^{ m/2}\right)}(\theta)\,\pi_0(\theta)\dd\theta\\
&=\int_{1/4}^{3/4}\mathds{1}_{\exp(-|\theta-\theta_0|^{-m})\le \exp\left(-\epsilon_t^{-m/2}  2^{m/2}\right)}(\theta)\,\frac{\exp(-|\theta-\theta_0|^{-m})}{Z_0}\dd\theta\\
&\le \frac{1}{Z_0}\exp\left(-\epsilon_t^{-m/2} 2^{m/2}\right)\\&\lesssim e^{-\gamma t}
\end{align*}
since $\epsilon_t^{-m/2}=(t^{-2/m})^{-m/2}=t$, and where $Z_0\in(0,\infty)$ is the normalising constant of the prior which is finite and positive. Hence, setting $\gamma=2^{m/2}>\mathfrak{h^\star}=0$ and $\alpha=Z_0^{-1}>0$, Lemma 5 of \cite{shalizi2009} applies and $P_{\theta_0}$--almost surely
$$\Pi_t^{(B)}(N_{\epsilon_t}|X_1^t)\to0.$$
To see this more concretely, we may apply our \cref{thm:contraction_bad-sets}, by appealing to \cref{thm:Bayes} and \cref{lem:uniform-convergence-on-sieve,lem:pointwise-a-s-convergence-to-h}, whose assumptions we have just established to hold, which leads us to the \as conclusion.
\end{proof}

\begin{remark}
It is also important to note that \cref{lem:shalizi-counterexample} does not rescue the failure of Theorem 6 of \cite{shalizi2009}, as the stopping time is random. Moreover, $\tau(N_{\epsilon_t}^c\cap G_t, \delta)$ says nothing about how fast $\delta$ may decay to $0$, which is precisely the $\eta_3$ failure discussed earlier.
\end{remark}

\end{document}